\documentclass[reqno, 11pt]{amsart}
\usepackage{amssymb,amsmath,amsfonts}
\usepackage[margin=1in]{geometry}
\usepackage{dsfont}
\usepackage{color}
\usepackage{graphicx,graphics}
\usepackage{latexsym}
\usepackage[dvips]{epsfig}
\usepackage{mathtools,mathrsfs}
\usepackage{stmaryrd,cite}
\usepackage{cases}
\usepackage{enumerate}
\usepackage[usenames,dvipsnames,svgnames]{xcolor}
\usepackage{todonotes}
\usepackage{marginnote}
\usepackage{listings}
\usepackage{cleveref}
\allowdisplaybreaks
\definecolor{refkey}{rgb}{1,0,0.5}
\definecolor{labelkey}{rgb}{0,0.4,1}
\numberwithin{equation}{section}
\newtheorem{thm}{Theorem}[section]
\newtheorem{coro}[thm]{Corollary}
\newtheorem{lem}[thm]{Lemma}
\newtheorem{prop}[thm]{Proposition}
\newtheorem{rmk}[thm]{Remark}

\newcommand{\ea}{\epsilon}
\newcommand{\vea}{\varepsilon}

\newcommand{\al}{\alpha}

\newcommand{\da}{\delta}

\newcommand{\na}{\nabla}

\newcommand{\ta}{\theta}

\newcommand{\sa}{\sigma}

\newcommand{\ga}{\gamma}
\newcommand{\ba}{\beta}

\newcommand{\za}{\zeta}

\newcommand{\oa}{\omega}

\newcommand{\iy}{\infty}
\newcommand{\pl}{\partial}

\newcommand{\lt}{\left}
\newcommand{\rt}{\right}
\newcommand{\be}{\begin{equation}}
\newcommand{\ee}{\end{equation}}
\newcommand{\bee}{\begin{equation*}}
\newcommand{\eee}{\end{equation*}}

\newcommand{\ef}{\eqref}
\newcommand{\f}{\frac}
\newcommand{\les}{\lesssim}

\begin{document}
\title[] {Almost Global Solutions to  the Three-dimensional Isentropic Inviscid Flows with Damping in Physical Vacuum Around Barenlatt Solutions }
\author{Huihui Zeng}
\maketitle

\begin{abstract} For the three-dimensional vacuum free boundary problem with physical singularity that the sound speed is  $C^{ {1}/{2}}$-H$\ddot{\rm o}$lder continuous across the vacuum boundary of the compressible  Euler equations with damping, without any symmetry assumptions, we prove the almost global existence of smooth solutions when the
initial data are small perturbations of the Barenblatt self-similar solutions to the corresponding porous media equations simplified via Darcy's law.
It is proved that if the initial perturbation is of the size of $\epsilon$, then the existing time for  smooth solutions is at least of the order of $\exp(\epsilon^{-2/3})$.  The key issue for the analysis is the slow {\em sub-linear} growth of vacuum boundaries of the order of $t^{1/(3\gamma-1)}$,  where $\gamma>1$ is the adiabatic exponent for the gas. This is in sharp contrast to the currently available global-in-time existence theory of expanding solutions to the vacuum free boundary problems with physical singularity of compressible Euler equations for which the expanding rate of vacuum boundaries is linear. The results obtained in this paper is closely related to the open question in multiple dimensions since T.-P. Liu's construction of particular solutions in 1996 .
\end{abstract}

	


\section{Introduction}
Consider the  following three-dimensional  vacuum free boundary problem for compressible  Euler equations with damping:
\begin{subequations}\label{2.1} \begin{align}
& \pl_t \rho   + {\rm div}(\rho   u ) = 0 &  {\rm in}& \ \ \Omega(t), \label{2.1a}\\
 &  \pl_t  (\rho   u )   + {\rm div}(\rho   u \otimes   u )+\nabla_{  x} p(\rho) = -\rho {  u}  & {\rm in}& \ \ \Omega(t),\label{2.1b}\\
 &\rho>0 &{\rm in }  & \ \ \Omega(t),\label{2.1c}\\
 & \rho=0    &    {\rm on}& \  \ \Gamma(t)=\pl \Omega(t), \label{2.1d}\\
 &    \mathcal{V}(\Gamma(t))={  u}\cdot \mathcal{N}, & &\label{2.1e}\\
&(\rho,{ u})=(\rho_0, { u}_0) & {\rm on} & \ \   \Omega(0), \label{2.1f}
 \end{align} \end{subequations}
where $( t,x)\in [0,\iy)\times\mathbb{R}^3 $,  $\rho $, ${ u} $, and $p$ denote, respectively, the time and space variable, density, velocity and  pressure; $\Omega(t)\subset \mathbb{R}^3$, $\Gamma(t)$, $\mathcal{V}(\Gamma(t))$ and $ \mathcal{N}$ represent, respectively, the changing volume occupied by the gas at time $t$, moving vacuum boundary, normal velocity of $\Gamma(t)$, and exterior unit normal vector to $\Gamma(t)$.  We are concerned with the polytropic gas for which the equation of state is  given by
$$ p(\rho)=\rho^{\gamma}, \  \ {\rm where} \ \  \gamma>1 {~\rm is~the~adiabatic~exponent}. $$
Let $c(\rho)=\sqrt{ p'(\rho)}$ be the sound speed, the condition
\be\label{physical vacuum} -\infty<\nabla_\mathcal{N}\lt(c^2(\rho)\rt)<0  \  \ {\rm on} \ \  \Gamma(t) \ee
defines {\em a  physical vacuum boundary}  (cf. \cite{7,10',16',23,24,25}), which is also called {\em a vacuum boundary with physical singularity} in contrast to the case that $ \nabla_\mathcal{N}\lt(c^2(\rho)\rt)=0$ on  $\Gamma(t)$.
The physical vacuum singularity plays the role of pushing vacuum boundaries, which can be seen by restricting the momentum equation \ef{2.1b} on $\Gamma(t)$: $D_t u \cdot \mathcal{N} = - (\ga-1)^{-1} \nabla_\mathcal{N}\lt(c^2(\rho)\rt) - u\cdot \mathcal{N}$, where $D_t u= (\pl_t + u\cdot \nabla_x)u$ is the acceleration of $\Gamma(t)$, and the term  $- (\ga-1)^{-1} \nabla_\mathcal{N}\lt(c^2(\rho)\rt)>0$ serves as a force due to the pressure effect to accelerate vacuum boundaries.
In order to capture this physical singularity, the initial density is supposed to satisfy
\be\label{initial density}\begin{split}
& \rho_0>0 \ \ {\rm in} \ \ \Omega(0), \ \  \rho_0=0 \ \ {\rm on} \ \ \Gamma(0),   \ \ \int_{\Omega(0)} \rho_0(x) dx =M,
\\
&  -\infty<\nabla_\mathcal{N}\lt(c^2(\rho_0)\rt)<0  \  \ {\rm on} \ \  \Gamma(0),
\end{split}\ee
where $M\in (0, \iy)$ is the initial total mass.

The compressible Euler equations of isentropic flows  with damping, \ef{2.1a}-\ef{2.1b}, is closely related to the
porous media equation (cf. \cite{HL, HMP, HPW, 23, LZ, HZ}):
\begin{equation}\label{pm}
\pl_t \rho  =\Delta p(\rho),
\end{equation}
when \ef{2.1b} is simplified to Darcy's law:
\begin{equation}\label{darcy}
\nabla_{ x} p(\rho)=- \rho { u}.
\end{equation}
(The equivalence can  be seen formally by the rescaling ${ x}'=\ea { x}, t'=\ea^2 t, { u}'= { u}/ \ea$.)
For \ef{pm}, basic understanding of the solution  with finite mass is provided by Barenblatt (cf. \cite{ba}), which is given by
\be\label{1.6}
\bar\rho(t,{ x}) =(1+t)^{-{3}/({3\ga-1})}\lt(\underline{ A }- \underline{B} ({1+ t})^{-{2}/({3\ga-1})}  |{ x}|^2 \rt)^{{1}/({\ga-1})} ,
   \ee
where $\underline{A}$ and $\underline{B}$ are positive constants determined by $\ga$ and the total mass $M$. Precisely,
$$
\underline{B}= \f{\ga-1}{2\ga(3\ga-1)} \ \ {\rm and} \ \  (\ga \underline{A})^{\frac{3\ga-1}{2(\ga-1)}}= \frac{1}{4\pi}M\ga^{\frac{1}{\ga-1}} (\ga \underline{B})^{\frac{3}{2}}
   \lt( \int_0^1 y^2 \lt(1-y^2\rt)^{\frac{1}{\ga-1}}dy \rt)^{-1}.
$$
The Barenblatt self-similar solution defined in  $\bar\Omega(t)=B_{\bar R(t)}(0)$, which  is  the ball  centered at the origin with the radius $\bar R(t) =
\sqrt{\underline{A}/\underline{B}}(1+t)^{ {1}/({3\ga-1}) }$, satisfies
$$\bar\rho> 0 \ \ {\rm in} \ \  \bar\Omega(t), \ \  \bar\rho=0 \ \ {\rm on} \ \ \pl \bar\Omega(t), \ \ {\rm and} \ \   \int_{\bar\Omega(t)}  \bar\rho(t,{ x} )d{ x}  =M  \ \ {\rm for} \ \
t\ge 0. $$
The corresponding  Barenblatt velocity $\bar { u}$ is defined  by
\begin{align*}
\bar { u}=- \f{ \nabla_{ x} p(\bar \rho) }{\bar\rho}=    \f{  x }{(3\ga-1)(1+t)}
 \ \ {\rm in} \ \  \bar\Omega(t).
\end{align*}
So, $(\bar\rho, \bar { u})$ defined in the region $\bar\Omega(t)$
 solves \ef{pm}-\ef{darcy}. There is only one
parameter, total mass $M$, when $\gamma$ is fixed, for the Barenblatt self-similar solution.
We assume that the initial total mass of problem \ef{2.1}  is the same as that for the Barenblatt solution.

It is apparent  that the vacuum boundary $\pl\bar\Omega(t) $  of the Barenblatt solution satisfies the physical vacuum condition, which is the major motivation to study problem \ef{2.1} with the initial condition \eqref{initial density}.  To this end,   a class of particular solutions to problem \ef{2.1} was constructed in  \cite{23} by T. P. Liu  using the following ansatz:
\be\label{liuexplicitsolution}
\Omega(t) =B_{R(t)}({ 0}) ,  \ \
  c^2({x}, t)={e(t)-b(t)r^2},  \  \ { u}({ x}, t)= ({{ x}}/{r} )u(r, t),
  \ee
where $r=| x|$,
$ R(t)= \sqrt {e(t)/b(t)} $    and  $  u(r, t)=a(t) r $.
In \cite{23}, a  system of ordinary differential equations for $(e,b,a)(t) $ was derived  with $e(t), b(t)>0$ for $t\ge 0$, and it was shown that this family of particular solutions  is  time-asymptotically equivalent to the Barenblatt self-similar solution  with the same total mass.
Indeed, the Barenblatt solution  of \ef{pm}-\ef{darcy} can be obtained by the same ansatz as \ef{liuexplicitsolution}:
$
 \bar c^2({ x}, t)=\bar e(t)-\bar b(t)r^2$ and  $  { u}({ x}, t)= \bar a(t) {{ x}}  $, and
it was proved in \cite{23}  that
  $$
  (a,\ b,\ e)(t)=(\bar a, \ \bar b, \ \bar e)(t)+ O(1)(1+t)^{-1}{\ln (1+t)} \ \ {\rm as}\  \ t\to\infty.$$
Since the construction of particular solutions to \ef{2.1} in \cite{23}, it has been an important open question whether there is still a long time existence theory for  problem \ef{2.1} capturing the physical vacuum singular behavior \ef{physical vacuum}, and if there is a
time-asymptotic equivalence of the solution to \ef{2.1} and the corresponding Barenblatt self-similar solution with the same total mass.
This question is answered, respectively,  in the one-dimensional case (cf. \cite{LZ}) and three-dimensional spherically symmetric case (cf. \cite{HZ}).  However, the problem for general  three-dimensional perturbations  without  symmetry assumptions keeps open. The aim of this paper is to investigate this problem.

It is quite challenging to extend the spherically symmetric results in  \cite{HZ} to the general three-dimensional motions.
Because one will have to deal with the intricate evolution of the vacuum boundary geometry and its thorny coupling with the interior solution, and to investigate the bounds for both vorticity and divergence of the velocity field.
Indeed, for the three-dimensional vacuum free boundary problem of compressible Euler equations with physical singularity, the general theory is mostly in the local-in-time nature (cf. \cite{7,10',16'}), and the currently available global-in-time results (cf. \cite{HaJa1,ShSi}) are for expanding solutions of which the expanding rate of vacuum boundaries is linear, $O(1+t)$, when the initial data are small perturbations of affine motions (cf. \cite{sideris1,sideris2}).
(See also \cite{CHJ,HaJa2,PHJ} for related results on global-in-time expanding solutions with linear expanding rate.)
The distinction of problem \ef{2.1} is that the expanding rate of vacuum boundaries for the corresponding Barenblatt solutions, which are the background approximate solutions  for \ef{2.1} in long time, is sub-linear, $O((1+t)^{1/(3\gamma-1)})$, which is less than  $O((1+t)^{{1}/{2}})$ for $\gamma>1$.
This slow expanding rate of  vacuum boundaries creates much severe difficulties in obtaining the long time existence of solutions to problem \ef{2.1}, due to the slow decay of various quantities.
Indeed, the stabilizing effect of fluid expansions also plays important role in the analysis in other context, for example, in general relativistic cosmological models (cf. \cite{Rodnianski, Oliynyk}).
In this article, we prove the almost global existence of solutions to  problem \ef{2.1} in the sense that the lower bound of the life span of solutions is at least $O(\exp\{\epsilon^{-2/3}\})$ if the size of the initial perturbation of the Barenblatt solution is $O(\epsilon)$.
The results obtained in the present work are the first ones for the long time dynamics of vacuum free boundary problems of compressible fluids with physical singularity at the {\em sub-linear} expanding rate of  vacuum boundaries in multi-dimensions.

We review  some previous related works before closing the introduction.
Theoretical study of vacuum states of gas dynamics dates back to 1980 when it was shown in \cite{LiuSmoller} that shock waves vanish at the vacuum.  Early study of well-posedness of smooth solutions with sound speed $c(\rho)$ smoother than $C^{ {1}/{2}}$-H$\ddot{\rm o}$lder continuous  near vacuum states for compressible inviscid fluids can be found in \cite{chemin1,chemin2,24,25,MUK,Makino,38, 39}.
For the physical vacuum singularity that $c(\rho)$ is $C^{ {1}/{2}}$-H$\ddot{\rm o}$lder continuous across  vacuum boundaries, the standard approach of symmetric hyperbolic systems (cf. \cite{Friedrichs, Kato,17}) do not apply. This makes the study of  well-posedness of such problems in compressible fluids extremely  challenging and interesting, even for  the local-in-time existence theory. For compressible isentropic Euler equations in physical vacuum,  the characteristic speeds become singular with infinite spatial derivatives
at  vacuum boundaries that creates much  difficulties in analyzing the regularity near boundaries, so that the local-in-time well-posedness theory is established recently in \cite{16,10, 7, 10', 16'}. (See also \cite{zhenlei,zhenlei1, LXZ,serre} for related works on the local theory.)
The phenomena of physical vacuum singularity arise naturally in several important situations  besides the above mentioned, for example, the equilibrium and dynamics of boundaries of gaseous stars (cf. \cite{6', cox, HaJa2, LXZ}). A paramount motivation in the study of physical vacuum is to understand the long time stability of some physically important explicit solutions with scaling invariance such as Barenblatt self-similar solutions and affine motions. This requires obtaining long time higher order regularity of solutions near vacuum boundaries.  Extending from local-in-time existence to  long-time ones  is of fundamental importance in nonlinear problems and  poses a great challenge due to strong degenerate nonlinear hyperbolic characters.
It should be pointed out, for the Cauchy problem of the one-dimensional compressible Euler equations with damping,   the $L^p$-convergence of
$L^{\infty}$-weak solutions   to Barenblatt  solutions of the porous media equations was given in \cite{HMP} with  $p=2$ if $1<\ga\le 2$ and $p=\ga$ if $\ga>2$ and in \cite{HPW} with $p=1$, respectively, using   entropy-type estimates for the solution itself without deriving  estimates for derivatives.
However, the interfaces separating gases and vacuum cannot be traced in the framework of $L^{\infty}$-weak solutions.

\section{Reformulation of the problem and main results}
\subsection{Lagrangian variables, ansatz, and perturbations}

The domains of gases for the free boundary problem \ef{2.1} and the Barenblatt solutions, $\Omega(t)$ and $\bar\Omega(t)$, are
generally  different. In order to compare solutions defined on different domains, we reduce the problems to the ones defined on a common fixed domain, the initial domain of the Barenblatt solution,  $\Omega= \bar\Omega(0)$, which is  the ball  centered at the origin with the radius $\bar R(0)=\sqrt{\underline{A}/\underline{B}} $.

We  define $x$ as the Lagrangian flow of the velocity $u$ by
\be\label{171017}
\pl_t x(t,y)= u(t, x(t,y)) \ \  {\rm for} \ \  t>0,  \ \ {\rm and} \ \  x(0,y)=x_0(y) \ \ {\rm for} \ \ y\in \Omega,
\ee
and set the Lagrangian density, the inverse of the Jacobian matrix, and the Jacobian determinant by
\begin{align*}
\varrho(t,y)=\rho(t, x(t,y)),  \ \ \ \
\mathscr{A}(t,y)=\lt(\frac{\pl x}{\pl y}\rt)^{-1}, \ \ \  \
\mathscr{J}(t,y)={\rm det} \lt(\frac{\pl x}{\pl y}\rt)  .
\end{align*}
Then  system \ef{2.1} can be written in Lagrangian coordinates as
\begin{subequations}\label{2.1n}\begin{align}
& \pl_t  \varrho  +   \varrho \mathscr{A}_i^k \pl_t\pl_k x^i = 0 &  {\rm in}& \ \ \Omega\times(0,T], \label{2.1na}\\
 &   \varrho \pl_{tt} x_i   + \mathscr{A}_i^k \pl_k ( { \varrho^\ga} )  = -  \varrho \pl_t x_i  & {\rm in}& \ \ \Omega\times(0,T],\label{2.1nb}\\
 & \varrho>0 &{\rm in }  & \ \ \Omega\times(0,T],\label{2.1nc}\\
 &  \varrho=0    &    {\rm on}& \  \pl \Omega\times(0,T], \label{2.1nd}\\
&( \varrho,   x, \pl_t x)=(\rho_0(x_0), x_0, u_0(x_0)) & {\rm on} & \ \  \Omega\times \{t=0\}, \label{2.1nf}
 \end{align} \end{subequations}
where $x^i=x_i$ and $\pl_k=\frac{\pl}{\pl y_k}$. It follows from  \ef{2.1na} and
$\pl_t \mathscr{J} = \mathscr{J} \mathscr{A}_i^k \pl_t\pl_k x^i   $ that
$$ \varrho(t,y)\mathscr{J}(t,y)= \varrho (0,y)\mathscr{J}(0,y)=\rho_0\lt(x_0(y)\rt) {\rm det} \lt(\frac{\pl x_0(y)}{\pl y}\rt) .$$
We choose $x_0(y)$ such that
$\rho_0\lt(x_0(y)\rt) {\rm det} \lt(\frac{\pl x_0(y)}{\pl y}\rt)=  \bar\rho_0(y) $,
where $\bar\rho_0(y)=\bar\rho(0,y)$ is the initial density of the Barenblatt solution  given by \ef{1.6}. The existence of such an $x_0$ follows from the Dacorogna-Moser theorem (cf. \cite{DM}) and \ef{initial density}. It means that the Lagrangian density can be expressed as
\be\label{1.6'}
\varrho=\bar\rho_0 \mathscr{J}^{-1}, \ \ {\rm where} \ \
\bar\rho_0(y) =\lt(\underline{ A }- \underline{B} |{y}|^2 \rt)^{{1}/({\ga-1})},
   \ee
and problem \ef{2.1n} reduces to
\begin{subequations}\label{system}\begin{align}
&\bar\rho_0 \pl_{tt}  x_i    + \mathscr{J}  \mathscr{A}_i^k\pl_k  \lt(\bar\rho_0^\ga \mathscr{J}^{-\ga}\rt) = -  \bar\rho_0  \pl_t x_i  &  {\rm in}& \ \ \Omega\times(0,T], \label{171018} \\
&(   x, \pl_t x)=(x_0, u_0(x_0)) & {\rm on} & \ \  \Omega\times \{t=0\}.
\end{align}\end{subequations}

We define $\bar x$ as the Lagrangian flow of the  Barenblatt velocity $\bar u$ by
$$\pl_t \bar x(t,y)= \bar u(t, \bar x(t,y))  \ \ {\rm for} \ \ t>0, \ \ {\rm and} \ \   \bar x (0,y) =y  \ \ {\rm for} \ \  y\in \Omega.$$
then
$$\bar x(t,y)=\nu(t) y \ \ {\rm for} \ \  y\in  \Omega , \ \  {\rm where} \ \ \nu(t)=(1+t)^{\frac{1}{3\ga-1}} . $$
Since $\bar x(t, y)$ does not solve equation \ef{171018}, as in \cite{LZ, HZ}, we introduce a correction $h(t)$ which is the  solution to the following initial value problem of ordinary differential equations:
\be\label{pomt}\begin{split}
& h_{tt} + h_t -  (3\ga-1)^{-1}  (\nu+h)^{2-3\ga}  +  \nu_{ tt}  + \nu_{t}    =0, \ \ t >0,  \\
& h(t=0)=h_t(t=0)=0.
 \end{split}\ee
It should be noted that $\ta=\nu+h$ behaves like $\nu$. Precisely, there exist positive constants
$K$ and $C(n)$ independent of time $t$ such that for all $t\ge 0$,
\begin{subequations}\label{decay}\begin{align}
&\lt(1 +   t \rt)^{ {1}/({3\ga-1})} \le \ta(t) \le K \lt(1 +   t \rt)^{ {1}/({3\ga-1})},  \ \   \ \   \ta_t(t) \ge 0 , \label{decaya} \\
&\lt|\f{d^k }{dt^k}\ta(t)\rt| \le C(n)\lt(1 +   t \rt)^{\frac{1}{3\ga-1}-k},   \ \ k=1, 2,  \cdots, n, \label{decayb}
 \end{align}\end{subequations}
whose proof can be found in \cite{HZ}.
The new ansatz is then given by
\be\label{teeta}
\tilde x(t,y )=\bar x (t,y )  +   h(t) y =\ta(t) y, \ \ {\rm where} \ \ \ta(t)=\nu(t)+h(t),
\ee
which satisfies
 \be\label{equeta} \begin{split}
\bar\rho_0  \pl_{tt} \tilde x_i   +  \tilde{\mathscr{J}}\tilde{\mathscr{A}}_i^k  \pl_k\lt(\bar\rho_0^\ga \tilde{\mathscr{J}}^{-\ga}\rt)  = -  \bar\rho_0  \pl_t \tilde x_i \ \ {\rm in} \ \  \Omega\times(0,\iy),
\end{split}
\ee
where $\tilde{\mathscr{J}}= {\rm det} \lt(\frac{\pl \tilde x}{\pl y}\rt)=\ta^3$ and $\tilde{\mathscr{A}}=\lt(\frac{\pl \tilde x}{\pl y}\rt)^{-1} = \ta^{-1}\mathbf{Id}$.

We define the perturbation $\oa$ by
\be\label{perturbation}
\oa(t,y) =\ta^{-1}(t) \lt( x(t,y)- \tilde{x}(t,y)\rt)= \eta(t,y)-y , \ \  {\rm where} \  \ \eta(t,y)=\ta^{-1}(t)x(t,y),
\ee
then    \ef{171018}  can be expressed as
 \be\label{3-1-2}
\ta \bar\rho_0   \pl_t^2 \oa_i +   (\ta+ 2\ta_t ) \bar\rho_0 \pl_t  \oa_i +  (3\ga-1)^{-1}\ta^{2-3\ga} \bar\rho_0 \eta_i     + \ta^{2-3\ga} J A_i^k  \pl_k \lt(\bar\rho_0^\ga {J}^{-\ga}\rt)  = 0,
\ee
 where
\begin{align*}
A(t,y)  =\lt(\frac{\pl \eta}{\pl y}\rt)^{-1} =\lt(\mathbf{Id}+ \frac{\pl \oa}{\pl y}\rt)^{-1}    \ {\rm and}  \
 J(t,y)  = {\rm det} \lt(\frac{\pl \eta }{\pl y}\rt)= {\rm det} \lt(\mathbf{Id} +\frac{\pl \oa }{\pl y}\rt) .
\end{align*}
Problem \ef{system}, hence problems \ef{2.1n} and \ef{2.1}, can be written as
\begin{subequations}\label{newsystem}\begin{align}
&\ta \bar\rho_0   \pl_t^2 \oa_i +   (\ta+ 2\ta_t) \bar\rho_0 \pl_t  \oa_i +   (3\ga-1)^{-1}\ta^{2-3\ga} \bar\rho_0 \oa_i && \notag\\
& \qquad \qquad \qquad \qquad + \ta^{2-3\ga} \pl_k \lt(  \bar\rho_0^\ga \lt(A_i^k {J}^{1-\ga} -  \da_i^k\rt) \rt)   = 0 \ \ &  {\rm in} & \ \ \Omega\times(0,T], \label{3-1-3}\\
& (  \oa , \ \pl_t \oa)=(\ta^{-1}(0) x_0-y,\ \ta^{-1}(0) u_0(x_0)- \ta^{-2}(0)\ta_t(0)x_0)  \ \  & {\rm on}  \ \  & \Omega\times \{t=0\},
\end{align}\end{subequations}
due to \ef{3-1-2},  $\bar\rho_0(y) =\lt(\underline{ A }- \underline{B} |{y}|^2 \rt)^{{1}/({\ga-1})}$ and the Piola identity $\pl_k(JA^k_i)=0$ $(i=1,2,3)$.
In fact, equations \ef{3-1-2} and \ef{3-1-3} are useful, respectively, for curl estimates and energy estimates.

\subsection{Notation and main results}
We let $\pl_k=\frac{\pl}{\pl y_k}$,  $\pl^\al=\pl_1^{\al_1}\pl_2^{\al_2}\pl_3^{\al_3}$ for multi-index  $\al=(\al_1,\al_2, \al_3)$ and $\pl^j=\sum_{|\al|=j}\pl^\al$ for nonnegative integer $j$. We use $(\bar\pl_1,\bar\pl_2,\bar\pl_3)=y\times ( \pl_1, \pl_2, \pl_3)$
to denote  the angular momentum derivative, and let, similarly,
 $\bar\pl^\al=\bar\pl_1^{\al_1}\bar\pl_2^{\al_2}\bar\pl_3^{\al_3}$ for multi-index  $\al=(\al_1,\al_2, \al_3)$ and $\bar\pl^j=\sum_{|\al|=j}\bar\pl^\al$ for nonnegative integer $j$.

 The divergence and the $i$-th component of the curl of a vector filed $F$ are
$$
   {\rm div} F=\da^{k}_i \pl_k F^i, \ \ {\rm and} \ \   [ {\rm curl} F ]_i = \epsilon^{ijk} \pl_j F_k, \ \ i=1,2,3,
$$
where  $\epsilon^{ijk}$ is the standard permutation symbol given by
\begin{align*}
\epsilon^{ijk}=
\begin{cases}
&1,    \ \  \ \  \textrm{even permutation of} \ \ \{1,2,3\},   \\
&-1,   \ \  \  \textrm{odd permutation of} \ \ \{1,2,3\}, \\
&0,    \ \ \  \  \textrm{otherwise}.
\end{cases}
\end{align*}
Indeed, the angular momentum derivative can be written as $\bar\pl_i = \epsilon^{ijk} y_j \pl_k$.

Along the flow map $\eta$, the $i$-th component of the gradient  of a function $f$ is
$$ \lt[\nabla_\eta f \rt]_i = A^k_i \pl_k f,$$
the divergence and the $i$-th component of the curl of a vector filed $F$ are
$${\rm div}_\eta F = A^k_i \pl_k F^i, \  \ {\rm and} \ \
[ {\rm curl}_\eta F ]_{i}=\epsilon^{ijk} \lt[\nabla_\eta F_k \rt]_j   =\epsilon^{ijk} A^r_j \pl_r F_k, \ \  i=1,2,3.$$

Let
$$
\iota= ({\ga-1})^{-1} \ \ {\rm   and } \  \  \sigma(y) = \bar\rho_0^{\ga-1}=\underline{ A }- \underline{B} |{y}|^2. $$
 We introduce,   for nonnegative integers $m,n,l,j$,
\begin{align}
\mathscr{E}_j(t)= & \sum_{m+n+l=j} \lt\{ (1+t)^{2m+1} \lt\| \sa^{\frac{\iota+n}{2}}     \pl_t^{m+1}\pl^n \bar\pl^l \oa\rt\|_{L^2(\Omega)}^2  \rt. \notag \\
 & \lt. +   (1+t)^{2m }  \lt( \lt\| \sa^{\frac{\iota+n}{2}}     \pl_t^{m}\pl^n \bar\pl^l \oa\rt\|_{L^2(\Omega)}^2 +\lt\| \sa^{\frac{\iota+n+1}{2}}     \pl_t^{m}\pl^{n+1} \bar\pl^l \oa\rt\|_{L^2(\Omega)}^2  \rt) \rt\} , \label{7.6}
\end{align}
and  define the higher order weighted Sobolev norm $\mathscr{E}$ by
\begin{align}
\mathscr{E}(t)=\sum_{0\le j\le [\iota]+7} \mathscr{E}_j(t) . \label{7.6-1}
\end{align}
In addition to \ef{7.6-1}, we also need the following Sobolev norm for curl:
\begin{align}
\mathfrak{V}_{add}(t)
 =  \sum_{ 0\le m\le 1, \  0\le  n+l \le [\iota]+7} (1+t)^{2m} \lt\|\sa^{\frac{\iota+n+1}{2}}   \pl^{n}\bar\pl^{l}   {\rm curl}_\eta  \pl_t^m \oa \rt\|_{L^2(\Omega)}^2 . \label{7.20}
\end{align}

We are  now ready to state the main result.

\begin{thm}\label{mainthm}
There exist positive constants  $\bar\epsilon$ and $\bar\da$ depending only on the adiabatic exponent $\ga$ and the initial total mass $M$ such that,  for $ \mathscr{E}(0) + \mathfrak{V}_{add}(0) \le \bar\epsilon$, the life span of the unique smooth solution to  problem \ef{system} (hence to problem \ef{2.1}) exceeds $T_\iy$, where
\be\label{10.3}
T_\iy = \exp\lt\{\min\lt\{ \lt(\frac{\bar\da}{\mathscr{E}(0)}\rt)^{1/2}, \  \ \lt(\frac{\bar\da}{ \mathfrak{V}_{add}(0)} \rt)^{1/3}  \rt\} \rt\}-1.
\ee
\end{thm}

\begin{rmk} There exist positive constants  $\hat\epsilon$ and $\hat\da$ depending on  $\ga$ and  $M$ such that if  $\mathscr{E}(0) \le \bar\epsilon$ and
\be\label{condition}
\mathscr{V}  (t)\le \hat\da \mathscr{E}  (t)   \ \  {\rm for} \ \ 0\le t<T,
\ee
then problem \ef{system} (hence problem \ef{2.1}) admits a unique smooth solution in $[0, T)$  with 
$$\sup_{0\le t<T}\mathscr{E}(t) \le C \mathscr{E}(0)$$ 
for a certain constant $C$ independent of time $t$, where
$$
\mathscr{V} (t)= \sum_{ m+n+l\le [\iota]+7}
   (1+t)^{2m} \min\lt\{  \lt\| \sa^{\frac{\iota+n+1}{2}}  {\rm curl}  \pl_t^m\pl^n \bar\pl^l   \oa \rt\|_{L^2}^2, \ \lt\| \sa^{\frac{\iota+n+1}{2}}  \pl_t^m\pl^n \bar\pl^l   {\rm curl} \oa \rt\|_{L^2}^2  \rt\}.
   $$
If  condition \ef{condition} holds for $T=\infty$, then we have the global-in-time existence of smooth solutions.
This conclusion can be derived mainly from the estimates in Corollary \ref{10.8}.
Clearly, condition \ef{condition} holds if ${\rm curl} \oa=0$ ($\pl_i \oa_j-\pl_j \oa_i =0 $, $i,j=1,2,3$), or equivalently, ${\rm curl} x=0$ ($\pl_i x_j-\pl_j x_i =0 $, $i,j=1,2,3$), in particular, this condition is true for the spherically symmetric perturbations for $0\le t<\infty$.  The results obtained in this paper coincide with the ones in \cite{HZ}.

\end{rmk}

\begin{rmk} It should be noted that
\begin{align}
(1+t)^{2m}\sum_{
  m+2n+l\le 4} \lt\| \pl_t^m\pl^{n}\bar\pl^{l} \oa \rt\|_{L^\iy}^2  \le C \mathscr{E}(t)
\end{align}
for constant $C$ depending only on $M$ and $\ga$. The high order norm $\mathscr{E}$ has been defined to have the fewest derivatives to ensure that $\pl_t^2 \pl \oa  $ is pointwise bounded, a requirement for the curl estimate, which is easy to see from the curl equation:
\be\label{rmk-curl}\ta {\rm curl}_\eta (\pl_t^2 \oa) +(2\ta_t + \ta) {\rm curl}_\eta    \pl_t \oa =0.
\ee
Because the regularity that $\pl_t \pl \oa  $ is pointwise bounded is needed  at least  to ensure that a solution to problem \ef{system} is also a solution to  problem \ef{2.1}.
\end{rmk}

\begin{rmk}
The reason why the perturbation is chosen as $\oa(t,y) = \ta^{-1}(t) \lt( x(t,y)- \tilde{x}(t,y)\rt)$, instead of, $\za(t,y)=x(t,y)- \tilde{x}(t,y)$, is as follows. The curl equation for $\za$ is
\be\label{9.27}
 {\rm curl}_x  \pl_t^2 \za  +   {\rm curl}_x  \pl_t  \za = ({3\ga-1})^{-1}\ta^{1-3\ga} {\rm curl}_x    \za ,
\ee
where $[ {\rm curl}_x F ]_{i}=\epsilon^{ijk} \mathscr{A}^r_j \pl_r F_k $  $(i=1,2,3)$ for any vector $F$. Since $\ta^{1-3\ga}(t)$ is equivalent to $(1+t)^{-1}$, the accumulation of the term on the right hand side of \ef{9.27} in time cannot be controlled easily. However, this bad term can be absorbed in the curl equation for $\oa$, \ef{rmk-curl}.

\end{rmk}

\section{A priori estimate}

The proof of Theorem \ref{mainthm} is based on the following a priori estimates, together with the local existence theory (cf. \cite{10',16'}).

\begin{thm}\label{thm3.1}
Let $\oa(t,y)$ be a solution to  problem \ef{newsystem} in the time interval $[0, T]$ satisfying the following a priori assumptions:
\begin{subequations}\label{initial}\begin{align}
& \mathscr{E}(t) \le \ea_0^2 , &\ \  t\in [0,T], \label{assume}\\
& \lt( \ln  (1+t) \rt)^2 \sup_{s\in [0,t]}  \mathscr{E}(s)  \le \ea_0^2 , & \ \  t\in [0,T],
\end{align}\end{subequations}
then
\begin{align}\label{energy}
  \mathscr{E}(t) + \mathfrak{V}_{add}(t) \le C \lt(   \mathscr{E}(0) + \mathfrak{V}_{add}(0) +  \ln(1+t) \mathfrak{V}_{add}(0) \rt), \ \  t\in [0,T],
\end{align}
where $C$ is a positive constant independent of $t$.
\end{thm}

To simplify the presentation, we introduce some notation. Throughout the rest of paper,   $C$  will denote a positive constant which  only depend on the parameters of the problem,
$\ga$ and $M$, but does not depend on the data. They are referred as universal and can change
from one inequality to another one. Also we use $C(\beta)$ to denote  a certain positive constant
depending on quantity $\beta$.
We will employ the notation $a\lesssim b$ to denote $a\le C b$, $a  \thicksim b$ to denote $C^{-1}b\le a\le Cb$ and $a \gtrsim b$ to denote $a \ge  C^{-1} b$,
where  $C$ is the universal constant  as defined
above.

Recall that $\sigma(y) = \bar\rho_0^{\ga-1}=\underline{ A }- \underline{B} |{y}|^2$, and $ \Omega= \bar\Omega(0)$ is the ball centered at the origin with the radius $ \sqrt{\underline{A}/\underline{B}} $. So, $\sa(y)$ is equivalent to $d(y,\pl\Omega)$, the distance function to the boundary of $\Omega$, that is, $\sa(y) \thicksim d(y,\pl\Omega)$.
We will use, in the rest of this work, the notation
$$ \int=\int_{\Omega}, \ \  {\rm and} \ \ \|\cdot\|_{W^{k,p}}=\|\cdot\|_{W^{k,p}(\Omega)}  $$
for $k\ge 0$ and $p\in[1,\infty]$.

\subsection{Basic inequalities I}
In this subsection, we will show the bounds derived from  $\mathscr{E}(t)$.
Indeed, it holds that
\begin{align*}
&\sum_{
  m+2n+l\le 4} \lt\| \pl_t^m\pl^{n}\bar\pl^{l} \oa \rt\|_{L^\iy(\Omega)}    + \sum_{\substack{
  m+2n+l =5 }}  \lt\|  \pl_t^m\pl^{n}\bar\pl^{l} \oa \rt\|_{H^1(\Omega)}  \notag\\
& + \sum_{\substack{
6 \le  m+2n+l \\
m+n+l\le  [\iota]+6 }}  \lt\|\sa^{\frac{m+2n+l-4}{2}}\pl_t^m\pl^{n}\bar\pl^{l} \oa \rt\|_{L^\iy(\Omega)} \les (1+t)^{-  m} \sqrt{\mathscr{E}(t)},
\end{align*}
whose proof will be given later in Lemma \ref{very}, based on the following Hardy inequalities and weighted Sobolev embeddings.

Let $k>-1$ be a given real number, $\da$ be a positive constant,  and $f$ be a function satisfying
$
\int_0^{\da} r^{k+2} \lt(f^2 + |f'|^2\rt) dr < \iy,
$
then it holds that
\be\label{hardy'}
\int_0^{\da} r^{k } f^2 dr \le C(\da,k)  \int_0^{\da} r^{k+2} \lt(f^2 + |f'|^2\rt) dr,
\ee
whose proof  can be found  in \cite{18'}. Indeed, \ef{hardy'} is a   general version of the  standard Hardy inequality:
$\int_0^\iy |r^{-1} f|^2 dr \le C \int_0^\iy |f'|^2 dr$.

As a consequence of \ef{hardy'}, we have the following  estimates.
\begin{lem}
Let $k>-1$ be a given real number,   and $f$ be a function satisfying $\int_\Omega \sa^{k+2} (f^2 + |\pl f|^2 )dy <\iy$, then it holds that
\begin{align}\label{hard}
\int_\Omega \sa^k f^2 dy \le C(k,\Omega ) \int_\Omega \sa^{k+2} (f^2 + |\pl f|^2) dy.
\end{align}
\end{lem}

{\em Proof}. Let  $\Omega^0$ be a ball centered at the origin with the radius $ \sqrt{\underline{A}/(4\underline{B})}$, and   $\Omega^b= \Omega \setminus \Omega^0$. In  $\Omega^0$, $\sa$ has  positive upper and lower  bounds so that
$$\int_{\Omega^0} \sa^k f^2 dy \le  C(k,\Omega ) \int_{\Omega^0} \sa^{k+2}  f^2  dy  \le  C(k,\Omega ) \int_{\Omega} \sa^{k+2}  f^2  dy.$$
Near the boundary, we may write the integral in the spherical coordinates $(r,\phi, \psi)$:
$$
\int_{\Omega^b} \sa^k(y) f^2(y) dy=\int_{\sqrt{\underline{A}/(4\underline{B})}}^{\sqrt{\underline{A}/ \underline{B} }} \sa^k(r)    F^2(r)  dr,
$$
where
$$
F^2 (r)=   \int_0^\pi \int_0^{2\pi}   f^2(r,\phi, \psi) r^2 \sin \phi  d\phi d\psi.
$$
This, together with \ef{hardy'},  the equivalence of $\sa$ and $\sqrt{\underline{A}/ \underline{B} }-r$,  and  the  H\"{o}lder  inequality, implies that  for some constants  $C= C(k,\Omega )$,
\begin{align*}
&\int_{\Omega^b} \sa^k(y) f^2(y) dy
 \le C \int_{\sqrt{\underline{A}/(4\underline{B})}}^{\sqrt{\underline{A}/ \underline{B} }} \sa^{k+2}(r)   ( F^2 + |F_r|^2 )(r)  dr \\
 \le & C \int_{\sqrt{\underline{A}/(4\underline{B})}}^{\sqrt{\underline{A}/ \underline{B} }} \sa^{k+2}(r)  \lt(  \int_0^\pi \int_0^{2\pi}  ( f^2 + |\pl_r f|^2 ) r^2 \sin \phi  d\phi d\psi \rt)  dr\\
\le & C \int_{\sqrt{\underline{A}/(4\underline{B})}}^{\sqrt{\underline{A}/ \underline{B} }} \sa^{k+2}(r)  \lt(  \int_0^\pi \int_0^{2\pi}  ( f^2 + |\pl  f|^2 ) r^2 \sin \phi  d\phi d\psi \rt)  dr.
\end{align*}
It proves \ef{hard} by writing the last integral in the coordinates $y=(y_1,y_2,y_3)$.
\hfill $\Box$

Let $\mathfrak{U}$ be a bounded smooth domain in $\mathbb{R}^3  $, and $d=d(y)=dist(y, \partial  \mathfrak{U})$ be a distance function to the boundary.
 For any  $a>0$ and nonnegative integer $b$,  we define the  weighted Sobolev space  $H^{a, b}(  \mathfrak{U})$   by
$$ H^{a, b}(\mathfrak{U}) = \lt\{   d^{a/2}f \in L^2(\mathfrak{U}): \ \  \int_\mathfrak{U}    d^a|\pl^k f|^2dy<\infty, \ \  0\le k\le b\rt\}$$
  with the norm
$ \|f\|^2_{H^{a, b}(\mathfrak{U})} = \sum_{k=0}^b \int_\mathfrak{U}    d^a|\pl^k f|^2dy$. Let $H^s( \mathfrak{U} )$ $(s\ge 0)$ be the standard Sobolev space, then  for $b\ge  {a}/{2}$, we have the following embedding of weighted Sobolev spaces (cf. \cite{18'}):
 $ H^{a, b}(\mathfrak{U} )\hookrightarrow H^{b- {a}/{2}}( \mathfrak{U})$
    with the estimate
  \be\label{wsv} \|f\|_{H^{b- {a}/{2}}( \mathfrak{U})} \le C(a,b,\mathfrak{U})  \|f\|_{H^{a, b}(\mathfrak{U} )} .\ee

As a conclusion of \ef{hard} and \ef{wsv}, we have the following  estimates.
\begin{lem}\label{very}
Let $m$ be nonnegative integers, $\al$ and $\ba$ be multi-indexes.
Suppose that $\mathscr{E}(t)$  is finite, then it holds that
\begin{align}
&\sum_{
  m+2|\al|+|\ba|\le 4} \lt\| \pl_t^m\pl^{\al}\bar\pl^{\ba} \oa \rt\|_{L^\iy(\Omega)}^2   + \sum_{\substack{
  m+2|\al|+|\ba| =5 }}  \lt\|  \pl_t^m\pl^{\al}\bar\pl^{\ba} \oa \rt\|_{H^1(\Omega)}^2  \notag\\
& + \sum_{\substack{
6 \le m+2|\al|+|\ba| \\
m+|\al|+|\ba|\le  [\iota]+6 }}  \lt\|\sa^{\frac{m+2|\al|+|\ba|-4}{2}}\pl_t^m\pl^{\al}\bar\pl^{\ba} \oa \rt\|_{L^\iy(\Omega)}^2
  \les (1+t)^{-2 m}   \mathscr{E}(t) . \label{verify}
\end{align}
 \end{lem}

{\em Proof}. When $m+2|\al|+|\ba|\le 4$, we have
\begin{align}
&(1+t)^{ 2m} \lt\|\pl_t^m \pl^\al \bar\pl^\ba \oa\rt\|_{H^{\iota+[\iota]+8-m-|\ba|, \ [\iota]+8- m - |\al|-|\ba| }(\Omega)}^2 \notag\\
 =  & (1+t)^{ 2m}  \sum_{0\le |h| \le [\iota]+8- m - |\al|-|\ba|}\lt\|\sa^{\frac{\iota+[\iota]+8-m-|\ba|}{2}} \pl^h \pl_t^m \pl^\al \bar\pl^\ba \oa \rt\|_{L^2(\Omega)}^2 \notag\\
\les &  (1+t)^{ 2m}  \sum_{0\le |h| \le [\iota]+8- m - |\al|-|\ba|} \lt\|\sa^{\frac{\iota+ |\al|+|h| }{2}} \pl_t^m \pl^{\al+h} \bar\pl^\ba \oa \rt\|_{L^2(\Omega)}^2  \notag\\
\les &  \mathscr{E}_{m+|\al|+|\ba|} + \sum_{1\le |h| \le [\iota]+8- m - |\al|-|\ba|}
\mathscr{E}_{m+|\al|+|h|-1+|\ba|}
\le    \mathscr{E}, \label{8.22d}
\end{align}
which, together with  \ef{wsv} and the fact that $H^q(\Omega)\hookrightarrow L^\infty(\Omega)$ for $q>3/2$, implies that
\begin{align*}
& \| \pl_t^m\pl^{\al}\bar\pl^{\ba} \oa \|_{L^\iy(\Omega)}^2  \les   \| \pl_t^m\pl^{\al}\bar\pl^{\ba} \oa \|_{H^{\frac{3}{2} + \frac{ (1 + [\iota] - \iota ) +  (4-m-2|\al|-|\ba| )  }{2} }(\Omega)}^2  \notag \\
\les  &  \lt\|\pl_t^m \pl^\al \bar\pl^\ba \oa\rt\|_{H^{\iota+[\iota]+8-m-|\ba|, \ [\iota]+8- m - |\al|-|\ba|}(\Omega)}^2
 \les  (1+t)^{-2m} \mathscr{E}.
\end{align*}
Similarly, we can obtain for $m+2|\al|+|\ba| \ge 6$ and $m+|\al|+|\ba|\le  [\iota]+6$,
\begin{align}
&\lt\|\sa^{\frac{m+2|\al|+|\ba|-4}{2}}\pl_t^m\pl^{\al}\bar\pl^{\ba} \oa \rt\|_{L^\iy (\Omega)}^2 \les \lt\|\sa^{\frac{m+2|\al|+|\ba|-4}{2}}\pl_t^m\pl^{\al}\bar\pl^{\ba} \oa \rt\|_{H^{\frac{3}{2} + \frac{1 + [\iota]-  \iota   }{4} } (\Omega)}^2 \notag \\
 \les &
\lt\|\sa^{\frac{m+2|\al|+|\ba|-4}{2}}\pl_t^m\pl^{\al}\bar\pl^{\ba} \oa \rt\|_{H^{ \frac{\iota+3[\iota]+25-4m -4|\al|-4|\ba|}{2}, \ [\iota]+8- m -|\al|-|\ba| } (\Omega)}^2
\les (1+t)^{-2m} \mathscr{E} . \label{6.22}
\end{align}
Indeed, the derivation of the last inequality in  \ef{6.22} is not trivial, which is based on \ef{hard}. We only  examine  the  difficult case where $m+2|\al|+|\ba|\ge 7$.
\begin{align*}
&\lt\|\sa^{\frac{m+2|\al|+|\ba|-4}{2}}\pl_t^m\pl^{\al}\bar\pl^{\ba} \oa \rt\|_{H^{ \frac{\iota+3[\iota]+25-4m -4|\al|-4|\ba|}{2}, \ [\iota]+8- m -|\al|-|\ba| } (\Omega)}^2
\\
\les
& \sum_{\substack{
0\le |h| \le [\iota]+8- m - |\al|-|\ba| \\
0 \le j\le |h|
}} \lt\|\sa^{ \frac{\iota+3[\iota]+17-2m -2|\ba|}{4}-j } \pl_t^m\pl^{|\al|+|h|-j}\bar\pl^{|\ba|} \oa  \rt\|_{L^2(\Omega)}^2
\\
\les
& \sum_{\substack{
0\le |h| \le [\iota]+8- m - |\al|-|\ba| \\
0 \le j\le |h|
}} \lt\|\sa^{ \frac{\iota+3[\iota]+17-2m -2|\ba|}{4} } \pl_t^m\pl^{|\al|+|h|-j}\bar\pl^{|\ba|} \oa  \rt\|_{L^2(\Omega)}^2
\\
\les
& \sum_{\substack{
0\le |h| \le [\iota]+8- m - |\al|-|\ba| \\
0 \le j\le |h|
}} \lt\|\sa^{\frac{\iota+ |\al|+|h| }{2} } \pl_t^m\pl^{|\al|+|h|-j}\bar\pl^{|\ba|} \oa \rt\|_{L^2(\Omega)}^2 \les  (1+t)^{-2m} \mathscr{E},
\end{align*}
where  \ef{hard} has been used $j$ times to derive the second inequality. When $m+2|\al|+|\ba|=5$, it follows from \ef{wsv} that
\begin{align*}
\lt\|  \pl_t^m\pl^{\al}\bar\pl^{\ba} \oa \rt\|_{H^1(\Omega)}^2
\les  \lt\|  \pl_t^m\pl^{\al}\bar\pl^{\ba} \oa \rt\|_{H^{1+\frac{1+[\iota]-\iota}{2}}(\Omega)}^2 =  \lt\|  \pl_t^m\pl^{\al}\bar\pl^{\ba} \oa \rt\|_{H^{ \frac{ [\iota]+8-\iota-(m+2|\al|+|\ba|)}{2}}(\Omega)}^2\\
\les \lt\|  \pl_t^m\pl^{\al}\bar\pl^{\ba} \oa \rt\|_{H^{\iota+[\iota]+8-m-|\ba|, \ [\iota]+8- m - |\al|-|\ba|}(\Omega)}^2  \les  (1+t)^{-2m} \mathscr{E},
\end{align*}
where the last inequality follows from the same derivation of \ef{8.22d}.
\hfill $\Box$

\subsection{Basic inequalities II}

Since $JA$ is the adjugate matrix of $(\frac{\pl \eta}{\pl y})$ and $\eta(t,y)=\oa(t,y) + y$, then
\be\label{6.7-2}
J A = \lt(\frac{\pl\eta}{\pl y} \rt)^* = \lt[\begin{split} \pl_2 \eta  \times \pl_3 \eta \\
\pl_3 \eta  \times  \pl_1 \eta  \\
\pl_1 \eta  \times \pl_2 \eta \end{split}\rt] =   \lt(1 +{\rm div}\oa\rt)\mathbf{Id} -  \lt(\frac{\pl\oa}{\pl y}\rt) +  b,
\ee
where $b$ is the adjugate matrix of $(\frac{\pl \oa}{\pl y})$ given by
\bee
b= \lt(\frac{\pl\oa}{\pl y} \rt)^* = \lt[\begin{split} \pl_2 \omega  \times \pl_3 \omega \\
\pl_3 \omega  \times  \pl_1\omega  \\
\pl_1 \omega  \times \pl_2 \omega \end{split}\rt]  .
\eee
This, together with the fact that $(\frac{\pl \eta}{\pl y})(\frac{\pl \eta}{\pl y})^*=J\mathbf{Id} $, implies that
 \begin{align}
 J=1+{\rm div} \oa + 2^{-1}\lt(|{\rm div} \oa|^2 + |{\rm curl} \oa|^2- |\pl \oa|^2\rt) + b^s_r \pl_s \oa^r . \label{3-19-2}
\end{align}

Due to \ef{assume} and \ef{verify},  we have for $t\in [0,T]$,
\be\label{a}
|\pl\oa(t,y)| \les \ea_0.
\ee
Thus, it follows from   \ef{3-19-2} and \ef{6.7-2}  that  for $t\in [0,T]$,
\be\label{7.9}
|J-1| \les |\pl \oa| \les \ea_0 \ \ {\rm and} \ \   \|A-{\bf Id}\|_{L^\iy} \les |\pl \oa| \les \ea_0,
\ee
which implies, with the aid of the smallness of $\ea_0$, that  for $t\in [0,T]$
 \be\label{6.7-1a}
  2^{-1}\le J \le 2   \ \  {\rm and} \ \
    \|A \|_{L^\iy} \le 2  .
    \ee
Indeed,    $2^{-1}\le J \le 2$ follows from $|J-1|   \les \ea_0$, $\|A-{\bf Id}\|_{L^\iy} \les |\pl \oa|$ follows from \ef{6.7-2}, $ 2^{-1}\le J$ and $|J-1| \les |\pl \oa|$.
Moreover,  we have for any function $f$
$$
   |[\na_\eta f]_i -\pl_i f| =|( A^r_i - \da^r_i) \pl_r f| \les \ea_0 |\pl f|,
$$
which means
 \be\label{6.7-1c}
   2^{-1} |\pl f| \le |\na_\eta f| \le 2 |\pl f|.
\ee

\section{energy estimates}
We let $m,n,l,j$ be nonnegative integers, and
introduce  the following  $j$-th order   energy functional $\mathfrak{E}_j$ and   dissipation functional $\mathfrak{D}_j$:
\begin{align*}
&    \mathfrak{E}_j(t)= \sum_{m+n+l=j} \mathfrak{E}^{m,n,l}(t) = \sum_{m+n+l=j}  \lt( \mathfrak{E}_{I}^{m,n,l}   + \mathfrak{E}_{II}^{m,n,l} \rt)(t) , \\
&   \mathfrak{D}_j (t) =  \sum_{m+n+l=j} \mathfrak{D}^{m,n,l}(t) = \sum_{m+n+l=j} \lt(  \mathfrak{E}_{I}^{m,n,l}  + (1+t)^{-1}  \mathfrak{E}_{II}^{m,n,l}\rt)(t) ,
\end{align*}
where
\begin{align*}
& \mathfrak{E}_{I}^{m,n,l}(t) = (1+t)^{2m+1} \lt\|  \sa^{\frac{\iota+n}{2}}   \pl_t^{m+1}\pl^n \bar\pl^l \oa \rt\|_{L^2}^2, \\
&\mathfrak{E}_{II}^{m,n,l}(t) = (1+t)^{2m} \lt( \lt\|  \sa^{\frac{\iota+n}{2}}   \pl_t^{m}\pl^n \bar\pl^l \oa \rt\|_{L^2}^2 + \lt\|  \sa^{\frac{\iota+n+1}{2}}   \na_\eta \pl_t^{m}\pl^n \bar\pl^l \oa \rt\|_{L^2}^2 \rt. \notag\\
& \qquad \qquad \qquad \qquad  \qquad \lt.
 + \iota^{-1}  \lt\|  \sa^{\frac{\iota+n+1}{2}}   {\rm div}_\eta  \pl_t^{m}\pl^n \bar\pl^l \oa \rt\|_{L^2}^2 \rt).
\end{align*}
Let $\oa(t,y)$ be a solution to problem \ef{newsystem} in the time interval $[0, T]$ satisfying \ef{assume}, then
it is easy to see the equivalence of the  weighted Sobolev norm $\mathscr{E}$ and the energy functional $\mathfrak{E}=\sum_{0\le j \le [\iota]+7} \mathfrak{E}_j $. Indeed, it follows from \ef{6.7-1c} that
\begin{align}\label{7.10}
\mathscr{E}_j  \thicksim \mathfrak{E}_j,\ \     j=0, 1, \cdots, [\iota]+7.
 \end{align}
In addition to $\mathfrak{E}_j$ and $\mathfrak{D}_j$, we introduce the $j$-th order  weighted Sobolev norm for curl:
$$
  \mathfrak{V}_j (t) =  \sum_{m+n+l=j} \mathfrak{V}^{m,n,l} (t) = \sum_{m+n+l=j}
   (1+t)^{2m}   \lt\| \sa^{\frac{\iota+n+1}{2}}   {\rm curl}_\eta \pl_t^m\pl^n \bar\pl^l \oa \rt\|_{L^2}^2.
$$

Now, we have the following estimates.

\begin{prop}\label{newedv} Let $\oa(t,y)$ be a solution to problem \ef{newsystem} in the time interval $[0, T]$ satisfying \ef{assume}. Then for $j=0, 1,2,\cdots, [\iota]+7 $,
\begin{align}\label{new-hig}
   \mathfrak{E}_{j}(t) + \int_0^t \mathfrak{D}_{j}(s)ds \les  \sum_{0\le k \le j }\lt( \mathfrak{E}_k  (0)+   \mathfrak{V}_k  (t)  + \int_0^t(1+s)^{-1} \mathfrak{V}_k (s)ds \rt), \  \  t\in [0,T].
 \end{align}
\end{prop}

The proof  consists of Lemmas \ref{est-low} and \ref{est-high}, which we will prove later in this section.
Based on Proposition \ref{newedv} and the fact that
\begin{align*}
&  \mathfrak{V}_j (t) \les \sum_{m+n+l=j}
   (1+t)^{2m}   \lt\| \sa^{\frac{\iota+n+1}{2}}  {\rm curl}  \pl_t^m\pl^n \bar\pl^l   \oa \rt\|_{L^2}^2 + \ea_0^2 \mathscr{E}_{j}(t) ,\\
&
  \mathfrak{V}_j (t) \les  \sum_{m+n+l=j}
   (1+t)^{2m}   \lt\| \sa^{\frac{\iota+n+1}{2}}  \pl_t^m\pl^n \bar\pl^l   {\rm curl} \oa \rt\|_{L^2}^2  + \ea_0^2 \mathscr{E}_{j}(t)+ \sum_{0\le k\le j-1} \mathscr{E}_{k}(t) ,
\end{align*}
due to \ef{7.9} and the commutator estimate \ef{commutator2} which will be proved later, we can use   \ef{7.10} and  the  mathematical induction to prove
\begin{coro}\label{10.8} Let $\oa(t,y)$ be a solution to problem \ef{newsystem} in the time interval $[0, T]$ satisfying \ef{assume}. Then for $j=0, 1,2,\cdots, [\iota]+7 $,
\begin{align*}
   & \mathscr{E}_{j}(t) + \int_0^t (1+s)^{-1}\mathscr{E}_{j}(s)ds \notag\\
    \les & \sum_{0\le k \le j }\lt( \mathscr{E}_k  (0)+   \mathscr{V}_k  (t)  + \int_0^t(1+s)^{-1} \mathscr{V}_k (s)ds \rt), \  \  t\in [0,T],
 \end{align*}
 where
$$
\mathscr{V}_k  (t)= \sum_{m+n+l=k}
   (1+t)^{2m} \min\lt\{  \lt\| \sa^{\frac{\iota+n+1}{2}}  {\rm curl}  \pl_t^m\pl^n \bar\pl^l   \oa \rt\|_{L^2}^2, \ \ \lt\| \sa^{\frac{\iota+n+1}{2}}  \pl_t^m\pl^n \bar\pl^l   {\rm curl} \oa \rt\|_{L^2}^2  \rt\}.
   $$
\end{coro}

\subsection{The zeroth order  estimate}\label{sec4.1}
In this subsection, we prove that
\begin{lem}\label{est-low} Let $\oa(t,y)$ be a solution to problem \ef{newsystem} in the time interval $[0, T]$ satisfying \ef{assume}. Then,
\begin{align}\label{low}
\mathfrak{E}_0 (t) + \int_0^t  \mathfrak D_0(s) ds \les & \mathfrak{E}_0(0)+  \mathfrak{V}_0(t)   + \int_0^t (1+s)^{-1} \mathfrak{V}_0(s) ds,  \  \  t\in [0,T].
\end{align}
\end{lem}

{\em Proof}. Multiply  \ef{3-1-3} by $\ta^{-1}$ and use $\bar\rho_0=\sa^\iota$ to obtain
\begin{align}
\sigma^\iota   \pl_t^2 \oa_i +   (1+2\ta^{-1}\ta_t) \sigma^\iota  \pl_t \oa_i  +   (3\ga-1)^{-1}\ta^{1-3\ga} \sigma^\iota \oa_i  \notag \\
  + \ta^{1-3\ga} \pl_k \lt( \sigma^{\iota+1}  (A_i^k {J}^{1-\ga} -  \da_i^k ) \rt)   = 0. \label{3.4}
\end{align}
Integrate the product  of \ef{3.4} and $\pl_t \oa^i$  over $\Omega$ and use $\pl_t J=J A^k_i \pl_t \pl_k  \oa^i $ to get
\begin{align}\label{3-19-1}
&\f{1}{2}\f{d}{dt} \int \lt\{  \sa^\iota |\pl_t \oa|^2 + \ta^{1-3\ga}\lt(   (3\ga-1)^{-1} \sa^\iota  |\oa|^2  +   2 \sa^{\iota+1} \mathcal{M}_0 \rt)  \rt\} dy
\notag \\
& +  (1+2\ta^{-1}\ta_t)  \int \sa^\iota  |\pl_t \oa|^2 dy
  =       \frac{1}{2} \lt( \ta^{ 1-3\ga} \rt)_t \int \lt(   (3\ga-1)^{-1} \sa^\iota  |\oa|^2  +   2 \sa^{\iota+1} \mathcal{M}_0 \rt) dy,
 \end{align}
where $\mathcal{M}_0= ({\ga-1})^{-1}  ( J^{1-\ga} -1 ) + {\rm div} \oa$. Due to \ef{3-19-2}, $\mathcal{M}_0$   can be rewritten as
$$
  \mathcal{M}_0=  2^{-1} \lt( |\pl \oa|^2 + (\ga-1)|{\rm div} \oa|^2 - |{\rm curl} \oa|^2\rt) +e_0, $$
 where $e_0$ represents the cubic term, given by
\begin{align*}
  e_0=  &({\ga-1})^{-1}\lt(J^{1-\ga}-1-(1-\ga)(J-1)+2^{-1}(1-\ga)\ga(J-1)^2\rt)\\
& + 2^{-1} {\ga} \lt((J-1)^2-|{\rm div} \oa|^2\rt)-b^s_r \pl_s \oa^r .
 \end{align*}
This, together with  \ef{3-19-2}, the Taylor expansion    and \ef{a}, implies that
\be\label{11-23}
 0\le  \f{1}{4} |\pl \oa|^2 + \f{\ga-1}{2}|{\rm div} \oa|^2   \le \mathcal{M}_0 +  \f{1}{2} |{\rm curl} \oa|^2 \le   |\pl \oa|^2 + \f{\ga-1}{2}|{\rm div} \oa|^2 .
\ee
It follows from \ef{3-19-1}, \ef{11-23}  and $\ta_t\ge 0$ that
\begin{align}\label{3.4-1}
 \f{1}{2}\f{d}{dt} \int \lt\{  \sa^\iota |\pl_t \oa|^2 + \ta^{1-3\ga}\lt(   (3\ga-1)^{-1} \sa^\iota  |\oa|^2  +   2 \sa^{\iota+1} \mathcal{M}_0 \rt)  \rt\} dy
 \notag \\
   +     \int \sa^\iota  |\pl_t \oa|^2 dy  \le   - \frac{1}{2} \lt( \ta^{ 1-3\ga} \rt)_t \int  \sa^{\iota+1}  |{\rm curl} \oa|^2 dy.
 \end{align}
Integrate \ef{3.4-1} over $[0,t]$ and use \ef{decay} and \ef{11-23} to obtain the basic estimate:
\begin{align}\label{3.4-2}
& (E_{0I}+E_{0II})(t)+\int_0^t E_{0I}(s) ds \notag\\
\les & (E_{0I}+E_{0II})(0) +V_0(t) + \int_0^t (1+s)^{-1} V_0(s) ds ,
\end{align}
where
\begin{align*}
& E_{0I}(t) = \int \sa^\iota  |\pl_t \oa|^2 dy,\  \ \  \   V_0(t)=(1+t)^{-1} \int  \sa^{\iota+1}  |{\rm curl} \oa |^2 dy,\\
& E_{0II}(t)=   (1+t)^{-1}  \int  \sa^{\iota }\lt(  |\oa|^2 + \sa  |\pl \oa|^2 + \iota^{-1}   \sa  |{\rm div}   \oa|^2  \rt)dy.
\end{align*}

To improve the  estimate \ef{3.4-2}, we integrate the product  of \ef{3.4} and $ \oa^i$  over $\Omega$ to get
\begin{align}
&\frac{d}{dt}\int \sa^\iota \lt( \oa^i \pl_t \oa_i + (2^{-1}+\ta^{-1}\ta_t)  | \oa|^2  \rt) dy
\notag \\
& +    \ta^{1-3\ga}   \int \lt\{ (3\ga-1)^{-1}\sa^\iota |\oa|^2 -  \sigma^{\iota+1} \lt(A_i^k {J}^{1-\ga} -  \da_i^k\rt) \pl_k \oa^i  \rt\} dy \notag\\
  =
  &\int \sa^\iota |\pl_t \oa|^2 dy + ( \ta^{-1}\ta_t)_t\int   \sa^\iota | \oa|^2 dy.   \label{3.4-3}
\end{align}
It follows from \ef{6.7-2}  that
\be\label{linear}
\da^k_i- A_i^k {J}^{1-\ga}  = \da^k_i- \lt(\lt(1 +{\rm div}\oa\rt)\da^k_i -  \pl_i\oa^k +  b^k_i\rt){J}^{-\ga}  = (\ga-1)     {\rm div} \oa \da^k_i   +   \pl_i  \oa^k  - Q^k_i,
 \ee
where $Q$ represents the quadratic term, given by
\begin{align*}
 & Q^k_i=(J^{-\ga} - 1 + \ga {\rm div} \oa )\da^k_i + ( J^{-\ga}-1)({\rm div} \oa \da^k_i - \pl_i \oa^k) + J^{-\ga} b^k_i  .
 \end{align*}
This, together with  \ef{3-19-2}, the Taylor expansion    and \ef{a}, gives
\begin{align*}
& - \int  \sigma^{\iota+1} \lt(A_i^k {J}^{1-\ga} -  \da_i^k\rt) \pl_k \oa^i  dy \\
  \ge
  & \int  \sigma^{\iota+1} \lt(2^{-1} |\pl \oa|^2 + (\ga-1) |{\rm div} \oa|^2 - |{\rm curl} \oa|^2\rt) dy.
\end{align*}
Due to the Cauchy inequality and \ef{decay} (especially, $\ta_t\ge 0$), we have
\begin{align*}
 \int \sa^\iota \lt(4^{-1}| \oa|^2 -   |\pl_t \oa |^2      \rt) dy \le &\int \sa^\iota \lt( \oa^i \pl_t \oa_i + (2^{-1}+\ta^{-1}\ta_t)  | \oa|^2  \rt) dy \\
\les & \int \sa^\iota \lt( |\pl_t \oa |^2  +  | \oa|^2  \rt) dy.
\end{align*}
Thus, integrating \ef{3.4-3} over $[0,t]$, and using \ef{decay} (especially, $ |( \ta^{-1}\ta_t)_t|\les (1+t)^{-2}$), \ef{3.4-2} and the Grownwall inequality, we obtain that
\begin{align}\label{3.4-4}
 \int \sa^\iota  |  \oa|^2 dy + \int_0^t E_{0II}(s) ds
  \les  (E_{0I}+E_{0II})(0) +   V_0(t) + \int_0^t V_0(s) ds.
\end{align}

Finally, we integrate the product of $1+t$ and \ef{3.4-1} over $[0,t]$, and use \ef{decay}, \ef{3.4-2} and   \ef{3.4-4} to get
\begin{align}
  (1+t) (E_{0I}+E_{0II})(t)+\int_0^t \lt( (1+s) E_{0I} +  E_{0II} \rt)(s) ds \notag\\
\les   (E_{0I}+E_{0II})(0) + (1+t) V_0(t) + \int_0^t   V_0(s) ds . \label{8.20a}
\end{align}
Due to \ef{7.9} and \ef{a}, we have
\begin{align*}
|{\rm curl}_\eta \oa -{\rm curl} \oa | + |{\rm div}_\eta \oa - {\rm div}   \oa | \les |\pl \oa|^2 \les \ea_0 |\pl \oa|,
\end{align*}
which, together with \ef{6.7-1c}, implies that
\begin{align*}
& V_0(t) \les (1+t)^{-1}\mathfrak{V}_0(t) + \ea_0 E_{0II}(t), \\
& \mathfrak{E}_0(t)  \thicksim  (1+t) (E_{0I}+E_{0II})(t),   \\
& \mathfrak{D}_0(t) \thicksim    (1+t) E_{0I}(t) +  E_{0II} (t).
\end{align*}
Substitute these into \ef{8.20a} to obtain \ef{low}. \hfill $\Box$

\subsection{Preliminaries for  the higher order estimates}\label{subsec-1}

\subsubsection{Basic identities} The following identities indicate how the higher order functional are constructed.
\begin{lem}
For any vector field  $F $ with $F^i=F_i$, we have
\begin{align}
& A^k_rA^s_i (\pl_s  F^r)  \pl_t \pl_k F^i  =2^{-1} \pl_t \lt( |\nabla_\eta F|^2 - |{\rm curl}_\eta F|^2 \rt) +    \lt[ \nabla_{\eta } F^r \rt]_i \lt[\nabla_\eta \pl_t \oa^s\rt]_r\lt[\nabla_\eta F^i\rt]_s, \label{nabt}\\
&A^k_r A^s_i ( \pl_s F^r ) \pl_k F^i  = |\nabla_\eta F|^2 - |{\rm curl}_\eta F|^2. \label{nab}
\end{align}
\end{lem}

{\em Proof}.  We commute $\nabla_\eta$ with $\pl_t $ and use $ \pl_t  A^k_i = - A^k_r A^s_i \pl_t\pl_s \oa^r$ to obtain
\begin{align}
\lt[ \nabla_{\eta } F^r \rt]_i   \lt[\nabla_{\eta}\pl_t F^i\rt]_r  =\lt[ \nabla_{\eta } F^r \rt]_i  \pl_t \lt[\nabla_{\eta} F^i\rt]_r - \lt[ \nabla_{\eta } F^r \rt]_i  \lt( \pl_t  A^k_r \rt) \pl_k F^i \notag\\
= \lt[ \nabla_{\eta } F^r \rt]_i  \pl_t \lt[\nabla_{\eta} F^i\rt]_r  +  \lt[ \nabla_{\eta } F^r \rt]_i \lt[\nabla_\eta \pl_t \oa^s\rt]_r\lt[\nabla_\eta F^i\rt]_s. \label{8.20b}
\end{align}
Simple calculation gives that
\begin{align*}
 \lt[ \nabla_{\eta } F^r \rt]_i \pl_t \lt[\nabla_{\eta}   F^i\rt]_r
=& \sum_{i,r} \lt\{  \lt[ \nabla_{\eta } F_r \rt]_i - \lt[ \nabla_{\eta } F_i \rt]_r   \rt\} \pl_t \lt\{ \lt[\nabla_{\eta}  F_i\rt]_r - \lt[\nabla_{\eta}   F_r\rt]_i  \rt\}
\\
 & +2 \sum_{i,r} \lt[ \nabla_{\eta } F_i \rt]_r \pl_t \lt[ \nabla_{\eta } F_i \rt]_r
 -\sum_{i,r} \lt[ \nabla_{\eta } F_i \rt]_r  \pl_t \lt[ \nabla_{\eta } F_r \rt]_i\\
 =& - \pl_t  |{\rm curl}_\eta F|^2 + \pl_t |\nabla_\eta F|^2 -\lt[ \nabla_{\eta } F^i \rt]_r  \pl_t \lt[ \nabla_{\eta } F^r \rt]_i,
\end{align*}
which implies that
$$\lt[ \nabla_{\eta } F^r \rt]_i \pl_t \lt[\nabla_{\eta}   F^i\rt]_r
 = 2^{-1} \lt( \pl_t |\nabla_\eta F|^2- \pl_t  |{\rm curl}_\eta F|^2 \rt).$$
Substitute this into \ef{8.20b} to get  \ef{nabt}.  \ef{nab} can be proved similarly.  \hfill $\Box$

\subsubsection{Commutators}
The following estimates are for commuting  $\pl $ and $\bar\pl $. We use the notation $e_1=(1,0,0)$, $e_2=(0,1,0)$ and $e_3=(0,0,1)$ here.

\begin{lem} For any   function $f$ and multi-indexes  $\alpha$ and $\beta$, we have
\begin{align}
& \big|[\bar\pl^\beta , \pl^\alpha ]f\big| \le C(\alpha, \beta) \sum_{0\le j  \le |\beta|-1} \lt|\pl^{|\alpha|} \bar\pl^j f \rt| . \label{commutator2}
\end{align}
\end{lem}

{\em Proof}. We use the mathematical induction to prove \ef{commutator2}, and  first show that
\begin{align}
&\sum_{|\al|=1}\big|[\bar\pl^\beta , \pl^{\al}]f\big| \le C( \beta) \sum_{0\le j  \le |\beta|-1} \lt|\pl  \bar\pl^j f \rt|   . \label{commutator1}
\end{align}
Clearly, \ef{commutator1} holds for
$|\beta|=1$, due to
$[\bar\pl_i, \pl_l]f= -\epsilon^{ijk}\da_{lj}\pl_k f $.
Suppose that \ef{commutator1} holds for $|\beta|=1, \cdots, L-1$,
and note that for $|\beta|=L$,
\begin{align*}
 \bar\pl^\beta   \pl_l f
  = & \bar\pl^{\beta-e_i} \bar\pl_i  \pl_l f   =\bar\pl^{\beta-e_i} [\bar\pl_i , \pl_l ] f
 + \bar\pl^{\beta-e_i}  \pl_l \bar\pl_i  f \\
  = &  -\sum_{i}\epsilon^{ijk}\da_{lj}\bar\pl^{\beta-e_i}  \pl_k f
   + [\bar\pl^{\beta-e_i},  \pl_l] \bar\pl_i  f + \pl_l \bar\pl^{\beta}  f \\
    = & -\sum_{i}\epsilon^{ijk}\da_{lj} ([\bar\pl^{\beta-e_i},  \pl_k] f + \pl_k\bar\pl^{\beta-e_i} f)  + [\bar\pl^{\beta-e_i},  \pl_l] \bar\pl_i  f + \pl_l \bar\pl^{\beta}  f.
\end{align*}
Then, \ef{commutator1} holds for $|\beta|=L$ using the   induction assumption. Similarly, we apply the mathematical induction to  $\alpha$ and obtain \ef{commutator2}.
\hfill
$\Box$

\begin{lem}  For any   function $f$, and multi-indexes  $\alpha$ and $\beta$, we have for $k=1,2,3$,
\begin{align}
& \lt| \pl^\alpha \bar\pl^\beta \lt(  \sa^{-\iota}\pl_k  ( \sa^{\iota+1} f  ) \rt)- \sa^{-\iota - |\alpha|} \pl_k \lt(\sa^{\iota+ |\alpha|+ 1 }  \pl^\alpha \bar\pl^\beta f   \rt) \rt| \notag\\
\le &   C \sum_{0\le j \le |\beta|-1} \lt(\sa \lt|\pl^{|\al|+1}\bar\pl^j f\rt| +   \lt|\pl^{|\al| }\bar\pl^j f \rt| \rt) +  C |\alpha| \sum_{0\le j \le |\beta|+1}
\lt| \pl^{|\alpha|-1} \bar\pl^j f \rt|,  \label{est1}
\end{align}
where $C=C(\al,\ba,\iota,\Omega)$.
\end{lem}

{\em Proof}. Recall that $\sigma(y)  =\underline{A}-\underline{B}|y|^2$, then we have $\bar\pl \sa=0$ and
 \begin{align}
& \pl^\alpha \bar\pl^\beta  \lt(  \sa^{-\iota}\pl_k  ( \sa^{\iota+1} f  ) \rt) =    \pl^\alpha \bar\pl^\beta \lt( (\iota+1) (\pl_k \sa) f  + \sa \pl_k f  \rt)\notag\\
= &  (\iota+1)   (\pl_k \sa) \pl^\alpha \bar\pl^\beta   f
  +  \sa  \pl_k   \pl^\alpha \bar\pl^\beta  f   +  \sum_{1\le j\le 3}  \alpha_i (\pl_i \sa) \pl_k  \pl^{\alpha-  e_i} \bar\pl^\beta  f   + \sum_{1\le j\le 3} I_{k,j}^{ \al,\ba} , \notag
\end{align}
where
\begin{align*}
&I_{k,1}^{ \al,\ba}=   (\iota+1) \lt(  \pl^\alpha \bar\pl^\beta ( (\pl_k \sa) f  ) -(\pl_k \sa)  \pl^\alpha \bar\pl^\beta   f  \rt), \\
& I_{k,2}^{ \al,\ba}=    \pl^\alpha  \lt( \sa \bar\pl^\beta \pl_k f  \rt)-\sa  \pl^\alpha \bar\pl^\beta  \pl_k f - \sum_{1\le i \le 3} \alpha_i (\pl_i \sa)  \pl^{\alpha-  e_i} \bar\pl^\beta \pl_k f, \\
& I_{k,3}^{ \al,\ba}= \sa \pl^\alpha[ \bar\pl^\beta ,  \pl_k] f   + \sum_{1\le i \le 3} \alpha_i (\pl_i \sa)      \pl^{\alpha-  e_i} [\bar\pl^\beta , \pl_k ]   f.
\end{align*}
This implies that
\begin{align}
  \pl^\alpha \bar\pl^\beta  \lt(  \sa^{-\iota}\pl_k  ( \sa^{\iota+1} f  ) \rt) -   \sa^{-\iota - |\alpha|} \pl_k \lt(\sa^{\iota+ |\alpha|+ 1 } \pl^\alpha \bar\pl^\beta  f    \rt) = \sum_{1\le j\le 4}  I_{k,j}^{ \al,\ba},  \label{4.20}
\end{align}
where
$$ I_{k,4}^{ \al,\ba}= \sum_{1\le i \le 3} \alpha_i (  (\pl_i \sa) \pl_k - (\pl_k \sa) \pl_i  )  \pl^{\alpha-  e_i} \bar\pl^\beta f  .$$
When $|\alpha|\ge 2$, $\pl^\alpha \bar\pl^\ba \pl_k \sa=0$ for any $\beta$ and $k$, and then
\begin{align*}
|I_{k,1}^{ \al,\ba}| \le   C(  \beta,\iota,\Omega) \sum_{0\le j  \le |\beta|-1}
\lt| \pl^{|\alpha| } \bar\pl^j f \rt| +   C(  \beta,\iota,\Omega) |\alpha| \sum_{0\le j \le |\beta|}
\lt| \pl^{|\alpha|-1} \bar\pl^j f \rt| .
\end{align*}
Due to \ef{commutator2}, $\pl_i \sa =- 2B y_i$, $\pl_i \pl_j \sa =- 2B \da_{ij}$ and   $\pl^\alpha \sa =0$ for $|\alpha|\ge 3$, we have
\begin{align*}
&|I_{k,2}^{ \al,\ba}| \le C(\al,\ba, \Omega)  \sum_{1\le i\le 3} \alpha_i(\alpha_i-1)  \sum_{0\le j \le |\beta| } \lt|\pl^{|\alpha|-1}\bar\pl^{j } F \rt|,
\\
&|I_{k,3}^{ \al,\ba}| \le  C(\al,\ba, \Omega)  \sum_{0\le j \le |\beta|-1} \lt(\sa \lt|\pl^{|\al|+1}\bar\pl^j F\rt| +  |\alpha| \lt|\pl^{|\al| }\bar\pl^j F \rt| \rt),\\
& |I_{k,4}^{ \al,\ba}| \le  C(\al,\ba )  |\alpha|  \lt(\big|\pl^{|\alpha|-1}\bar\pl^{|\beta|+1}  F \big| +   \big|\pl^{|\alpha|-1}\bar\pl^{|\beta| } F \big| \rt).
\end{align*}
These estimates, together with \ef{4.20},  prove \ef{est1}.
\hfill $\Box$

\subsubsection{Derivatives of $A$ and $J$}
The differentiation formulae for $A$ and $J$ are
\begin{subequations}\label{7.12}\begin{align}
&\pl_j J=J A^s_r \pl_{j }\pl_{s } \oa^{r},     \ \  \ \  \bar\pl_j J=J A^s_r \bar\pl_j \pl_s  \oa^{r},  \  \  \ \  \pl_t J=J A^s_r \pl_t \pl_s \oa^r,  \label{7.12-1}\\
&\pl_j A^k_i = - A^k_r   A^s_i  \pl_{j }\pl_{s } \oa^r,  \ \  \ \   \bar\pl_j A^k_i = - A^k_r   A^s_i \bar\pl_{j }\pl_{s } \oa^r,  \ \  \ \
  \pl_t  A^k_i = - A^k_r A^s_i \pl_t\pl_s \oa^r,\label{7.12-2}
\end{align}\end{subequations}
which, together with \ef{6.7-1a}, implies that for any polynomial function $\mathscr{P}$,
\begin{subequations}\label{7.12h}\begin{align}
 &|\pl  \mathscr{P}(J) | +   |\pl  \mathscr{P}(A) | \les |\pl^2 \oa |,
  \\
&|\bar\pl  \mathscr{P}(J) | +   |\bar\pl  \mathscr{P}(A) | \les |\bar\pl \pl  \oa |,
\\
&|\pl_t  \mathscr{P}(J) | +  |\pl_t  \mathscr{P}(A) | \les |\pl_t \pl  \oa |. \label{8.23a}
\end{align}\end{subequations}
Moreover, we can use the mathematical induction to obtain that for  any polynomial function $\mathscr{P}$, nonnegative integers $m$, and multi-indexes $\al$ and $\ba$,
\begin{align}
\lt|\pl_t^m \pl^\alpha \bar\pl^\ba \mathscr{P} (A )  \rt| + \lt|\pl_t^m \pl^\alpha \bar\pl^\ba \mathscr{P}(J) \rt|
\les  \mathcal{I}^{m,|\alpha|,|\beta|},  \label{5.29}
\end{align}
where $\mathcal{I}^{m,|\alpha|,|\beta|}$ are defined inductively as follows:
\begin{subequations}\label{5.30}\begin{align}
& \mathcal{I}^{0,0,0}=1, \\
& \mathcal{I}^{m,|\alpha|,|\beta|} =\lt|\pl_t^m \pl^{|\alpha|} \bar\pl^{|\ba|} \pl \oa \rt| +    \sum_ {\substack{0\le i\le m, \ 0\le j\le |\al|, \ 0\le k \le |\ba|
\\ 1\le i+j+k \le m+|\al|+|\ba|-1
  }}  \mathcal{I}^{i,j,k}
\lt|\pl_t^{m-i}\pl^{|\al|-j}\bar\pl^{|\ba|-k}\pl\oa\rt|. \label{5.30b}
 \end{align}\end{subequations}
We use the notation  $\widetilde{\mathcal{I}}^{m,|\alpha|,|\beta|}$ to denote
the lower order terms in $\mathcal{I}^{m,|\alpha|,|\beta|}$, that is,
\begin{align}\label{5.30-1}
\widetilde{\mathcal{I}}^{m,|\alpha|,|\beta|}
= \sum_{\substack{0\le i\le m, \ 0\le j\le |\al|, \ 0\le k \le |\ba|
\\ 1\le i+j+k \le m+|\al|+|\ba|-1
  }}  \mathcal{I}^{i,j,k}
\lt|\pl_t^{m-i}\pl^{|\al|-j}\bar\pl^{|\ba|-k}\pl\oa\rt|.
\end{align}
Then, we have the following estimates.
\begin{lem}\label{lem-non}
For any $m$, $\al$ and $\ba$ satisfying $2\le m+|\al|+|\ba|\le [\iota]+8$, we have
\begin{align}\label{5.30-2}
(1+t)^{ 2m   }  \int  \sigma^{\iota+ |\alpha|+1 } \big| \widetilde{\mathcal{I}}^{m,|\al|,|\ba|}\big|^2      dy \les \mathscr{E}(t) \sum_{0\le j\le m+|\al|+|\ba|-1}\mathscr{E}_j(t),
\end{align}
provided that $\mathscr{E}(t)$ is small.
 \end{lem}

{\em Proof}. It follows from \ef{5.30} and \ef{commutator2} that
\begin{align*}
 (1+t)^{2m}  \int  \sigma^{\iota+ |\alpha|+1 } | \widetilde{\mathcal{I}}^{m,|\al|,|\ba|}|^2      dy
 \les
 \sum_{(i,j,k)\in S} \mathcal{P}_{i,j,k}^{m,|\al|,|\ba|}  +  l.o.t. ,
\end{align*}
where $l.o.t.$ represents the lower order terms, and
\begin{align*}
&S=\lt\{(i,j,k)\in \mathbb{Z}^3  \  \big|  \     0\le i\le m, \ 0\le j\le |\al|, \ 0\le k \le |\ba|, \rt.  \\
& \lt. \qquad \qquad \qquad \qquad \qquad 1\le i+j+k \le m+|\al|+|\ba|-1\rt\},
 \\
&\mathcal{P}_{i,j,k}^{m, |\al|,|\ba|}= (1+t)^{ 2m   } \int  \sigma^{\iota+ |\alpha|+1 }
  \lt|\pl_t^{ i}\pl^{ j+1}\bar\pl^{ k} \oa\rt|^2
  \lt|\pl_t^{m-i}\pl^{|\al|-j+1}\bar\pl^{|\ba|-k} \oa\rt|^2 dy.
\end{align*}
It suffices to prove that
\begin{align}
 \sum_{(i,j,k)\in S} \mathcal{P}_{i,j,k}^{m,|\al|,|\ba|} \les  \mathscr{E}(t)   \sum_{0\le j\le m+|\al|+|\ba|-1}\mathscr{E}_j(t), \label{5.31-4}
\end{align}
since $l.o.t.$ can be bounded similarly. To prove \ef{5.31-4}, it is enough to consider the case of $i+2j+k \le 2^{-1}(m+2|\al|+|\ba|)$, since the other case can be dealt with analogously. In what follows, we assume $(i,j,k)\in S$ and $i+2j+k \le 2^{-1}(m+2|\al|+|\ba|)$, which implies that $i+j+k\le [\iota]+6$.

When $   i+2j+k \le 2$, it follows from \ef{verify} that
\begin{align}
 \mathcal{P}_{i,j,k}^{m, |\al|,|\ba|}
\les & \mathscr{E}(t)   (1+t)^{ 2m-2i }\int \sigma^{\iota+ |\alpha|+1 }
\lt|\pl_t^{m-i}\pl^{|\al|-j+1}\bar\pl^{|\ba|-k} \oa\rt|^2 dy
\notag\\
\les & \mathscr{E}(t) \mathscr{E}_{m+|\al|+|\ba|-i-j-k}(t)
 . \label{5.30-4}
\end{align}

When $ i+2j+k = 3$, it follows from the H\"{o}lder inequality, \ef{verify} and the fact that $H^1(\Omega)\hookrightarrow L^6(\Omega)$ and $H^{1/2}(\Omega)\hookrightarrow L^3(\Omega)$ that
\begin{align*}
 \mathcal{P}_{i,j,k}^{m, |\al|,|\ba|}
\les   (1+t)^{2i} \lt\|\pl_t^{ i}\pl^{ j+1}\bar\pl^{ k} \oa\rt\|_{L^6}^2 (1+t)^{ 2m-2i } \lt\| \sigma^{\frac{\iota+ |\alpha|+1}{2} }
 \pl_t^{m-i}\pl^{|\al|-j+1}\bar\pl^{|\ba|-k} \oa\rt\|_{L^3}^2 \\
 \les \mathscr{E}(t)   (1+t)^{ 2m-2i } \lt\| \sigma^{\frac{\iota+ |\alpha|+1}{2} }
 \pl_t^{m-i}\pl^{|\al|-j+1}\bar\pl^{|\ba|-k} \oa\rt\|_{H^{1/2}}^2.
\end{align*}
Due to \ef{wsv}  and \ef{hard}, we have
\begin{align*}
&  \lt\| \sigma^{\frac{\iota+ |\alpha|+1}{2} }
 \pl_t^{m-i}\pl^{|\al|-j+1}\bar\pl^{|\ba|-k} \oa\rt\|_{H^{1/2}}^2 \les \lt\| \sigma^{\frac{\iota+ |\alpha|+1}{2} }
 \pl_t^{m-i}\pl^{|\al|-j+1}\bar\pl^{|\ba|-k} \oa\rt\|_{H^{1,1}}^2 \\
 \les & \int \sigma^{\iota+ |\alpha|  } \lt(
\lt|\pl_t^{m-i}\pl^{|\al|-j+1}\bar\pl^{|\ba|-k} \oa\rt|^2
+ \lt|\sa\pl_t^{m-i}\pl^{|\al|-j+2}\bar\pl^{|\ba|-k} \oa\rt|^2   \rt)dy \\
\les & \int \sigma^{\iota+ |\alpha| +2 } \lt(
\lt|\pl_t^{m-i}\pl^{|\al|-j+1}\bar\pl^{|\ba|-k} \oa\rt|^2
+ \lt| \pl_t^{m-i}\pl^{|\al|-j+2}\bar\pl^{|\ba|-k} \oa\rt|^2   \rt)dy \\
\les & (1+t)^{2i-2m} \sum_{0\le h\le 1} \mathscr{E}_{m+|\al|+|\ba|+h-(i+j+k)},
\end{align*}
which, together with $2\le i+j+k \le 3$, implies that
\begin{align}
 \mathcal{P}_{i,j,k}^{m, |\al|,|\ba|}
 \les  \mathscr{E}(t) \sum_{0\le h\le 1} \mathscr{E}_{m+|\al|+|\ba|-i-j-k+h}(t)\le  \mathscr{E}(t)   \sum_{0\le j\le m+|\al|+|\ba|-1}\mathscr{E}_j(t) . \label{5.30-6}
\end{align}

When $ i+2j+k \ge 4$, it follows from \ef{verify} that
\begin{align}
  \mathcal{P}_{i,j,k}^{m, |\al|,|\ba|}
\les   \mathscr{E}(t)  (1+t)^{ 2m-2i }\int \sigma^{\iota+ |\alpha|+3-i-2j-k }
\lt|\pl_t^{m-i}\pl^{|\al|-j+1}\bar\pl^{|\ba|-k} \oa\rt|^2 dy . \label{5.30-8}
\end{align}
To apply \ef{hard} to the right hand side of \ef{5.30-8}, we need $\iota+ |\alpha|+3-i-2j-k >-1$, which is the case for $m+|\al|+|\ba|\le [\iota]+8$, due to
\begin{align*}
 \iota+ |\alpha|+3-(i+2j+k)
 \ge   \iota+ |\alpha|+3 - 2^{-1}(m+2|\al|+|\ba|) \\
 =  \iota +3 - 2^{-1}(m +|\ba|)
 \ge 2^{-1} \iota - 1 >-1.
\end{align*}
So, we can apply $i+j+k-2$ times of \ef{hard} to the right hand side of \ef{5.30-8}  to get
\begin{align}
 \mathcal{P}_{i,j,k}^{m, |\al|,|\ba|}
\les & \mathscr{E}(t)  (1+t)^{ 2m-2i }\int \sigma^{\iota+ |\alpha| +i+k-1 }
\sum_{0\le h\le i+j+k-2}\lt|\pl_t^{m-i}\pl^{|\al|-j+1+h}\bar\pl^{|\ba|-k} \oa\rt|^2 dy \notag \\
\les & \mathscr{E}(t)   \sum_{0\le j\le m+|\al|+|\ba|-2}\mathscr{E}_j(t).  \label{5.30-9}
\end{align}

Finally, \ef{5.31-4} is a consequence of \ef{5.30-4}, \ef{5.30-6} and \ef{5.30-9}. This finishes the proof of \ef{5.30-2}. \hfill $\Box$

In addition to the estimates stated in Lemma \ref{lem-non}, we  need the following estimate to perform the curl estimates in Section \ref{sec-curl}.

\begin{lem}\label{lem-non-new}
For any $m$, $\al$ and $\ba$ satisfying $2\le m+|\al|+|\ba|\le [\iota]+9$, we have
\begin{align}
\sum_{\substack{0\le i\le m, \ 0\le j\le |\al|, \ 0\le k \le |\ba|
\\  i+j+k \le m+|\al|+|\ba|-1, \  4\le i+2j+k \le m+2|\al|+|\ba|-4
  }}    \int  \sigma^{\iota+ |\alpha|+1 } \big|\mathcal{I}^{i,j,k}
 \pl_t^{m-i}\pl^{|\al|-j}\bar\pl^{|\ba|-k}\pl\oa  \big|^2      dy \notag \\
\les (1+t)^{- 2m   } \mathscr{E}(t) \sum_{0\le j\le m+|\al|+|\ba|-2}\mathscr{E}_j(t), \label{7.14'}
\end{align}
provided that $\mathscr{E}(t)$ is small.
\end{lem}

{\em Proof}. In the spirit of the proof of  \ef{5.30-2},  it suffices to prove that
\be\label{8.20c}
 \mathcal{P}_{i,j,k}^{m, |\al|,|\ba|}
\les  \mathscr{E}(t)   \sum_{0\le j\le m+|\al|+|\ba|-2}\mathscr{E}_j(t)
\ee
for $m+|\al|+|\ba| \le [\iota]+9$ and  $4\le  i+2j+k \le  2^{-1}(m+2|\al|+|\ba|) $.
In a similar way to the derivation of \ef{5.30-9}, we can show that \ef{8.20c} holds except for a bad case of $\iota=1$, $|\al|=j=0$ and $i+k=2^{-1}(m+|\ba|)=5$, where
\begin{align*}
\mathcal{P}_{i,j,k}^{m, |\al|,|\ba|}= & (1+t)^{ 2m   }  \lt\|  \sigma
 |\pl_t^{ i}\pl \bar\pl^{ k} \oa |
  |\pl_t^{m-i}\pl \bar\pl^{|\ba|-k} \oa | \rt\|_{L^2}^2 \\
  \le &  (1+t)^{ 2i   }  \lt\|  \sigma
 \pl_t^{ i}\pl \bar\pl^{ k} \oa   \rt\|_{L^6}^2 (1+t)^{ 2m-2i  }  \lt\|
  \pl_t^{m-i}\pl \bar\pl^{|\ba|-k} \oa   \rt\|_{L^3}^2,
\end{align*}
due to the H\"{o}lder inequality.
It follows from \ef{wsv} and \ef{hard} that
\begin{align*}
& \|\sa \pl_t^i \pl \bar\pl^k \oa\|_{H^1}^2 \les  \int \lt( |  \pl_t^i \pl \bar\pl^k \oa|^2 +   \sa^2  |\pl \pl_t^i  \pl \bar\pl^k \oa|^2   \rt)dy\\
\les & \sum_{1\le l \le 3} \int     \sa^4  | \pl_t^i  \pl^l \bar\pl^k \oa|^2  dy
\le (1+t)^{-2i} \sum_{1\le l \le 3} \mathscr{E}_{4+l}(t) \le (1+t)^{-2i}   \mathscr{E} (t)
\end{align*}
and
\begin{align*}
& \| \pl_t^{m-i} \pl \bar\pl^{|\ba|-k} \oa\|_{H^{1/2}}^2 \les \| \pl_t^{m-i} \pl \bar\pl^{|\ba|-k} \oa\|_{H^{1, 1}}^2\\
 = &
 \int \sa ( |\pl_t^{m-i} \pl \bar\pl^{|\ba|-k} \oa|^2 +|\pl \pl_t^{m-i} \pl \bar\pl^{|\ba|-k} \oa|^2   )dy \\
 \les & \sum_{1\le l\le 4} \int \sa^5  |\pl_t^{m-i} \pl^l \bar\pl^{|\ba|-k} \oa|^2  dy \les (1+t)^{2i-2m} \sum_{4\le j\le 8}\mathscr{E}_j(t) .
\end{align*}
This, together with the fact that $H^1(\Omega)\hookrightarrow L^6(\Omega)$ and $H^{1/2}(\Omega)\hookrightarrow L^3(\Omega)$, implies that  \ef{8.20c} holds for the bad case.
\hfill $\Box$

\subsection{The higher order estimates}
In this subsection, we prove that

\begin{lem}\label{est-high} Let $\oa(t,y)$ be a solution to problem \ef{newsystem} in the time interval $[0, T]$ satisfying \ef{assume}. Then
for $j=1,2,\cdots, [\iota]+7 $,
\begin{align}\label{hig}
   \mathfrak{E}_{j}(t) + \int_0^t \mathfrak{D}_{j}(s)ds \les  \sum_{0\le k \le j }\lt( \mathfrak{E}_k  (0)+   \mathfrak{V}_k  (t)  + \int_0^t(1+s)^{-1} \mathfrak{V}_k (s)ds \rt), \ \ t\in [0,T].
 \end{align}
\end{lem}

{\em Proof}.
Apply $\pl_t^m\pl^\alpha \bar\pl^\beta$ to the product of  $\ta^{3\ga-1}\sa^{-\iota}$ and  \ef{3.4}, and multiply the resulting equation by  $\ta^{1-3\ga}$ to  obtain
\begin{align}
&    \pl_t^{m+2} \pl^\alpha \bar\pl^\beta\oa_i +   \lt(1 +\lt(2+ m(3\ga-1) \rt) \ta^{-1}\ta_t\rt)     \pl_t^{m+1} \pl^\alpha \bar\pl^\beta \oa_i
 + \lt(({ 3\ga-1 })^{-1} \ta^{1-3\ga} \rt. \notag \\
& \lt. + m(3\ga-1)  \ta^{-1}\ta_t   \rt)\pl_t^{m} \pl^\alpha \bar\pl^\beta  \oa_i
 +\ta^{1-3\ga} { \sigma^{ - \iota - |\alpha|}}  \pl_k \lt( \sigma^{\iota+ |\alpha|+1} ( \mathcal{R}_{1, i}^{m,\al,\ba,k} \rt.
 \notag \\
&   \lt. -  {J}^{1-\ga}  (   A^k_r  A^s_i  \pl_s \pl_t^{m} \pl^\alpha \bar\pl^\beta   \omega^r
 + \iota^{-1} A_i^k {\rm div}_\eta   \pl_t^{m} \pl^\alpha \bar\pl^\beta   \omega    )) \rt)
  = \ta^{1-3\ga} \sum_{j=2,3}  \mathcal{R}_{j,i}^{m,\al,\ba}, \label{5.1}
\end{align}
where
\begin{align*}
&\mathcal{R}_{1,i}^{m,\al,\ba,k} =       \pl_t^m\pl^\alpha\bar\pl^\beta ( A_i^k {J}^{1-\ga}  -  \da_i^k ) \\
&\qquad \qquad \quad  + {J}^{1-\ga} (   A^k_r  A^s_i  \pl_s \pl_t^{m} \pl^\alpha \bar\pl^\beta   \omega^r   + \iota^{-1} A_i^k {\rm div}_\eta   \pl_t^{m} \pl^\alpha \bar\pl^\beta   \omega    ),
 \\
& \mathcal{R}_{2,i}^{m,\al,\ba}=   { \sigma^{ - \iota - |\alpha|}} \pl_k \lt( \sigma^{\iota+ |\alpha|+1}   \pl_t^m\pl^\alpha\bar\pl^\beta ( A_i^k {J}^{1-\ga}  -  \da_i^k  ) \rt)
\\
  & \qquad \qquad \quad  -\pl_t^m\pl^\alpha \bar\pl^\beta \lt( \sigma^{-\iota } \pl_k  \lt( \sigma^{\iota+1}  (A_i^k {J}^{1-\ga} -  \da_i^k ) \rt) \rt),
  \\
  & \mathcal{R}_{3,i}^{m,\al,\ba}= -\sum_{2\le k \le m} C_m^k  (\pl_t^k   \ta^{3\ga-1})  \lt( \pl_t^{m-k+2} \pl^\alpha \bar\pl^\beta\oa_i+ \pl_t^{m-k+1} \pl^\alpha \bar\pl^\beta\oa_i\rt)
 \\
&\qquad \qquad\quad  -2 ({ 3\ga-1 })^{-1} \sum_{1\le k\le m} C_m^k  (\pl_t^{k+1}   \ta^{3\ga-1} )    \pl_t^{m-k+1} \pl^\alpha \bar\pl^\beta\oa_i.
\end{align*}
It should be noted that the terms on the right hand side of \ef{5.1} are not principal ones. So, we will first analyze the principal terms on the left hand side and then do the others.

{\em Step 1}. In this step, we will focus on  the left hand side of equation \ef{5.1} and show where the functionals $\mathfrak{E}_{j}$ and $\mathfrak{D}_{j}$ come from. We integrate the product of $\sigma^{  \iota + |\alpha|}\pl_t^{m+1} \pl^\alpha \bar\pl^\beta\oa^i$ and \ef{5.1}   over $\Omega$  and use \ef{nabt}  to  get
\begin{align}\label{5.2}
\frac{d}{dt} \mathcal{E}_1^{m,\al,\ba}(t)+  \mathcal{D}_1^{m,\al,\ba}(t)  =  \mathcal{H}^{m,\al,\ba}(t) =\sum_{1\le j \le 4} \mathcal{H}_j^{m,\al,\ba}(t),
\end{align}
where
\begin{align*}
&  \mathcal{E}_1^{m,\al,\ba}(t)   =    \frac{1}{2} \int \sigma^{  \iota + |\alpha|} \lt| \pl_t^{m+1} \pl^\alpha \bar\pl^\beta\oa \rt|^2 dy \\
& \quad    +  \frac{1}{2}\lt(({ 3\ga-1 })^{-1} \ta^{1-3\ga} + m(3\ga-1)  \ta^{-1}\ta_t   \rt)
 \int  \sigma^{  \iota + |\alpha|}\lt| \pl_t^{m} \pl^\alpha \bar\pl^\beta\oa \rt|^2 dy
\\
&\quad +\frac{1}{2}\ta^{1-3\ga}  \int  \sigma^{\iota+ |\alpha|+1} {J}^{1-\ga}  (   |\na_\eta \pl_t^m\pl^\alpha \bar\pl^\beta \oa|^2  - |{\rm curl}_\eta \pl_t^m\pl^\alpha \bar\pl^\beta \oa|^2
\\
    &\quad  + \iota^{-1}    |{\rm div}_\eta   \pl_t^m\pl^\alpha \bar\pl^\beta \oa|^2      ) dy
 -   \ta^{1-3\ga} \int   \sigma^{\iota+ |\alpha|+1} \mathcal{R}_{1,i}^{m,\al,\ba,k}  \pl_t^m \pl_k \pl^\alpha \bar\pl^\beta \oa^i dy,
\\
&  \mathcal{D}_1^{m,\al,\ba} (t)  =   \lt(1 +\lt(2+ m(3\ga-1) \rt) \ta^{-1}\ta_t\rt)    \int  \sigma^{  \iota + |\alpha|} \lt|  \pl_t^{m+1} \pl^\alpha \bar\pl^\beta \oa \rt|^2 dy  ,
\\
& \mathcal{H}_1^{m,\al,\ba}(t)   =   - \int  \sigma^{\iota+ |\alpha|+1}  \pl_t ( \ta^{1-3\ga} \mathcal{R}_{1,i}^{m,\al,\ba,k} )  \pl_t^m \pl_k \pl^\alpha \bar\pl^\beta \oa^i dy,
\\
& \mathcal{H}_j^{m,\al,\ba}(t)  =  \ta^{1-3\ga} \int  \sigma^{  \iota + |\alpha|}  \mathcal{R}_{j,i}^{m,\al,\ba}  \pl_t^{m+1} \pl^\alpha \bar\pl^\beta\oa^i dy , \ \ \ \  j=2,3,
\\
& \mathcal{H}_4^{m,\al,\ba}(t) =   \frac{1}{2} \lt(({ 3\ga-1 })^{-1} \ta^{1-3\ga} + m(3\ga-1)  \ta^{-1}\ta_t   \rt)_t \int  \sigma^{  \iota + |\alpha|}\lt| \pl_t^{m} \pl^\alpha \bar\pl^\beta\oa \rt|^2 dy
\\
&\quad+ \frac{1}{2}    \int  \sigma^{\iota+ |\alpha|+1}\pl_t \lt( \ta^{1-3\ga}{J}^{1-\ga} \rt)  (   |\na_\eta \pl_t^m\pl^\alpha \bar\pl^\beta \oa|^2 - |{\rm curl}_\eta \pl_t^m\pl^\alpha \bar\pl^\beta \oa|^2
\\
& \quad+ \iota^{-1}    |{\rm div}_\eta   \pl_t^m\pl^\alpha \bar\pl^\beta \oa|^2     ) dy
  - \ta^{1-3\ga} \int  \sigma^{\iota+ |\alpha|+1}  {J}^{1-\ga}  (
[\na_\eta \pl_t^{m } \pl^\alpha \bar\pl^\beta\oa^r]_i[\na_\eta \pl_t \oa^s]_r\\
&\quad \times [\na_\eta \pl_t^{m } \pl^\alpha \bar\pl^\beta\oa^i]_s
 -\iota^{-1} (\pl_t A_i^k) ( {\rm div}_\eta   \pl_t^{m} \pl^\alpha \bar\pl^\beta   \omega ) \pl_k \pl_t^{m } \pl^\alpha \bar\pl^\beta\oa^i ) dy.
\end{align*}
It follows from the Cauchy inequality, \ef{decay} (especially, $\ta_t\ge 0$), \ef{6.7-1a} and \ef{6.7-1c} that
\begin{subequations}\label{6.9-1}\begin{align}
&\mathfrak{E}^{m,n,l}(t)- C \mathcal{L}^{m,n,l}(t)  \les
(1+t)^{2m+1} \sum_{|\al|=n, \ |\ba|=l} \mathcal{E}_1^{m,\al,\ba}(t) \les
\mathfrak{E}^{m,n,l}(t) + C \mathcal{L}^{m,n,l}(t),
\label{5.4-1}
\\ & \mathfrak{E}_{I}^{m,n,l}(t) \le  (1+t)^{2m+1} \sum_{|\al|=n, \ |\ba|=l} \mathcal{D}_1^{m,\al,\ba} (t)  \label{5.4-2},
\end{align}\end{subequations}
where
$$\mathcal{L}^{m,n,l} = \mathfrak{V}^{m,n,l}(t) + (1+t) \sum_{|\al|=n, \ |\ba|=l} \mathcal{Q}_1^{m,\al,\ba}(t). $$
Here $\mathcal{Q}_1^{m,\al,\ba}$ defined in \ef{8.22a} is a lower order term shown later in \ef{5.31-1}.
\ef{5.4-1} implies that the  bound of $\mathfrak{E}^{m,n,l}$   can be achieved by integrating the product of  \ef{5.2}  and $(1+t)^{2m+1}$ over $[0,t]$, which needs the bound of
$\int_0^t (1+s)^{2m} \mathcal{E}_1^{m,\al,\ba}  ds$ whose principal part is $\int_0^t (1+s)^{-1} \mathfrak{E}^{m,n,l} ds$. Due to \ef{5.4-2}, the problem turns to estimating $\int_0^t (1+s)^{-1} \mathfrak{E}_{II}^{m,n,l} ds$.

For this purpose, we integrate the product of $ \sigma^{  \iota + |\alpha|}\pl_t^{m} \pl^\alpha \bar\pl^\beta\oa^i$ and \ef{5.1} over $\Omega$  and  use \ef{nab} to give
\bee\label{5.3}
 \frac{d}{dt} \mathcal{E}_2^{m,\al,\ba}(t)+  \mathcal{D}_2^{m,\al,\ba}(t)  =  \mathcal{F}^{m,\al,\ba}(t)  =  \sum_{1\le j \le 4} \mathcal{F}_j^{m,\al,\ba}(t),
\eee
where
\begin{align*}
& \mathcal{E}_2^{m,\al,\ba}(t)  =    \int   \sigma^{  \iota + |\alpha|} ( \pl_t^{m} \pl^\alpha \bar\pl^\beta\oa^i )  \pl_t^{m+1} \pl^\alpha \bar\pl^\beta\oa_i dy \\
 &\quad + \frac{1}{2} \lt(1 +\lt(2+ m(3\ga-1) \rt) \ta^{-1}\ta_t\rt)  \int   \sigma^{  \iota + |\alpha|} | \pl_t^{m} \pl^\alpha \bar\pl^\beta\oa|^2 dy,\\
& \mathcal{D}_2^{m,\al,\ba}(t)  =     \lt(({ 3\ga-1 })^{-1} \ta^{1-3\ga} + m(3\ga-1)  \ta^{-1}\ta_t   \rt) \int  \sigma^{  \iota + |\alpha|} \lt|  \pl_t^{m} \pl^\alpha \bar\pl^\beta \oa \rt|^2 dy \\
 & \quad+   \ta^{1-3\ga}  \int  \sigma^{\iota+ |\alpha|+1} {J}^{1-\ga}  (   |\na_\eta \pl_t^m\pl^\alpha \bar\pl^\beta \oa|^2
    + \iota^{-1}    |{\rm div}_\eta   \pl_t^m\pl^\alpha \bar\pl^\beta \oa|^2      ) dy \\
    & \quad -  \int  \sigma^{  \iota + |\alpha|} \lt|  \pl_t^{m+1} \pl^\alpha \bar\pl^\beta \oa \rt|^2 dy,
  \\
&\mathcal{F}_1^{m,\al,\ba}(t)  =    \ta^{1-3\ga }  \int  \sigma^{\iota+ |\alpha|+1}   \mathcal{R}_{1,i}^{m,\al,\ba,k}    \pl_t^m \pl_k \pl^\alpha \bar\pl^\beta \oa^i dy, \\
& \mathcal{F}_j^{m,\al,\ba}(t)  =   \ta^{1-3\ga} \int  \sigma^{  \iota + |\alpha|}  \mathcal{R}_{j,i}^{m,\al,\ba}  \pl_t^{m} \pl^\alpha \bar\pl^\beta\oa^i dy , \ \  \ \ j=2,3, \\
& \mathcal{F}_4^{m,\al,\ba}(t) =   \frac{1}{2}  \lt(2+ m(3\ga-1) \rt) \lt( \ta^{-1}\ta_t\rt)_t  \int   \sigma^{  \iota + |\alpha|} | \pl_t^{m} \pl^\alpha \bar\pl^\beta\oa|^2 dy \\
 &\quad + \ta^{1-3\ga}  \int  \sigma^{\iota+ |\alpha|+1} {J}^{1-\ga}    |{\rm curl}_\eta \pl_t^m\pl^\alpha \bar\pl^\beta \oa|^2 dy.
\end{align*}
Clearly, it follows from \ef{decay} (especially, $\ta_t\ge 0$), \ef{6.7-1a} and \ef{5.4-2} that
\begin{align}
(1+t)^{-1} \mathfrak{E}_{II}^{m,n,l}(t)
\les & (1+t)^{2m} \sum_{|\al|=n, \  |\ba|=l}    \mathcal{D}_2^{m,\al,\ba}(t) + (1+t)^{-1} \mathfrak{E}_{I}^{m,n,l}(t) \notag \\
\les & (1+t)^{2m} \sum_{|\al|=n, \ |\ba|=l, \ j=1,2}    \mathcal{D}_j^{m,\al,\ba}(t) ,  \label{5.4-3}
\end{align}
which implies
\begin{align}
 \mathfrak{D}^{m,n,l}(t) & \les   \sum_{|\al|=n, \  |\ba|=l} \lt\{  (1+t)^{2m+1} \mathcal{D}_1^{m,\al,\ba} (t) +(1+t)^{2m} ( 4\mathcal{D}_1^{m,\al,\ba} + \mathcal{D}_2^{m,\al,\ba}     )(t)  \rt\}. \label{6.8-1}
 \end{align}
Note that
\begin{subequations}\label{5.4-4}\begin{align}
&   (1+t)^{2m} \sum_{|\al|=n, \  |\ba|=l} (4\mathcal{E}_1^{m,\al,\ba}+ \mathcal{E}_2^{m,\al,\ba} )(t)
   \gtrsim    (1+t)^{-1} \mathfrak{E}^{m,n,l}(t)
\notag \\
  & \quad + (1+t)^{2m} \sum_{|\al|=n, \  |\ba|=l} \int  \sigma^{  \iota + |\alpha|} \lt|  \pl_t^{m} \pl^\alpha \bar\pl^\beta \oa \rt|^2 dy- C (1+t)^{-1} \mathcal{L}^{m,n,l}, \label{5.4-4a}
\\
&  (1+t)^{2m} \sum_{|\al|=n, \  |\ba|=l} (4\mathcal{E}_1^{m,\al,\ba}+ \mathcal{E}_2^{m,\al,\ba} )(t) \les
(1+t)^{-1} \mathfrak{E}^{m,n,l}(t) \notag \\
&\quad  + (1+t)^{2m} \sum_{|\al|=n, \  |\ba|=l} \int  \sigma^{  \iota + |\alpha|} \lt|  \pl_t^{m} \pl^\alpha \bar\pl^\beta \oa \rt|^2 dy + C (1+t)^{-1}\mathcal{L}^{m,n,l}(t),  \label{5.4-4b}
 \end{align}\end{subequations}
due to the Cauchy inequality, \ef{decay}, \ef{6.7-1a} and \ef{6.7-1c}. So, we multiply  the following equation
\begin{align}\label{5.4}
\frac{d}{dt} (4\mathcal{E}_1^{m,\al,\ba} + \mathcal{E}_2^{m,\al,\ba})(t)+  (4\mathcal{D}_1^{m,\al,\ba}+ \mathcal{D}_2^{m,\al,\ba})(t)    =   4\mathcal{H}^{m,\al,\ba}(t)  + \mathcal{F}^{m,\al,\ba} (t)
\end{align}
by $(1+t)^k$ and integrate the resulting  equation over $[0,t]$ from $k=0$ to $k=2m$ step by step, and then  integrate the product of $(1+t)^{2m+1}$ and \ef{5.2} over $[0,t]$ to get the  desired higher order estimates \ef{hig} for $1\le j\le [\iota]+7$.

During the process, it occurs some difficulties in dealing with the first term on the second line of \ef{5.4-4b} in the case of  $m\ge 1$.   For example, in the step $k=2m$, the dissipation we could expect is $(1+t)^{-1} \mathfrak{E}^{m,n,l}(t)$, due to \ef{5.4-2} and \ef{5.4-3}, which should be bounded by  $  (1+t)^{2m-1} (  4\mathcal{E}_1^{m,\al,\ba} + \mathcal{E}_2^{m,\al,\ba} )(t)  $, due to \ef{5.4}, whose principal part contains
$ (1+t)^{2m-1}  \int   \sigma^{  \iota + |\alpha|} \lt|  \pl_t^{m} \pl^\alpha \bar\pl^\beta \oa \rt|^2 $, due to \ef{5.4-4a},  which is a part of  $(1+t)^{-1} \mathfrak{E}_{II}^{m,|\al|,|\ba|}  $, such that nothing could be obtained.  To overcome the difficulty,
we may regard $\int  \sigma^{  \iota + |\alpha|} \lt|  \pl_t^{m} \pl^\alpha \bar\pl^\beta \oa \rt|^2 dy$ as  $\int  \sigma^{  \iota + |\alpha|} \lt| \pl_t \pl_t^{m-1} \pl^\alpha \bar\pl^\beta \oa \rt|^2 dy$, since the latter one can be bounded by $(1+t)^{1-2m}\mathfrak{E}_{I}^{m-1,|\al|,|\ba|}$  which is a lower order term. The technique will be frequently used in dealing with the reminder terms $\mathcal{H}^{m,\al,\ba}$ and  $\mathcal{F}^{m,\al,\ba}$, see \ef{5.3-1} for instance.

{\em Step 2}. In this step, we prove that for any $\vea\in (0,1)$,
\begin{subequations}\label{4-1}\begin{align}
& (1+t)^{  2m + 1 } \mathcal{H}^{m,\al,\ba}(t)  \les \mathcal{H}_g^{m,\al,\ba}(t)  +  \mathcal{H}_b^{m,\al,\ba}(t) ,  \label{4.2-1} \\
& (1+t)^{  2m  } \mathcal{F}^{m,\al,\ba} (t)  \les  \mathcal{F}_g^{m,\al,\ba}(t)  +  \mathcal{F}_b^{m,\al,\ba} (t)   ,\label{4.2-2}
\end{align}\end{subequations}
where
\begin{subequations}\label{8.23c}\begin{align}
& \mathcal{H}_g^{m,\al,\ba}=
  (\vea + \ea_0 +  \vea^{-1} \ea_0^2) \mathfrak{D}_{m+|\al|+|\ba|}
 +   \vea^{-1}   \sum_{0\le j \le m+|\al|+|\ba|-1} \mathfrak{D}_j+ (1+t)^{-1} \mathfrak{V}^{m, |\al|,|\ba|},  \\
& \mathcal{H}_b^{m,\al,\ba} =  \begin{cases}
&  \vea^{-1} |\ba|\lt(\mathfrak{D}^{m+1,|\al|,|\ba|-1}  + \mathfrak{D}^{m,|\al|+1, |\ba|-1} \rt) +  \vea^{-1} |\al|  \mathfrak{D}^{m ,|\al|-1,|\ba|+1},   \\
& \vea^{-1} \lt(|\ba|\mathfrak{D}^{m+1,|\al|,|\ba|-1}  +    \mathfrak{D}^{m+1, |\al|-1,|\ba|} \rt)
  , \ \  \ \   |\al| \ge 1,
\end{cases}  \\
& \mathcal{F}_g^{m,\al,\ba}=   \mathcal{H}_g^{m,\al,\ba}  + (1+t)^{2m-2} \int   \sigma^{  \iota + |\alpha|} | \pl_t^{m} \pl^\alpha \bar\pl^\beta\oa|^2 dy  , \\
& \mathcal{F}_b^{m,\al,\ba} =  \begin{cases}
&  \vea^{-1} |\ba| \mathfrak{D}^{m,|\al|+1, |\ba|-1}  +  \vea^{-1} |\al|  \mathfrak{D}^{m ,|\al|-1,|\ba|+1},   \\
& 0
     , \ \  \ \   |\al| \ge 1.
\end{cases}
\end{align}\end{subequations}
It should be noted that $\mathcal{H}_g^{m,\al,\ba}$ and $\mathcal{F}_g^{m,\al,\ba}$ represent the good terms which can be dealt with easily, in particular, the second term of $\mathcal{F}_g^{m,\al,\ba}$ can be bounded using the Grownwall inequality.
However, we have to use different methods to deal with tangential derivatives ($|\al|=0$) and normal derivatives  ($|\al|\ge 1$), see for instance in $\mathcal{H}_b^{m,\al,\ba}$ and $\mathcal{F}_b^{m,\al,\ba}$. (Indeed, the difference comes from estimates \ef{4.15b} and \ef{5.3-1} which devote to controlling $\mathcal{H}_2^{m,\al,\ba}$.) An example will be given in the next step to illustrate why we have to distinguish these two cases.

First, we prove \ef{4.2-1}.
It follows from \ef{decay} that
\begin{align*}
&(1+t)^{2m+1} \int  \sigma^{  \iota + |\alpha|}  |\mathcal{R}_{3,i }^{m,\al,\ba}|^2  dy \\
\les &  \sum_{1\le k \le m }   (1+t)^{2m-2k+1} \int  \sigma^{  \iota + |\alpha|}   |\pl_t^{m-k+1} \pl^\alpha \bar\pl^\beta\oa |^2 dy\\
&   + \sum_{2\le k \le m } (1+t)^{2m-2k+3} \int  \sigma^{  \iota + |\alpha|}  \lt( | \pl_t^{m-k+2} \pl^\alpha \bar\pl^\beta\oa|^2 + |\pl_t^{m-k+1} \pl^\alpha \bar\pl^\beta\oa|^2 \rt)dy \\
   \les  &   \sum_{1\le k \le m }  \mathfrak{D}^{{m-k} , |\al|, |\ba| } +  \sum_{2\le k \le m }\lt(\mathfrak{D}^{{m-k+1} , |\al|, |\ba| } + (1+t)^2 \mathfrak{D}^{{m-k} , |\al|, |\ba| } \rt)   ,
\end{align*}
which, together with the  Cauchy inequality, implies that for any $\vea\in (0,1)$,
\begin{align}
\mathcal{H}_3^{m,\al,\ba}(t) \les & \vea \int  \sigma^{  \iota + |\alpha|}    |\pl_t^{m+1} \pl^\alpha \bar\pl^\beta\oa|^2  dy +   \vea^{-1} (1+t)^{-2} \sum_i\int  \sigma^{  \iota + |\alpha|}  |\mathcal{R}_{3,i }^{m,\al,\ba}|^2  dy   \notag \\
\les &  (1+t)^{-2m-1}   \lt(\vea \mathfrak{D}^{m , |\al|, |\ba|  } + \vea^{-1}  \sum_{1\le k \le m}   \mathfrak{D}^{{m-k} , |\al|, |\ba| }
  \rt) .\label{4.15c}
\end{align}

It follows from  \ef{decay} (especially, $\ta_t\ge 0$), \ef{6.7-1a}, \ef{6.7-1c}, \ef{8.23a}, and  $|\pl_t \pl \oa| \les \ea_0$ (which is due to \ef{verify} and \ef{assume})   that for any $\vea\in (0,1)$,
\begin{align}
   \mathcal{H}_4^{m,\al,\ba}(t) \les &   (1+t)^{-2} \int  \sigma^{  \iota + |\alpha|} \lt(m\lt| \pl_t^{m} \pl^\alpha \bar\pl^\beta\oa \rt|^2 + \sa |{\rm curl}_\eta \pl_t^m\pl^\alpha \bar\pl^\beta \oa|^2  \rt)dy \notag  \\
 & +   (1+t)^{-1} \|\pl_t \pl \oa\|_{L^\iy} \int   \sigma^{\iota+ |\alpha|+1}  | \pl \pl_t^m\pl^\alpha \bar\pl^\beta \oa|^2 dy \notag\\
\les & (1+t)^{-2m-1}  ( m \mathfrak{D}^{m-1,|\al|,|\ba|} +   \ea_0  \mathfrak{D}^{m ,|\al|,|\ba|})   +(1+t)^{-2m-2} \mathfrak{V}^{m, |\al|,|\ba|}.  \label{4.15d}
\end{align}

It follows from the  Cauchy inequality, \ef{decay}  and   \ef{6.7-1c}   that for any $\vea\in (0,1)$,
\begin{subequations}\label{8.23b}\begin{align}
\mathcal{H}_1^{m,\al,\ba}(t) \les &  \vea (1+t)^{-2} \int  \sigma^{\iota+ |\alpha|+1} |\pl  \pl_t^m  \pl^\alpha \bar\pl^\beta \oa|^2 dy \notag\\
&+ \vea^{-1} (1+t)^2 \sum_{i,k}  \int  \sigma^{\iota+ |\alpha|+1}  |\pl_t ( \ta^{1-3\ga} \mathcal{R}_{1,i}^{m,\al,\ba,k} ) |^2  dy \notag\\
\les &  (1+t)^{-2m -1 }  \lt(\vea   \mathfrak{D}^{m , |\al|, |\ba|  }  + \vea^{-1} \mathcal{Q}_1^{m,\al,\ba} + \vea^{-1}  \mathcal{Q}_2^{m,\al,\ba} \rt),\label{4.15a}\\
 \mathcal{H}_2^{m,\al,\ba}(t) \les &  \vea  \int  \sigma^{\iota+ |\alpha| } | \pl_t^{m+1}  \pl^\alpha \bar\pl^\beta \oa|^2 dy
 + \vea^{-1} (1+t)^{-2} \sum_{i }  \int  \sigma^{\iota+ |\alpha| }  |  \mathcal{R}_{2,i}^{m,\al,\ba }   |^2  dy \notag\\
\les & (1+t)^{-2m -1 }  \lt( \vea    \mathfrak{D}^{m , |\al|, |\ba|  }
+ \vea^{-1}   \mathcal{Q}_3^{m,\al,\ba} \rt),  \label{4.15b}\\
\mathcal{H}_2^{m,\al,\ba}(t) \les & \vea^{-1}  \int  \sigma^{\iota+ |\alpha| } | \pl_t^{m+1}  \pl^\alpha \bar\pl^\beta \oa|^2 dy
+ \vea  (1+t)^{-2} \sum_{i }  \int  \sigma^{\iota+ |\alpha| }  |  \mathcal{R}_{2,i}^{m,\al,\ba }   |^2  dy  \notag \\
\les & (1+t)^{-2m -1 }  \lt( \vea^{-1}   \mathfrak{D}^{m+1 , |\al|-1, |\ba|  }
+ \vea  \mathcal{Q}_3^{m,\al,\ba}  \rt),  \ \   |\al| \ge 1, \label{5.3-1}
\end{align}\end{subequations}
where
\begin{align}
  \mathcal{Q}_1^{m,\al,\ba}= & (1+t)^{ 2m -1 }\sum_{i,k} \int  \sigma^{\iota+ |\alpha|+1 }    | \mathcal{R}_{1,i}^{m,\al,\ba,k}  |^2    dy, \label{8.22a}\\
  \mathcal{Q}_2^{m,\al,\ba}=& (1+t)^{ 2m+1 } \sum_{i,k}  \int  \sigma^{\iota+ |\alpha|+1}   |\pl_t  \mathcal{R}_{1,i}^{m,\al,\ba,k}   |^2    dy, \notag\\
  \mathcal{Q}_3^{m,\al,\ba}= & (1+t)^{ 2m -1 }  \sum_{i }  \int  \sigma^{\iota+ |\alpha| }  | \mathcal{R}_{2,i}^{m,\al,\ba }  |^2    dy
 .\notag
\end{align}

In what follows, we will use the estimates stated in Section \ref{subsec-1} to control $\mathcal{Q}_k^{m,\al,\ba}$ $(k=1,2,3)$. Due to \ef{6.7-1a}, \ef{commutator2}, \ef{5.29},  and
\begin{align*}
 \mathcal{R}_{1,i}^{m,\al,\ba,k} =  &      \pl_t^m\pl^\alpha\bar\pl^\beta ( A_i^k {J}^{1-\ga}  -  \da_i^k )  + {J}^{1-\ga} (   A^k_r  A^s_i   \pl_t^{m} \pl^\alpha \bar\pl^\beta  \pl_s \omega^r  \\
&   + \iota^{-1} A_i^k A^s_r    \pl_t^{m} \pl^\alpha \bar\pl^\beta   \pl_s \omega^r    )     + {J}^{1-\ga} (   A^k_r  A^s_i \pl_t^{m} \pl^\alpha  [\pl_s,  \bar\pl^\beta ]  \omega^r  \\
& + \iota^{-1} A_i^k A^s_r    \pl_t^{m} \pl^\alpha [ \pl_s, \bar\pl^\beta ]   \omega^r    ),
\end{align*}
we have
\begin{align*}
|\mathcal{R}_{1,i}^{m,\al,\ba,k} | \les   \widetilde{\mathcal{I}}^{m,|\al|,|\ba|}  + \sum_{0\le j  \le |\ba|-1}|\pl_t^{m} \pl^{|\alpha|+1} \bar\pl^{j} \oa|,
\end{align*}
where $\widetilde{\mathcal{I}}^{m,|\al|,|\ba|}$ is defined in \ef{5.30-1}.
This, together with  \ef{5.30-2},  \ef{assume} and \ef{7.10}, implies that
\begin{align}
\mathcal{Q}_1^{m,\al,\ba}
\les  &    (1+t)^{-1}  \mathscr{E}(t) \sum_{0\le j \le m+|\al|+|\ba|-1} \mathscr{E}_j (t)+  (1+t)^{-1}  \sum_{0\le j \le |\ba|-1}  \mathscr{E}_{m+|\al|+j}(t) \notag \\
 \les & (1+t)^{-1} (\ea_0^2 +1 ) \sum_{0\le j \le m+|\al|+|\ba|-1} \mathfrak{E}_j (t) \le  \sum_{0\le j \le m+|\al|+|\ba|-1} \mathfrak{D}_j. \label{5.31-1}
\end{align}
Similarly, we have
$$
|\pl_t \mathcal{R}_{1,i}^{m,\al,\ba,k} | \les     \widetilde{\mathcal{I}}^{m+1,|\al|,|\ba|} + \sum_{0\le j  \le |\ba|-1}\lt( |\pl_t\pl\oa||\pl_t^{m} \pl^{|\alpha|+1} \bar\pl^{j} \oa|  +|\pl_t^{m+1} \pl^{|\alpha|+1} \bar\pl^{j} \oa|\rt),
$$
so that
\begin{align}
\mathcal{Q}_2^{m,\al,\ba} \les  &  (1+t)^{-1} \lt(  \mathscr{E}  \sum_{0\le j \le m+|\al|+|\ba|} \mathscr{E}_j   + \sum_{0\le j \le |\ba|-1} \lt( \mathscr{E} \mathscr{E}_{m+|\al|+j} + \mathfrak{E}^{m+1, |\al|,j}_{II}  \rt) \rt)(t) \notag\\
\les & |
 \ba|\mathfrak{D}^{m+1,|\al|,|\ba|-1} +   \ea_0^2 \mathfrak{D}_{m+|\al|+|\ba|} +  \sum_{0\le j \le m+|\al|+|\ba|-1} \mathfrak{D}_j.      \label{5.31-2}
\end{align}
It needs more works to bound $\mathcal{Q}_3^{m,\al,\ba}$.  In view of \ef{est1}, we see  that
\begin{align*}
|\mathcal{R}_{2,i}^{m,\al,\ba} | \les & \sum_{0\le j\le |\ba|-1} \sum_{k} ( \sa |\pl_t^m \pl^{|\al|+1}\bar\pl^j H^k_i| +|\pl_t^m \pl^{|\al|}\bar\pl^j H^k_i| ) \\
& +|\al|\sum_{0\le j\le |\ba|+1}\sum_{k}|\pl_t^m \pl^{|\al|-1}\bar\pl^j H^k_i|,
 \end{align*}
where $H^k_i=A^k_i J^{1-\ga}-\da^k_i$.  This, together with   \ef{7.9}, \ef{6.7-1a}, \ef{5.29} and \ef{5.30-1}, implies that
\begin{align}
 |\mathcal{R}_{2,i}^{m,\al,\ba} |
   \les & \sum_{0\le j\le |\ba|-1} \lt( \sa  |\pl_t^m\pl^{|\al| +1 } \bar\pl^j \pl \oa |+ \sa \widetilde{\mathcal{I}}^{m ,|\al|+1,j}  +  |\pl_t^m\pl^{|\al|  } \bar\pl^j \pl \oa | + \widetilde{\mathcal{I}}^{m ,|\al| ,j} \rt)  \notag\\
  &+  |\al| \sum_{0\le j\le |\ba|+1} \lt(|\pl_t^m\pl^{|\al|-1  } \bar\pl^j \pl \oa | + \widetilde{\mathcal{I}}^{m ,|\al|-1 ,j}  \rt) . \label{5.18-1}
\end{align}
Due to \ef{commutator2}, \ef{hard} and \ef{6.7-1c}, one has that
\begin{align*}
&\int \sa^{\iota+|\al|}\sum_{0\le j\le |\ba|-1} \lt(|\sa \pl_t^m\pl^{|\al|+1}\bar\pl^{j}\pl\oa|  + |\pl_t^m\pl^{|\al|}\bar\pl^{j}\pl\oa| \rt)^2 dy\notag\\
\les &   \int \sa^{\iota+|\al|}\sum_{0\le j\le |\ba|-1}  \sum_{0\le l\le j}\lt(|\sa \pl_t^m\pl^{|\al|+2}\bar\pl^{l}\oa|^2  + |\pl_t^m\pl^{|\al|+1}\bar\pl^{l}\oa|^2 \rt) dy
\notag\\
\les &   \int \sa^{\iota+|\al|+2}\sum_{0\le j\le |\ba|-1}  \sum_{0\le l\le j}\lt(|\pl_t^m\pl^{|\al|+2}\bar\pl^{l}\oa|^2  + |\pl_t^m\pl^{|\al|+1}\bar\pl^{l}\oa|^2 \rt) dy\\
\le & (1+t)^{-2m}\sum_{0\le j\le |\ba|-1}  \sum_{0\le l\le j} \lt(\mathfrak{E}_{II}^{m,|\al|+1,l} + \mathfrak{E}_{II}^{m,|\al|,l}\rt)(t),
\end{align*}
which means
\begin{align*}
&(1+t)^{2m-1}\int \sa^{\iota+|\al|}\sum_{0\le j\le |\ba|-1} \lt(|\sa \pl_t^m\pl^{|\al|+1}\bar\pl^{j}\pl\oa|  + |\pl_t^m\pl^{|\al|}\bar\pl^{j}\pl\oa| \rt)^2 dy \notag \\
\les & |\ba|\mathfrak{D}^{m,|\al|+1,|\ba|-1}(t)+|\ba|\sum_{0\le j\le m+|\al|+|\ba|-1}\mathfrak{D}_j(t).
\end{align*}
Notice that for $2\le m+|\al|+|\ba|\le [\iota]+6$, or $ m+|\al|+|\ba|=[\iota]+7$ with $|\al|\ge 1$,
$$
(1+t)^{ 2m   }  \int  \sigma^{\iota+ |\alpha|} \big| \widetilde{\mathcal{I}}^{m,|\al|,|\ba|}\big|^2      dy \les  \mathscr{E}(t) \sum_{0\le j\le m+|\al|+|\ba|}\mathscr{E}_j(t),
$$
which can be proved in a similar way to deriving \ef{5.30-2}.
This, together with \ef{5.30-2}, \ef{assume} and \ef{7.10}, implies that
\begin{align*}
&  (1+t)^{2m-1}\int \sa^{\iota+|\al|}\sum_{0\le j\le |\ba|-1} \lt( | \sa \widetilde{\mathcal{I}}^{m ,|\al|+1,j}|^2 +  | \widetilde{\mathcal{I}}^{m ,|\al| ,j}|^2 \rt)dy \notag \\
 \les &  (1+t)^{-1} \sum_{0\le j\le |\ba|-1} \mathscr{E}(t)\sum_{0\le l \le m+|\al|+j }\mathscr{E}_l(t)\notag\\
\les &  \sum_{0\le j\le |\ba|-1} \ea_0^2  \sum_{0\le l \le m+|\al|+j } \mathfrak{D}_{l}  (t) \les \ea_0^2 |\ba|  \sum_{0\le j\le m+|\al|+|\ba|-1}\mathfrak{D}_j(t).
\end{align*}
Similarly, we can deal with the second line of \ef{5.18-1}, and  obtain
\begin{align}
 \mathcal{Q}_3^{m,\al,\ba} \les
 |\al| \mathfrak{D}^{m,|\al|-1, |\ba|+ 1} + |\ba| \mathfrak{D}^{m,|\al|+1, |\ba|-1} + (  |\al| +|\ba| )\sum_{ 0\le j \le m+|\al|+|\ba|-1}\mathfrak{D}_j . \label{5.31-3}
\end{align}

Now, it is easy to see that \ef{4.2-1}  is a conclusion of \ef{4.15c}, \ef{4.15d}, \ef{8.23b}, \ef{5.31-1}, \ef{5.31-2} and \ef{5.31-4}. In fact, \ef{4.2-2} can be obtained similarly so that we omit the detail of its proof.

{\em Step 3}.
In this step,   we prove \ef{hig} for $j=1$ and take the proof as an example to explain why we   deal with the tangential derivatives and normal ones using different estimates.
Indeed,
\begin{align}\label{step1'}
   \mathfrak{E}_{1}(t) + \int_0^t \mathfrak{D}_{1}(s)ds \les    \sum_{k=0,1} \lt( \mathfrak{E}_k  (0)+   \mathfrak{V}_k  (t)  + \int_0^t(1+s)^{-1} \mathfrak{V}_k (s)ds \rt) =  \mathcal{X}(t)
   \end{align}
is a consequence of the following estimates: for any $ \vea\in(0,1)$,
\begin{subequations}\label{step1}\begin{align}
   \mathfrak{E}^{1,0,0}(t) + \int_0^t \mathfrak{D}^{1,0,0}(s)ds \les &  ( \mathfrak{E}_0 + \mathfrak{E}^{1,0,0}  ) (0)+  (\mathfrak{V}_0 + \mathfrak{V}^{1,0,0} ) (t) \notag\\
&     + \int_0^t(1+s)^{-1}(\mathfrak{V}_0 + \mathfrak{V}^{1,0,0} ) (s)ds, \label{step1.1}\\
  \mathfrak{E}^{0,1,0}(t) + \int_0^t \mathfrak{D}^{0,1,0}(s)ds \les   & \vea^{-1}\mathcal{X}(t)   + \vea \int_0^t \mathfrak{D}^{0,0,1}(s)ds , \label{step1.3}\\
   \mathfrak{E}^{0,0,1}(t) + \int_0^t \mathfrak{D}^{0,0,1}(s) ds \les  & \mathcal{X}(t)  + \int_0^t \mathfrak{D}^{0,1,0}(s)ds. \label{step1.2}
 \end{align}\end{subequations}

When $|\al|=|\ba|=0$ and $m=1$, we have
$\mathcal{R}_{1,i}^{m,\al,\ba,k}=0$ and  $\mathcal{R}_{2,i}^{m,\al,\ba}=0$, so that $$\mathcal{Q}_i^{m,\al,\ba}=0 \ \  {\rm for} \ \  i=1,2,3, $$
and  \ef{step1.1} can be obtained easily by use of  \ef{low}.

When $m=|\ba|=0$ and $|\al|=1$, we have that
$\mathcal{R}_{1,i}^{m,\al,\ba,k}=0$ and  $|\mathcal{R}_{2,i}^{m,\al,\ba}|\les |\pl \bar\pl \oa| + |\pl \oa|$, due to \ef{est1}, \ef{7.12h} and \ef{commutator2}, so that
$$\mathcal{Q}_i^{m,\al,\ba} =0   \ \ {\rm for} \ \  i=1,2, \ \  {\rm and } \ \   \mathcal{Q}_3^{m,\al,\ba}\les \mathfrak{D}_0+\mathfrak{D}^{0,0,1},$$
due to \ef{6.7-1c}.
 This, together with \ef{5.3-1}, \ef{low} and \ef{step1.1}, implies \ef{step1.3}.

When $m= |\al|=0$ and $|\ba|=1$,
it follows from \ef{7.12}, \ef{est1}, \ef{7.12h} and \ef{commutator2} that  $|\mathcal{R}_{1,i}^{m,\al,\ba,k}|\les |\pl \oa|$, $|\pl_t\mathcal{R}_{1,i}^{m,\al,\ba,k}|\les |\pl_t \pl \oa|$ and  $|\mathcal{R}_{2,i}^{m,\al,\ba}|\les \sa |\pl^2 \oa| +  |\pl \oa|$,   which implies
$$ \mathcal{Q}_1^{m,\al,\ba}\les (1+t)^{-1} \mathfrak{E}_0(t) \le  \mathfrak{D}_0, \ \ \mathcal{Q}_2^{m,\al,\ba}\les  \mathfrak{D}^{1,0,0} \ \ {\rm and} \ \  \mathcal{Q}_3^{m,\al,\ba}\les  \mathfrak{D}^{0,1,0},$$
due to \ef{6.7-1c}. Indeed, $\mathcal{Q}_3^{m,\al,\ba}$ follows from
\begin{align}\label{5.18}
  \int \sa^\iota   |\pl \oa|^2  dy
\les  \int \sa^{\iota+2} (  |\pl^2 \oa|^2 + |\pl \oa|^2 )dy ,
\end{align}
due to \ef{hard}. This, together with  \ef{low} and \ef{step1.1}, proves \ef{step1.2}.

If we used \ef{4.15b}, instead of \ef{5.3-1}, to bound $\mathcal{H}_2^{m,\al,\ba}$  in the case of $m=|\ba|=0$ and $|\al|=1$, we would get
\begin{align}
  \mathfrak{E}^{0,1,0}(t) + \int_0^t \mathfrak{D}^{0,1,0}(s)ds \les  \mathcal{X}(t)+   \int_0^t \mathfrak{D}^{0,0,1}(s)ds, \label{4.14}
 \end{align}
instead of \ef{step1.3}. Apparently,  \ef{step1'}  cannot follow from \ef{step1.1}, \ef{step1.2} and \ef{4.14}. This simple case explains why  \ef{5.3-1}, instead of \ef{4.15b}, is needed to deal with normal derivatives.

{\em Step 4}.
 We use the mathematical induction to prove \ef{hig}.  Clearly,  \ef{hig} holds for $j=0,1$, due to \ef{low} and \ef{step1'}.  Suppose that \ef{hig} holds for $j=0,\cdots, l-1$, that is,
\begin{align}
  & \mathfrak{E}_{j}(t) + \int_0^t \mathfrak{D}_{j}(s)ds \notag \\
    \les & \sum_{0\le k \le j} \lt( \mathfrak{E}_k  (0)+   \mathfrak{V}_k  (t)  + \int_0^t(1+s)^{-1} \mathfrak{V}_k (s)ds \rt), \ \  j=0,1,\cdots, l-1.\label{step2}
 \end{align}
It suffices to prove \ef{step2} holds for $j=l$.

It follows from \ef{5.4}, \ef{5.2} and \ef{4-1}  that for any $\vea\in(0,1)$,
\begin{subequations}\label{h12}\begin{align}
& \frac{d}{dt} \mathcal{E}_1^{m,\al,\ba}(t)+  \mathcal{D}_1^{m,\al,\ba}(t)   \les (1+t)^{ - 2m -1 } \lt( \mathcal{H}_g^{m,\al,\ba} +  \mathcal{H}_b^{m,\al,\ba} \rt) , \label{6.5a}\\
& \frac{d}{dt} (4\mathcal{E}_1^{m,\al,\ba} + \mathcal{E}_2^{m,\al,\ba})(t)+  (4\mathcal{D}_1^{m,\al,\ba}+ \mathcal{D}_2^{m,\al,\ba})(t) \notag \\
&\qquad \les (1+t)^{ - 2m -1 } \mathcal{H}_b^{m,\al,\ba}   + (1+t)^{ - 2m } \lt( \mathcal{F}_g^{m,\al,\ba} +  \mathcal{F}_b^{m,\al,\ba} \rt), \label{6.5b}
\end{align}\end{subequations}
where $\mathcal{H}_g^{m,\al,\ba}$, $\mathcal{H}_b^{m,\al,\ba}$, $\mathcal{F}_g^{m,\al,\ba}$ and $\mathcal{F}_b^{m,\al,\ba}$ are  defined in \ef{8.23c}.
 Integrate the product of \ef{6.5b} and $(1+t)^k$ over $[0,t]$ from $k=0$ to $k=2m$ step by step, and then integrate the product of \ef{6.5a} and $(1+t)^{2m+1}$ over $[0,t]$ to obtain that for any $\vea\in (0,1)$,
\begin{align}\label{h6.9}
\mathfrak{E}^{m,|\al|,|\ba|}(t) + \int_0^t \mathfrak{D}^{m,|\al|,|\ba|}(s)ds \les  \mathfrak{P}_l(t)  \ \ {\rm when} \ \  m+|\al|+|\ba|=l,
\end{align}
where
\begin{align*}
\mathfrak{P}_l(t) =  & (\vea + \ea_0 +  \vea^{-1} \ea_0^2) \int_0^t \mathfrak{D}_{l}(s) ds\\
 & + \vea^{-1}\sum_{0\le k \le l} \lt( \mathfrak{E}_k  (0)+   \mathfrak{V}_k  (t)  + \int_0^t(1+s)^{-1} \mathfrak{V}_k (s)ds \rt) .
\end{align*}
Here \ef{6.9-1}, \ef{5.4-3}, \ef{6.8-1}, \ef{5.4-4}, the Grownwall inequality and the induction assumption \ef{step2} have been used to derive  \ef{h6.9}.

Indeed, the mathematical induction on $m$ has been used to prove \ef{h6.9}.
Clearly, \ef{h6.9} holds for $m=l$, since $\mathcal{H}_b^{l,0,0}=\mathcal{F}_b^{l,0,0}=0$. When $m=l-1$, we have
\begin{align*}
& \mathfrak{E}^{l-1,1,0}(t) + \int_0^t \mathfrak{D}^{l-1,1,0}(s)ds \les \mathfrak{P}_l(t) +\vea^{-1} \int_0^t \mathfrak{D}^{l,0,0}(s)\les \mathfrak{P}_l(t), \\
& \mathfrak{E}^{l-1,0,1}(t) + \int_0^t \mathfrak{D}^{l-1,0,1}(s)ds \les \mathfrak{P}_l(t) + \vea^{-1} \int_0^t (\mathfrak{D}^{l,0,0}+ \mathfrak{D}^{l-1,1,0} )(s)\les \mathfrak{P}_l(t),
\end{align*}
which implies that \ef{h6.9} holds for $m=l-1$. Suppose that
\begin{align}
\sum_{|\al|+|\ba|=j} \mathfrak{E}^{l-j,|\al|,|\ba|}(t) + \sum_{|\al|+|\ba|=j} \int_0^t \mathfrak{D}^{l-j,|\al|,|\ba|}(s)ds \les \mathfrak{P}_l(t) , \ \  j=0,1,2,\cdots,k-1. \label{6-9-1}
\end{align}
It is enough to prove \ef{6-9-1} holds for $j=k$. For $j=k$, we have
\begin{align*}
& \sum_{|\al|+|\ba|=k, \  |\al| \ge 1} \lt\{ \mathfrak{E}^{l-k,|\al|,|\ba|}(t) + \int_0^t \mathfrak{D}^{l-k,|\al|,|\ba|}(s)ds \rt\}  \notag \\
& \quad \les   \mathfrak{P}_l(t) + \vea^{-1}\sum_{|\al|+|\ba|=k-1}\int_0^t \mathfrak{D}^{l-k+1,|\al|,|\ba|}(s)ds \les  \mathfrak{P}_l(t),\\
& \mathfrak{E}^{l-k,0,k}(t) + \int_0^t \mathfrak{D}^{l-k,0,k}(s)ds \notag \\
& \quad  \les \mathfrak{P}_l(t) + \vea^{-1}\int_0^t (\mathfrak{D}^{l-k+1,0,k-1} + \mathfrak{D}^{l-k,1,k-1} ) (s)ds \les  \mathfrak{P}_l(t).
\end{align*}
So, \ef{6-9-1} holds for $j=k$, and we obtain \ef{h6.9}.

It follows from \ef{h6.9} that
\begin{align*}
 \mathfrak{E}_{l}(t) & + \int_0^t \mathfrak{D}_{l}(s)ds
 \les     (\vea + \ea_0 +  \vea^{-1} \ea_0^2) \int_0^t \mathfrak{D}_{l}(s) ds\\
 & + \vea^{-1}\sum_{0\le k \le l} \lt( \mathfrak{E}_k  (0)+   \mathfrak{V}_k  (t)  + \int_0^t(1+s)^{-1} \mathfrak{V}_k (s)ds \rt),
\end{align*}
which implies that \ef{step2} holds for $j=l$, by choosing $\vea=\ea_0$ and using the smallness of $\ea_0$. This finishes the proof of this Lemma. \hfill $\Box$

\section{Curl  estimates}\label{sec-curl}
This section devotes to performing the estimate for sobolev norms of curl, $\mathfrak{V}^{m,n,l}$, which is needed to bound the energy functional as we see in Proposition \ref{newedv}.

\begin{prop}\label{prop-curl}  Let $\oa(t,y)$ be a solution to problem \ef{newsystem} in the time interval $[0, T]$ satisfying \ef{assume}.  Let $m,n,l$ be nonnegative integers satisfying $m+n+l\le [\iota]+7$,  then  for $t\in [0,T]$,
\begin{subequations}\label{curl}\begin{align}
& \mathfrak{V}^{0,n,l}(t)
 \les \sum_{i=0,1}  \lt\|\sa^{\frac{\iota+n+1}{2}}  \pl^{n}\bar\pl^{l}{\rm curl}_\eta \pl_t^i \oa \big|_{t=0} \rt\|_{L^2}^2  + \sum_{0\le k\le l-1} \mathfrak{E}_{II}^{0,n,k}(t)  \notag\\
 & \quad   +   \sup_{s\in [0,t]}  \mathscr{E}(s)   \sum_{0\le j \le n+l } \lt(\sup_{s\in [0,t]} \mathscr{E}_{j}(s) + \ln (1+t)   \int_0^t (1+s)^{-1}   \mathscr{E}_{j}(s) ds \rt),  \ \ m=0,  \label{curl-a}\\
&  \mathfrak{V}^{m,n,l}(t)
  \les (1+t)^{-2}  \lt\|\sa^{\frac{\iota+n+1}{2}}  \pl^{n}\bar\pl^{l}{\rm curl}_\eta \pl_t\oa \big|_{t=0}\rt\|_{L^2}^2   +  \mathscr{E}(t) \sum_{0\le j\le m+n+l-1}\mathscr{E}_j(t)  \notag \\
 &\quad + (1+t)^{-2} \sup_{s\in [0,t]}  \mathscr{E}(s)   \sum_{0\le j \le m+ n+ l }  \sup_{s\in [0,t]} \mathscr{E}_{j}(s)   +  \sum_{0\le k\le l-1} \mathfrak{E}_{II}^{m,n,k}(t)   , \ \   m\ge 1. \label{curl-b}
\end{align}\end{subequations}
Moreover, we have for $n+ l= [\iota]+7$,
\begin{align}
& \mathfrak{V}^{1,n,l}(t)
 \les
   \lt\|\sa^{\frac{\iota+n+1}{2}}  \pl^{n}\bar\pl^{l}{\rm curl}_\eta \pl_t\oa \big|_{t=0}\rt\|_{L^2}^2    \notag \\
   & \qquad +   \sup_{s\in [0,t]}  \mathscr{E}(s)   \sum_{0\le j \le  n+ l }  \sup_{s\in [0,t]} \mathscr{E}_{j}(s)
     +  \sum_{0\le k\le l-1} \mathfrak{E}_{II}^{1,n,k}(t) , \ \  t\in [0,T] . \label{curl-c}
\end{align}
\end{prop}

 {\em Proof}.
Equation \ef{3-1-2} can be rewritten in the form of
\bee\label{8.15}
\ta     \pl_t^2 \oa +   (\ta+ 2\ta_t )   \pl_t  \oa  +  (3\ga-1)^{-1}\ta^{2-3\ga} \eta     + \frac{\ga}{\ga-1} \ta^{2-3\ga} \nabla_\eta \lt(\bar\rho_0^{\ga-1} {J}^{1-\ga}\rt)  = 0,
\eee
Let ${\rm curl}_\eta$ act on it, and use the fact ${\rm curl}_\eta \eta =0$ and ${\rm curl}_\eta \nabla_\eta =0$ to give
\bee\label{3-5-2}
\ta {\rm curl}_\eta \pl_t^2 \oa  +(2\ta_t + \ta) {\rm curl}_\eta    \pl_t \oa =0.
\eee
Commuting $\pl_t$ with ${\rm curl}_\eta$ and noting the integrating-factor $ \ta^2$, we have
\be\label{3-5-4}
{\rm curl}_\eta    \pl_t \oa  = \lt\{  \ta^2(0) {\rm curl}_\eta    \pl_t \oa \big|_{t=0}   +   \int_0^t  e^{\tau} \ta^2(\tau) \lt[\pl_\tau , {\rm curl}_\eta  \rt]   \pl_\tau \oa d\tau \rt\} e^{-t} \ta^{-2}(t).
\ee
Commute $\pl_t$ with ${\rm curl}_\eta$ again, and integrate the resulting equation over time to obtain
\begin{align}
& {\rm curl}_\eta \oa =  {\rm curl}_\eta \oa \big|_{t=0} +    \ta^2 (0) {\rm curl}_\eta    \pl_t \oa \big|_{t=0} \int_0^t  e^{-s} \ta^{-2}(s) ds \notag \\
& \quad + \int_0^t \lt[\pl_s , {\rm curl}_\eta  \rt]   \oa  ds  + \int_0^t   e^{-s} \ta^{-2}(s) \int_0^s e^{\tau} \ta^2(\tau) \lt[\pl_\tau , {\rm curl}_\eta  \rt]   \pl_\tau \oa  d\tau ds
 .\label{3-5-6}
\end{align}
In what follows, we use the formulae \ef{3-5-4} and \ef{3-5-6} to prove the estimates \ef{curl} and \ef{curl-c}.

{\em Step 1}.
In this step, we prove that for $|\al|+|\ba|\le [\iota]+7$,
\begin{align}
& \lt\|\sa^{\frac{\iota+|\al|+1}{2}}  \pl^{\al}\bar\pl^{\ba}{\rm curl}_\eta\oa \rt\|_{L^2}^2
 \les
\sum_{i=0,1}  \lt\|\sa^{\frac{\iota+|\al|+1}{2}}  \pl^{\al}\bar\pl^{\ba}{\rm curl}_\eta \pl_t^i \oa  \rt\|_{L^2}^2(t=0) \notag \\
& \quad   +   \sup_{s\in [0,t]}  \mathscr{E}(s)   \sum_{0\le j \le |\al|+|\ba| } \lt(\sup_{s\in [0,t]} \mathscr{E}_{j}(s) + \ln (1+t)   \int_0^t (1+s)^{-1}   \mathscr{E}_{j}(s) ds \rt). \label{7.26a}
\end{align}
Take $\pl^{\al}\bar\pl^{\ba}$ onto \ef{3-5-6} to obtain
 \begin{align}
&\pl^{\al}\bar\pl^{\ba}{\rm curl}_\eta \oa =  \pl^{\al}\bar\pl^{\ba}{\rm curl}_\eta \oa\big|_{t=0}  +    \ta^2(0) \pl^{\al}\bar\pl^{\ba} {\rm curl}_\eta    \pl_t \oa \big|_{t=0} \int_0^t  e^{-s} \ta^{-2}(s) ds \notag \\
&\quad + \int_0^t \pl^{\al}\bar\pl^{\ba} \lt[\pl_s , {\rm curl}_\eta  \rt]   \oa  ds   + \int_0^t   e^{-s} \ta^{-2}(s) \int_0^s e^{\tau} \ta^2(\tau) \pl^{\al}\bar\pl^{\ba} \lt[\pl_\tau , {\rm curl}_\eta  \rt]   \pl_\tau \oa  d\tau ds. \label{7.13}
\end{align}
Clearly, \ef{7.26a} holds if the second line of \ef{7.13} can be bounded by
 \begin{subequations}\label{8.8a}\begin{align}
&  \lt\|  \sa^{ \frac{\iota+|\al|+1}{2}}   \int_0^t \pl^{\al}\bar\pl^{\ba} \lt[\pl_s , {\rm curl}_\eta  \rt]   \oa  ds   \rt\|_{L^2}^2
\les   \mathscr{E}(0) \sum_{0\le j \le |\al|+|\ba|  } \mathscr{E}_{j}(0)
 + \mathscr{E}(t) \sum_{0\le j \le |\al|+|\ba|  } \mathscr{E}_{j}(t)
\notag  \\
& \qquad\qquad \qquad + \ln (1+t) \sup_{s\in [0,t]}\mathscr{E}(s) \sum_{0\le j \le |\al|+|\ba|  } \int_0^t (1+s)^{-1}   \mathscr{E}_{j}(s) ds, \label{8.1c}\\
& \lt\|  \sa^{ \frac{\iota+|\al|+1}{2}}  \int_0^t   e^{-s} \ta^{-2}(s) \int_0^s e^{\tau} \ta^2(\tau) \pl^{\al}\bar\pl^{\ba} \lt[\pl_\tau , {\rm curl}_\eta  \rt]   \pl_\tau \oa  d\tau ds  \rt\|_{L^2}^2 \notag \\
&\qquad \qquad \qquad  \les     \sup_{s\in [0,t]} \mathscr{E}(s)  \sum_{0\le j \le |\al|+|\ba| }\sup_{s\in [0,t]} \mathscr{E}_{j}(s). \label{7.15-2}
\end{align}\end{subequations}

We first prove \ef{8.1c}. It follows from \ef{5.29}   that
\begin{align*}
\pl^{ \al }\bar\pl^{ \ba } \lt(  \pl_t [{\rm curl}_\eta  \oa]_l  -   [{\rm curl}_\eta \pl_t  \oa  ]_l   \rt)
=\pl^{ \al }\bar\pl^{ \ba }(\epsilon^{ljk} ( \pl_r  \oa_k)  \pl_t A^r_j  )
  =   \mathcal{Y}^{\al, \ba}_{1, l} +  \mathcal{Y}^{\al, \ba}_{2, l} ,
\end{align*}
where
\begin{align*}
& \mathcal{Y}^{\al, \ba}_{1, l} =   \pl_t \lt(  \epsilon^{ljk}  (\pl_r   \oa_k) \pl^{ \al }\bar\pl^{ \ba }( A^r_j-\da^r_j) \rt)
- \epsilon^{ljk} (\pl_t  \pl_r    \oa_k)  \pl^{ \al }\bar\pl^{ \ba }( A^r_j-\da^r_j)  , \\
& \lt|\mathcal{Y}^{\al, \ba}_{2, l}\rt|\les  \sum_{\substack{ 0\le j\le |\al|, \ 0\le k \le |\ba|,
\  j+k\le |\al|+|\ba|-1
  }}  \mathcal{I}^{1,j,k}   \lt|  \pl^{|\al| -j}\bar\pl^{|\ba|-k} \pl \oa\rt|.
\end{align*}
Clearly,
  \begin{align}
& \lt\|  \sa^{ \frac{\iota+|\al|+1}{2}}  \int_0^t   \mathcal{Y}^{\al, \ba}_{1, l}   ds  \rt\|_{L^2}\les
\lt\|   \sa^{ \frac{\iota+|\al|+1}{2}}   |\pl \oa| |   \pl^{ \al }\bar\pl^{ \ba }( A^r_j-\da^r_j) |  \rt\|_{L^2}(0)
  \notag\\
&   +
  \lt\|   \sa^{ \frac{\iota+|\al|+1}{2}}  |\pl \oa| | \pl^{ \al }\bar\pl^{ \ba }( A^r_j-\da^r_j)|  \rt\|_{L^2}(t)
  + \int_0^t   \lt\|   \sa^{ \frac{\iota+|\al|+1}{2}}   |\pl_s \pl \oa| | \pl^{ \al }\bar\pl^{ \ba }( A^r_j-\da^r_j)|  \rt\|_{L^2} ds  . \label{8.1a}
\end{align}
Due to \ef{7.9}, \ef{5.29} and \ef{5.30-1}, one has
$|   \pl^{ \al }\bar\pl^{ \ba }( A^r_j-\da^r_j) | \les   |\pl^{|\al|}\bar\pl^{|\ba|}\pl \oa |+\widetilde{\mathcal{I}}^{0,|\al|,|\ba|}$. This, together with  \ef{a}, \ef{verify},  \ef{commutator2}    and \ef{5.30-2}, implies that
\begin{align*}
\lt\|   \sa^{ \frac{\iota+|\al|+1}{2}}   |\pl \oa| |   \pl^{ \al }\bar\pl^{ \ba }( A^r_j-\da^r_j) |  \rt\|_{L^2}^2 \les  \|\pl \oa\|_{L^\iy}^2 \lt\|   \sa^{ \frac{\iota+|\al|+1}{2}}    |\pl^{|\al|}\bar\pl^{|\ba|}\pl \oa |  \rt\|_{L^2}^2 \\
+ \ea_0^2 \lt\|   \sa^{ \frac{\iota+|\al|+1}{2}}   \widetilde{\mathcal{I}}^{0,|\al|,|\ba|}  \rt\|_{L^2}^2
\les  \mathscr{E}(t) \sum_{0\le j \le |\al|+|\ba|  } \mathscr{E}_{j}(t),
\end{align*}
which gives the bounds for the first two terms on the right hand side of \ef{8.1a}.
It follows from \ef{5.29}, \ef{5.30-1} and \ef{5.30-2} that for  $|\al|+|\ba|\ge 1$,
\begin{align*}
 \lt\|   \sa^{ \frac{\iota+|\al|+1}{2}}   |\pl_t \pl \oa| | \pl^{ \al }\bar\pl^{ \ba }( A^r_j-\da^r_j)|  \rt\|_{L^2}
 \les
  \lt\|   \sa^{ \frac{\iota+|\al|+1}{2}} |\pl_t\pl \oa| \mathcal{I}^{0,|\al|,|\ba|} \rt\|_{L^2}^2
   \\
   \le    \lt\|   \sa^{ \frac{\iota+|\al|+1}{2}}  \widetilde{\mathcal{I}}^{1,|\al|,|\ba|}\rt\|_{L^2}^2
\les   (1+t)^{-2} \mathscr{E}(t) \sum_{0\le j \le |\al|+|\ba|  } \mathscr{E}_{j}(t),
\end{align*}
which, together with  the H$\ddot{o}$lder inequality, implies that for  $|\al|+|\ba|\ge 1$,
 \begin{align}
& \lt(\int_0^t   \lt\|   \sa^{ \frac{\iota+|\al|+1}{2}}   |\pl_s \pl \oa| | \pl^{ \al }\bar\pl^{ \ba }( A^r_j-\da^r_j)|  \rt\|_{L^2} ds\rt)^2  \notag\\
  \les
 &\int_0^t (1+s)^{-1} ds  \int_0^t (1+s)^{-1} \lt( (1+s)^2    \lt\|   \sa^{ \frac{\iota+|\al|+1}{2}} |\pl_s \pl \oa | |  \pl^{ \al }\bar\pl^{ \ba }( A^r_j-\da^r_j) |  \rt\|_{L^2}^2  \rt)  ds \notag\\
\les & \ln (1+t) \int_0^t (1+s)^{-1} \mathscr{E}(s) \sum_{0\le j \le |\al|+|\ba|  } \mathscr{E}_{j}(s)   ds  \notag\\
\le &  \ln (1+t) \sup_{s\in [0,t]}\mathscr{E}(s)  \sum_{0\le j \le |\al|+|\ba|  } \int_0^t (1+s)^{-1}   \mathscr{E}_{j}(s) ds. \label{8.14a}
\end{align}
It is easy to show that \ef{8.14a} also holds for $|\al|=|\ba|=0$, so we obtain the bound  for the last term on the right hand side of \ef{8.1a}.
Similarly, we have for $|\al|+|\ba|\le [\iota]+7$,
\begin{align*}
\lt\|  \sa^{ \frac{\iota+|\al|+1}{2}}     \mathcal{Y}^{\al, \ba}_{2, l} \rt\|_{L^2}^2
\les   \lt\|   \sa^{ \frac{\iota+|\al|+1}{2}}  \widetilde{\mathcal{I}}^{1,|\al|,|\ba|}\rt\|_{L^2}^2
\les   (1+t)^{-2} \mathscr{E}(t) \sum_{0\le j \le |\al|+|\ba|  } \mathscr{E}_{j}(t),
\end{align*}
and
 \begin{align*}
  \lt\|  \sa^{ \frac{\iota+|\al|+1}{2}}  \int_0^t   \mathcal{Y}^{\al, \ba}_{2, l}   ds  \rt\|_{L^2}^2
\les       \ln (1+t) \sup_{s\in [0,t]}\mathscr{E}(s)  \sum_{0\le j \le |\al|+|\ba|  } \int_0^t (1+s)^{-1}   \mathscr{E}_{j}(s) ds.
\end{align*}
This finishes the proof of \ef{8.1c}.

Next, we prove \ef{7.15-2}.  It follows from \ef{5.29} that
\begin{align}\label{curl-12}
& \pl^{ \al }\bar\pl^{ \ba } \lt(  \pl_t [{\rm curl}_\eta \pl_t \oa]_l  -   [{\rm curl}_\eta \pl_t^2 \oa  ]_l   \rt)
=\pl^{ \al }\bar\pl^{ \ba }(\epsilon^{ljk} (\pl_t A^r_j) \pl_t \pl_r  \oa_k)
=  \mathcal{Z}^{\al, \ba}_{1,l} +\mathcal{Z}^{\al, \ba}_{2,l}   ,
\end{align}
where
\begin{align*}
& \mathcal{Z}^{\al, \ba}_{1,l}=\epsilon^{ljk} \sum_{(|h|,|g|)\in S_2 \cup S_3  } C(\al,\ba,h,g)   (\pl_t\pl^h \bar\pl^g  A^r_j) \pl_t \pl^{ \al-h }\bar\pl^{ \ba -g }  \pl_r  \oa_k     ,\\
& \lt|\mathcal{Z}^{\al, \ba}_{2,l}\rt| \les     \sum_{(j,k)\in S_1 \setminus(S_2\cup S_3)}   \mathcal{I}^{1,j,k}   \lt| \pl_t \pl^{|\al| -j}\bar\pl^{|\ba|-k}\pl \oa\rt|.
\end{align*}
Here $S_1= \{(j,k)\in \mathbb{Z}^3    \big|      0\le j\le |\al|, \ 0\le k \le |\ba|  \}$, $S_2 =  \{(j,k)\in S_1    \big|    j=k=0,  \ j=1 \ {\rm and} \ k=0, \ j=0 \ {\rm and} \ k=2     \}$ and $S_3 =  \{(j,k)\in S_1    \big|    j=|\al| \ {\rm and} \ k=|\ba|,  \ j=|\al|-1 \ {\rm and} \ k=|\ba|, \ j=|\al| \ {\rm and} \ k=|\ba|-2     \}$.

It follows from Lemmas \ref{lem-non} and \ref{lem-non-new} that
\begin{align}
 (1+t)^{2}\lt\|  \sa^{ \frac{\iota+|\al|+1}{2}}     \mathcal{Z}^{\al, \ba}_{2, l}  \rt\|_{L^2}
\les  \sqrt{\mathscr{E}(t)} \sum_{0\le j \le |\al|+|\ba|  }  \sqrt{\mathscr{E}_{j}(t)}. \label{8.4a}
\end{align}
Indeed, the case of  $3\le  2j+k \le  2|\al|+|\ba|-3$ follows from \ef{7.14'}; the case of $2j+k=1$ (with $j=0,k=1$)  follows from the same derivation of  \ef{5.30-4}  by noting $i=1$ and $m=2$; and the case of $2j+k= 2|\al|+|\ba|-1$ is the same as that of $2j+k=1$.
Notice that for $k\ge 1$,
\begin{align}
& e^{-t} \ta^{-2}(t)  \int_0^t  e^{\tau} \ta^2(\tau)  (1+\tau)^{-k} d\tau \le e^{-t}  \int_0^t  e^{\tau}   (1+\tau)^{-k} d\tau \notag\\
\le & e^{-t/2}  \int_0^{t/2}     (1+\tau)^{-k} d\tau + e^{-t}    (1+t/2)^{-k} \int_{t/2}^t  e^{\tau}  d\tau \notag\\
   \les & e^{-t/2} \ln (1+t/2) + (1+t/2)^{-k} \les (1+t)^{-k}, \label{7.16-1}
\end{align}
where $\ta_t\ge 0$ has been used to derive the first inequality. Then, we have
\begin{align}
&  \lt\|  \sa^{ \frac{\iota+|\al|+1}{2}}  \int_0^t   e^{-s} \ta^{-2}(s) \int_0^s e^{\tau} \ta^2(\tau)  \mathcal{Z}^{\al, \ba}_{2,l}  d\tau ds  \rt\|_{L^2} \notag \\
 \les   &  \sup_{\tau\in [0,t]} (1+\tau)^{2}\lt\|  \sa^{ \frac{\iota+|\al|+1}{2}}     \mathcal{Z}^{\al, \ba}_{2, l}  \rt\|_{L^2} \int_0^t (1+s)^{-2} ds \notag\\
  \les &  \sup_{\tau \in [0,t]} \sqrt{\mathscr{E}(\tau)}  \sum_{0\le j \le |\al|+|\ba|  }\sup_{\tau \in [0,t]} \sqrt{\mathscr{E}_{j}(\tau)} .\label{7.18}
\end{align}

It needs more careful works to deal with  $\mathcal{Z}^{\al, \ba}_{1,l}$.
When $(|h|,|g|)\in S_2$ and $|g|=0$,  we  integrate by parts over time to get
\begin{align}
& \int_0^s e^{\tau} \ta^2(\tau) (\pl_\tau \pl^h A^r_j) \pl_\tau  \pl^{ \al-h  }\bar\pl^{ \ba  }  \pl_r  \oa_k      d\tau
   =     \lt(  e^{\tau} \ta^2(\tau) (\pl_\tau \pl^h A^r_j)  \pl^{ \al -h }\bar\pl^{ \ba  }  \pl_r  \oa_k  \rt)\big|_{\tau=0}^s
\notag \\
  & \quad
      -\int_0^s e^{\tau} \ta^2(\tau) (\pl_\tau^2 \pl^h A^r_j)    \pl^{ \al-h  }\bar\pl^{ \ba  }  \pl_r  \oa_k      d\tau -   \int_0^s    e^{\tau} (\ta^2)_\tau(\tau) (\pl_\tau \pl^h A^r_j)  \pl^{ \al-h }\bar\pl^{ \ba  }  \pl_r  \oa_k     d\tau
      \notag \\
    &\quad  -   \int_0^s   e^{\tau} \ta^2(\tau)     (\pl_\tau \pl^h A^r_j)  \pl^{ \al-h }\bar\pl^{ \ba  }  \pl_r  \oa_k     d\tau =  \sum_{1\le r \le 4} I^{\al,\ba,h}_{r, j,k } (s). \label{8.4b}
\end{align}
Note that
\begin{align*}
&\|   \sa^{ \frac{\iota+|\al|+1}{2}} \mathcal{I}^{l,0,0}\pl^{|\al| }\bar\pl^{|\ba|}\pl \oa  \|_{L^2}\le \|\mathcal{I}^{l,0,0} \|_{L^\iy} \|   \sa^{ \frac{\iota+|\al|+1}{2}} \pl^{|\al| }\bar\pl^{|\ba|}\pl \oa  \|_{L^2} , \ \  l=1,2,\\
& \|   \sa^{ \frac{\iota+|\al|+1}{2}} \mathcal{I}^{l,1,0}\pl^{|\al|-1 }\bar\pl^{|\ba|}\pl \oa  \|_{L^2}\le \|\sa\mathcal{I}^{l,1,0} \|_{L^\iy} \|   \sa^{ \frac{\iota+|\al|-1}{2}} \pl^{|\al|-1 }\bar\pl^{|\ba|}\pl \oa  \|_{L^2}, \ \  l=1,2.
\end{align*}
Then, we can use \ef{hard} to obtain
\begin{align}
&(1+t) \lt\|   \sa^{ \frac{\iota+|\al|+1}{2}} \mathcal{I}^{1,|h|,0}\pl^{|\al|-|h| }\bar\pl^{|\ba|}\pl \oa    \rt\|_{L^2}  \notag \\
&+ (1+t)^2 \lt\|   \sa^{ \frac{\iota+|\al|+1}{2}} \mathcal{I}^{2,|h|,0}\pl^{|\al|-|h| }\bar\pl^{|\ba|}\pl \oa \rt\|_{L^2}
 \les  \sqrt{\mathscr{E}(t)} \sum_{0\le j \le |\al|+|\ba|  }  \sqrt{\mathscr{E}_{j}(t)}. \label{8.5a}
\end{align}
This, together with  \ef{5.29}, \ef{7.16-1} and $(\ta^2)_\tau \les (1+\tau)^{-1} \ta^2$, implies that
\begin{align*}
&\lt\|  \sa^{ \frac{\iota+|\al|+1}{2}}  \int_0^t   e^{-s} \ta^{-2}(s)   (I^{\al,\ba,h}_{2,j,k} + I^{\al,\ba,h}_{3,j,k})(s)   ds \rt\|_{L^2}\notag \\
\les &
\sup_{\tau\in [0,t]} \lt\{ (1+\tau)^2 \lt\|   \sa^{ \frac{\iota+|\al|+1}{2}} \mathcal{I}^{2,|h|,0}\pl^{|\al|-|h| }\bar\pl^{|\ba|}\pl \oa \rt\|_{L^2} \rt.\notag\\
 & \lt. + (1+\tau) \lt\|   \sa^{ \frac{\iota+|\al|+1}{2}} \mathcal{I}^{1,|h|,0}\pl^{|\al|-|h| }\bar\pl^{|\ba|}\pl \oa    \rt\|_{L^2} \rt\} \int_0^t (1+s)^{-2}ds \notag\\
  \les &  \sup_{\tau \in [0,t]} \sqrt{\mathscr{E}(\tau)}  \sum_{0\le j \le |\al|+|\ba|  }\sup_{\tau \in [0,t]} \sqrt{\mathscr{E}_{j}(\tau)} .
\end{align*}
Integrate by parts over time to obtain
\begin{align*}
& \int_0^t   e^{-s} \ta^{-2}(s)    I^{\al,\ba,h}_{4,j,k}(s) ds = - \int_0^t    \ta^{-2}(s)    I^{\al,\ba,h}_{4,j,k}(s) d e^{-s}\\
= &
-  e^{-s} \ta^{-2}(s)    I^{\al,\ba,h}_{4,j,k}(s)  \big|_{s=0}^t + \int_0^t e^{-s} \lt(   (\ta^{-2})_s (s)    I^{\al,\ba,h}_{4,j,k}(s) + \ta^{-2}(s)   \pl_s I^{\al,\ba,h}_{4,j,k}(s) \rt) ds,
\end{align*}
which  implies that
\begin{align*}
&\int_0^t   e^{-s} \ta^{-2}(s)    \lt( I^{\al,\ba,h}_{1,j,k} +I^{\al,\ba,h}_{4,j,k}  \rt)(s)  ds
\\
=& - \ta^2(0)    \{ (\pl_\tau \pl^h A^r_j)  \pl^{ \al-h }\bar\pl^{ \ba  }  \pl_r  \oa_k \}\big|_{\tau=0} \int_0^t e^{-s} \ta^{-2}(s) ds
\\
&+ e^{-t} \ta^{-2}(t) \int_0^t   e^{\tau} \ta^2(\tau)     (\pl_\tau \pl^h A^r_j)  \pl^{ \al-h }\bar\pl^{ \ba  }  \pl_r  \oa_k     d\tau
\\
&- \int_0^t e^{-s}   (\ta^{-2})_s (s) \int_0^s   e^{\tau} \ta^2(\tau)     (\pl_\tau \pl^h A^r_j)  \pl^{ \al-h }\bar\pl^{ \ba  }  \pl_r  \oa_k     d\tau ds.
\end{align*}
Then, we  use \ef{5.29}, \ef{7.16-1}, $-(\ta^{-2})_s \les (1+s)^{-1} \ta^{-2}$ and \ef{8.5a} to get
\begin{align*}
&\lt\|  \sa^{ \frac{\iota+|\al|+1}{2}}  \int_0^t   e^{-s} \ta^{-2}(s)   (I^{\al,\ba,h}_{1,j,k} + I^{\al,\ba,h}_{4,j,k})(s)   ds \rt\|_{L^2} \notag\\
&\les
\lt\|   \sa^{ \frac{\iota+|\al|+1}{2}} \mathcal{I}^{1,|h|,0}\pl^{|\al|-|h| }\bar\pl^{|\ba|}\pl \oa    \rt\|_{L^2}(\tau=0) \notag\\
&  +  \sup_{\tau \in [0,t]} \lt\{(1+\tau) \lt\|   \sa^{ \frac{\iota+|\al|+1}{2}} \mathcal{I}^{1,|h|,0}\pl^{|\al|-|h| }\bar\pl^{|\ba|}\pl \oa    \rt\|_{L^2}\rt\}
\lt( (1+t)^{-1} + \int_0^t (1+s)^{-2} ds \rt) \notag\\
& \les  \sup_{\tau \in [0,t]} \sqrt{\mathscr{E}(\tau)}  \sum_{0\le j \le |\al|+|\ba|  }\sup_{\tau \in [0,t]} \sqrt{\mathscr{E}_{j}(\tau)}  .
\end{align*}
 When $(|h|,|g|)\in S_2$ and $|g|\neq 0$ which means $|h|=0$ and $|g|=2$, we can obtain the same bounds  by noting that
\begin{align*}
& \int_0^s e^{\tau} \ta^2(\tau) (\pl_\tau \bar\pl^g A^r_j) \pl_\tau  \pl^{ \al  }\bar\pl^{ \ba -g  }  \pl_r  \oa_k      d\tau
   =     \lt(  e^{\tau} \ta^2(\tau) ( \bar\pl^g A^r_j) \pl_\tau \pl^{\al } \bar\pl^{ \ba -g  }  \pl_r  \oa_k  \rt)\big|_{\tau=0}^s
 \\
  & \quad
      -\int_0^s e^{\tau} \ta^2(\tau) ( \bar\pl^g A^r_j) \pl_\tau^2   \pl^{ \al   }\bar\pl^{ \ba -g  }  \pl_r  \oa_k      d\tau -   \int_0^s \pl_\tau \lt( e^{\tau} \ta^2(\tau) \rt)   ( \bar\pl^g A^r_j) \pl_\tau  \pl^{ \al }\bar\pl^{ \ba-g  }  \pl_r  \oa_k     d\tau.
\end{align*}
The case of $(|h|,|g|)\in S_3$ can be bounded similarly as that of $(|h|,|g|)\in S_2$, so we can obtain the estimate involving $\mathcal{Z}^{\al, \ba}_{1,l}$, which, together with  \ef{7.18}, proves \ef{7.15-2}.

{\em Step 2}.
In this step, we prove that  for $m\ge 1$ and $m+|\al|+|\ba|\le [\iota]+7$,
 \begin{align}
& (1+t)^{2(m+1)} \lt\|\sa^{\frac{\iota+|\al|+1}{2}} \pl_t^{m-1} \pl^{\al}\bar\pl^{\ba}{\rm curl}_\eta \pl_t\oa \rt\|_{L^2}^2 \notag \\
  \les &  \lt\|\sa^{\frac{\iota+|\al|+1}{2}}  \pl^{\al}\bar\pl^{\ba}{\rm curl}_\eta \pl_t\oa \rt\|_{L^2}^2(t=0)   +   \sup_{\tau \in [0,t]}  \mathscr{E}(\tau)   \sum_{0\le j \le m+ |\al|+|\ba|  }\sup_{\tau \in [0,t]}  \mathscr{E}_{j}(\tau) .\label{7.26c}
\end{align}

When $m=1$, apply $\pl^{\al}\bar\pl^{\ba}$  to  \ef{3-5-4} to get
\begin{align}
\pl^{\al}\bar\pl^{\ba} {\rm curl}_\eta    \pl_t \oa  =  \ta^2(0) \pl^{\al}\bar\pl^{\ba} {\rm curl}_\eta    \pl_t \oa \big|_{t=0} e^{-t} \ta^{-2}(t) \notag\\
 +    e^{-t} \ta^{-2}(t) \int_0^t  e^{\tau} \ta^2(\tau) \pl^{\al}\bar\pl^{\ba} \lt[\pl_\tau , {\rm curl}_\eta  \rt]   \pl_\tau \oa d\tau , \label{8.4c}
\end{align}
which, together with  \ef{7.16-1} and the following estimate:
\begin{align}
  \lt\|\sa^{\frac{\iota+|\al|+1}{2}}  \pl^{\al}\bar\pl^{\ba} \lt[\pl_t , {\rm curl}_\eta  \rt]   \pl_t \oa \rt\|_{L^2}  \les \sum_ {\substack{ 0\le j\le |\al|, \ 0\le k \le |\ba|
  }} \lt\|\sa^{\frac{\iota+|\al|+1}{2}}  \mathcal{I}^{1,j,k}|\pl_t  \pl^{|\al|-j}\bar\pl^{|\ba|-k} \pl \oa | \rt\|_{L^2}
 \notag  \\
\les \lt\|\sa^{\frac{\iota+|\al|+1}{2}}  \widetilde{\mathcal{I}}^{2,|\al|,|\ba|}  \rt\|_{L^2}
\les (1+t)^{-2}\sqrt{\mathscr{E}(t)} \sum_{0\le j \le  1+ |\al|+|\ba|  }  \sqrt{\mathscr{E}_{j}(t)}, \label{8.6c}
\end{align}
proves \ef{7.26c} for $m=1$.  Here \ef{5.29}, \ef{5.30-2} and \ef{5.30-1} have been used to derive \ef{8.6c}.

When $m=2$, take  $\pl_t$ onto  \ef{8.4c} and   integrate by parts over time to obtain
\begin{align}
& \pl_t \pl^{\al}\bar\pl^{\ba} {\rm curl}_\eta    \pl_t \oa  =    \ta^2(0) \pl^{\al}\bar\pl^{\ba} {\rm curl}_\eta    \pl_t \oa \big|_{t=0}\lt( e^{-t} \ta^{-2}(t) \rt)_t  \notag\\
&\qquad +\ta^2(0) \pl^{\al}\bar\pl^{\ba} \lt[\pl_t , {\rm curl}_\eta  \rt]   \pl_t \oa\big|_{t=0}  e^{-t} \ta^{-2}(t)\notag \\
&\qquad +  e^{-t} \ta^{-2}(t) \int_0^t  e^{\tau} \pl_\tau \lt( \ta^2(\tau) \pl^{\al}\bar\pl^{\ba} \lt[\pl_\tau , {\rm curl}_\eta  \rt]   \pl_\tau \oa  \rt) d\tau  \notag \\
&\qquad +  e^{-t}  (\ta^{-2})_t(t) \int_0^t  e^{\tau} \ta^2(\tau) \pl^{\al}\bar\pl^{\ba} \lt[\pl_\tau , {\rm curl}_\eta  \rt]   \pl_\tau \oa d\tau, \label{8.5b}
\end{align}
due to
\begin{align*}
&-e^{-t} \ta^{-2}(t) \int_0^t  e^{\tau} \ta^2(\tau) \pl^{\al}\bar\pl^{\ba} \lt[\pl_\tau , {\rm curl}_\eta  \rt]   \pl_\tau \oa d\tau \\
=&-e^{-t} \ta^{-2}(t) \int_0^t   \ta^2(\tau) \pl^{\al}\bar\pl^{\ba} \lt[\pl_\tau , {\rm curl}_\eta  \rt]   \pl_\tau \oa d  e^{\tau}\\
= &-\pl^{\al}\bar\pl^{\ba} \lt[\pl_t , {\rm curl}_\eta  \rt]   \pl_t \oa +\ta^2(0)   \pl^{\al}\bar\pl^{\ba} \lt[\pl_\tau , {\rm curl}_\eta  \rt]   \pl_\tau \oa \big|_{\tau=0}
e^{-t} \ta^{-2}(t)\\
&+  e^{-t} \ta^{-2}(t) \int_0^t  e^{\tau} \pl_\tau \lt( \ta^2(\tau) \pl^{\al}\bar\pl^{\ba} \lt[\pl_\tau , {\rm curl}_\eta  \rt]   \pl_\tau \oa \rt)  d\tau.
\end{align*}
In view of  \ef{5.29}, \ef{5.30-1} and \ef{5.30-2}, we see that
\begin{align}
& \lt\|\sa^{\frac{\iota+|\al|+1}{2}} \pl_t \pl^{\al}\bar\pl^{\ba} \lt[\pl_t , {\rm curl}_\eta  \rt]   \pl_t \oa \rt\|_{L^2} \notag \\
  \les & \sum_ {\substack{0\le i\le 1, \  0\le j\le |\al|,  \ 0\le k \le |\ba|
  }}  \lt\|\sa^{\frac{\iota+|\al|+1}{2}}  \mathcal{I}^{1+i,j,k}|\pl_t^{2-i} \pl^{|\al|-j}\bar\pl^{|\ba|-k} \pl \oa| \rt\|_{L^2}
 \notag  \\
\les & \lt\|\sa^{\frac{\iota+|\al|+1}{2}}  \widetilde{\mathcal{I}}^{3,|\al|,|\ba|}  \rt\|_{L^2}
\les (1+t)^{-3}\sqrt{\mathscr{E}(t)} \sum_{0\le j \le  2+ |\al|+|\ba|  }  \sqrt{\mathscr{E}_{j}(t)},  \label{8.6a}
\end{align}
which, together with  \ef{8.5b}, \ef{7.16-1}, $(\ta^2)_\tau \les (1+\tau)^{-1} \ta^2$, $-(\ta^{-2})_t \les (1+t)^{-1} \ta^{-2}$ and \ef{8.6c}, implies that
\begin{align*}
& (1+t)^{3} \lt\|\sa^{\frac{\iota+|\al|+1}{2}} \pl_t  \pl^{\al}\bar\pl^{\ba}{\rm curl}_\eta \pl_t\oa \rt\|_{L^2} \notag \\
 \les & \lt( \lt\|\sa^{\frac{\iota+|\al|+1}{2}}  \pl^{\al}\bar\pl^{\ba}{\rm curl}_\eta \pl_t \oa  \rt\|_{L^2}  +  \lt\|\sa^{\frac{\iota+|\al|+1}{2}}  \pl^{\al}\bar\pl^{\ba} \lt[\pl_t , {\rm curl}_\eta  \rt]   \pl_t \oa \rt\|_{L^2} \rt)(t=0) \notag \\
 &  +  \sup_{\tau \in [0,t]}  \lt\{ (1+\tau)^2\lt\|\sa^{\frac{\iota+|\al|+1}{2}}  \pl^{\al}\bar\pl^{\ba} \lt[\pl_\tau , {\rm curl}_\eta  \rt]   \pl_\tau \oa \rt\|_{L^2} \rt\} \notag\\
  & + \sup_{\tau \in [0,t]}  \lt\{(1+\tau)^3 \lt\|\sa^{\frac{\iota+|\al|+1}{2}} \pl_\tau \pl^{\al}\bar\pl^{\ba} \lt[\pl_\tau , {\rm curl}_\eta  \rt]   \pl_\tau \oa \rt\|_{L^2}\rt\} \\
 \les &\lt\|\sa^{\frac{\iota+|\al|+1}{2}}  \pl^{\al}\bar\pl^{\ba}{\rm curl}_\eta \pl_t \oa  \rt\|_{L^2}(t=0) + \sup_{\tau \in [0,t]} \sqrt{ \mathscr{E}(\tau) }   \sum_{0\le j \le 2+ |\al|+|\ba|  }\sup_{\tau \in [0,t]} \sqrt{ \mathscr{E}_{j}(\tau) } .
\end{align*}

In a similar way to deriving \ef{8.5b}, we have for  $m\ge 3$,
\begin{align*}
 \pl_t^{m-1} \pl^{\al}\bar\pl^{\ba} {\rm curl}_\eta    \pl_t \oa  =  I.D.    + \sum_{0\le i\le m-1}   \frac{(m-1)!}{i!(m-1-i)!}  e^{-t} \lt( \frac{d^{i}}{dt^{i}}  \ta^{-2}(t) \rt) \notag \\
  \times \int_0^t  e^{\tau} \pl_\tau^{m-1-i} \lt( \ta^2(\tau) \pl^{\al}\bar\pl^{\ba} \lt[\pl_\tau , {\rm curl}_\eta  \rt]   \pl_\tau \oa  \rt) d\tau  ,
\end{align*}
where
$$ e^{t} \|\sa^{\frac{\iota+|\al|+1}{2}}  I.D.  \|\les \lt\|\sa^{\frac{\iota+|\al|+1}{2}}  \pl^{\al}\bar\pl^{\ba}{\rm curl}_\eta \pl_t \oa  \rt\|_{L^2}(t=0)  +  \sqrt{\mathscr{E}(0)} \sum_{0\le j \le  m+ |\al|+|\ba|  }  \sqrt{\mathscr{E}_{j}(0)}.$$
This, together with \ef{decay},  \ef{7.16-1}, \ef{5.29}, \ef{5.30-1} and \ef{5.30-2}, proves \ef{7.26c} for  $m\ge 3$.

{\em Step 3}.
In this step, we prove that for $|\al|+|\ba|= [\iota]+7$,
\begin{align}
 (1+t)^2  \lt\|\sa^{\frac{\iota+|\al|+1}{2}} \pl^{\al}\bar\pl^{\ba}{\rm curl}_\eta\pl_t  \oa \rt\|_{L^2}^2
 \les
 \lt\|\sa^{\frac{\iota+|\al|+1}{2}} \pl^{\al}\bar\pl^{\ba}{\rm curl}_\eta\pl_t  \oa \rt\|_{L^2}^2 (t=0) \notag \\
   +    \sup_{\tau \in [0,t]}  \mathscr{E}(\tau)   \sum_{0\le j \le |\al|+|\ba|  }\sup_{\tau \in [0,t]}  \mathscr{E}_{j}(\tau). \label{8.3c}
\end{align}
It follows from \ef{8.4c} and
  \ef{curl-12} that
\begin{align}
&  \lt\|\sa^{\frac{\iota+|\al|+1}{2}} \pl^{\al}\bar\pl^{\ba}{\rm curl}_\eta\pl_t  \oa \rt\|_{L^2}
\les  \lt\|\sa^{\frac{\iota+|\al|+1}{2}} \pl^{\al}\bar\pl^{\ba}{\rm curl}_\eta\pl_t  \oa \big|_{t=0} \rt\|_{L^2}^2 e^{-t} \ta^{-2}(t) \notag\\
&\qquad   + \sum_{1\le l \le 3}  \lt\|\sa^{\frac{\iota+|\al|+1}{2}} e^{-t} \ta^{-2}(t) \int_0^t  e^{\tau} \ta^2(\tau)  \mathcal{Z}^{\al, \ba}_{1,l}    d\tau   \rt\|_{L^2} \notag \\
&\qquad    +  e^{-t} \ta^{-2}(t) \sum_{1\le l \le 3}    \int_0^t  e^{\tau} \ta^2(\tau)   \lt\| \sa^{\frac{\iota+|\al|+1}{2}}\mathcal{Z}^{\al, \ba}_{2,l}  \rt\|_{L^2}  d\tau, \label{8.3d}
\end{align}
which proves \ef{8.3c} by use of the following estimates:
\begin{align}
&e^{-t} \ta^{-2}(t) \sum_{1\le l \le 3}    \int_0^t  e^{\tau} \ta^2(\tau)   \lt\| \sa^{\frac{\iota+|\al|+1}{2}}\mathcal{Z}^{\al, \ba}_{2,l}  \rt\|_{L^2}  d\tau
 \notag\\
\les & (1+t)^{-2} \sup_{\tau \in [0,t]} \sqrt{\mathscr{E}(\tau)}  \sum_{0\le j \le |\al|+|\ba|  }\sup_{\tau \in [0,t]} \sqrt{\mathscr{E}_{j}(\tau)}, \label{8.15a}
\end{align}
and
\begin{align}
& \sum_{1\le l \le 3}  \lt\|\sa^{\frac{\iota+|\al|+1}{2}} e^{-t} \ta^{-2}(t) \int_0^t  e^{\tau} \ta^2(\tau)  \mathcal{Z}^{\al, \ba}_{1,l}    d\tau   \rt\|_{L^2}
\notag\\
\les &(1+t)^{-1} \sup_{\tau \in [0,t]} \sqrt{\mathscr{E}(\tau)}  \sum_{0\le j \le |\al|+|\ba|  }\sup_{\tau \in [0,t]} \sqrt{\mathscr{E}_{j}(\tau)}.\label{8.15b}
\end{align}
Indeed, \ef{8.15a} follows from \ef{8.4a} and \ef{7.16-1}, and \ef{8.15b} follows from \ef{8.4b} and \ef{8.5a}. For example, in the case of $(|h|,|g|)\in S_2$ and $|g|=0$, we have that
\begin{align*}
&   \lt\|\sa^{\frac{\iota+|\al|+1}{2}}  e^{-t} \ta^{-2}(t) \int_0^t e^{\tau} \ta^2(\tau) \epsilon^{ljk} (\pl_\tau \pl^h A^r_j) \pl_\tau  \pl^{ \al-h  }\bar\pl^{ \ba  }  \pl_r  \oa_k      d\tau   \rt\|_{L^2}\\
 \les &
 \lt\|   \sa^{ \frac{\iota+|\al|+1}{2}} \mathcal{I}^{1,|h|,0}\pl^{|\al|-|h| }\bar\pl^{|\ba|}\pl \oa    \rt\|_{L^2} (\tau =t) \\
& +  e^{-t} \ta^{-2}(t)  \lt\|   \sa^{ \frac{\iota+|\al|+1}{2}} \mathcal{I}^{1,|h|,0}\pl^{|\al|-|h| }\bar\pl^{|\ba|}\pl \oa    \rt\|_{L^2}(\tau=0) \\
& + (1+t)^{-2} \sup_{\tau\in [0,t]} \lt\{  (1+\tau)^2 \lt\|   \sa^{ \frac{\iota+|\al|+1}{2}} \mathcal{I}^{2,|h|,0}\pl^{|\al|-|h| }\bar\pl^{|\ba|}\pl \oa \rt\|_{L^2} \rt\} \\
 &   +(1+t)^{-1} \sup_{\tau\in [0,t]}  \lt\{  (1+\tau) \lt\|   \sa^{ \frac{\iota+|\al|+1}{2}} \mathcal{I}^{1,|h|,0}\pl^{|\al|-|h| }\bar\pl^{|\ba|}\pl \oa    \rt\|_{L^2} \rt\}\\
 \les & (1+t)^{-1} \sup_{\tau \in [0,t]} \sqrt{\mathscr{E}(\tau)}  \sum_{0\le j \le |\al|+|\ba|  }\sup_{\tau \in [0,t]} \sqrt{\mathscr{E}_{j}(\tau)}.
\end{align*}
The other cases of $(|h|,|g|)\in S_2\cup S_3$ can be done analogously.

{\em Step 4}.
Based on the estimates obtained in {\em Step 1-3}, we can prove \ef{curl} and \ef{curl-c} by use of the following commutator estimates.
In view of \ef{commutator2}, \ef{5.29}, \ef{5.30-1} and \ef{6.7-1c}, we see that for  $m \ge 1$,
\begin{align*}
& \lt| {\rm curl}_\eta   \pl_t^{m } \pl^{\al}\bar\pl^{\ba}\oa -\pl_t^{m-1} \pl^{\al}\bar\pl^{\ba}{\rm curl}_\eta \pl_t\oa \rt|   \\
\les & |\pl_t^m \pl^{ \al }[\pl, \bar\pl^{ \ba }]\oa|
+\sum_{\substack{0\le i\le m-1, \ 0\le j\le |\al|, \ 0\le k \le |\ba|
\\ 1\le i+j+k
  }}  \mathcal{I}^{i,j,k}
\lt|\pl_t^{m-i}\pl^{|\al|-j}\bar\pl^{|\ba|-k}\pl\oa\rt| \\
\les & \sum_{0\le k\le |\ba|-1} \lt| \pl_t^m \pl^{|\al|+1}\bar\pl^{k}\oa   \rt|
+   \widetilde{\mathcal{I}}^{m,|\al|,|\ba|}
\les   \sum_{0\le k\le |\ba|-1} \lt|\nabla_\eta \pl_t^m \pl^{|\al|}\bar\pl^{k}\oa   \rt|
+   \widetilde{\mathcal{I}}^{m,|\al|,|\ba|},
\end{align*}
which, together with \ef{5.30-2}, implies that for  $m \ge 1$,
 \begin{align}
& (1+t)^{2m}\lt\|\sa^{\frac{\iota+|\al|+1}{2}} \lt( {\rm curl}_\eta   \pl_t^{m } \pl^{\al}\bar\pl^{\ba}\oa -\pl_t^{m-1} \pl^{\al}\bar\pl^{\ba}{\rm curl}_\eta \pl_t\oa  \rt) \rt\|_{L^2}^2 \notag \\
& \qquad
\les \sum_{0\le k\le |\ba|-1} \mathfrak{E}_{II}^{m,|\al|,k}(t) +  \mathscr{E}(t) \sum_{0\le j\le m+|\al|+|\ba|-1}\mathscr{E}_j(t) . \label{8.8b}
\end{align}
Similarly, we have
\begin{align*}
& \lt|{\rm curl}_\eta  \pl^{\al}\bar\pl^{\ba} \oa - \pl^{\al}\bar\pl^{\ba} {\rm curl}_\eta  \oa \rt| \\
\les &  | \pl^{ \al }[\pl, \bar\pl^{ \ba }]\oa|
+\sum_{\substack{  0\le j\le |\al|, \ 0\le k \le |\ba|,
\  1\le  j+k
  }}  \mathcal{I}^{0,j,k}
\lt| \pl^{|\al|-j}\bar\pl^{|\ba|-k}\pl\oa\rt| \\
\les &    | \pl^{ \al }[\pl, \bar\pl^{ \ba }]\oa|
+|\pl\oa|\lt|   \pl^{|\al|}\bar\pl^{|\ba|}\pl \oa   \rt| +   \widetilde{\mathcal{I}}^{0,|\al|,|\ba|}\\
\les &   \sum_{0\le k\le |\ba|-1} \lt| \nabla_\eta \pl^{|\al| }\bar\pl^{k}\oa   \rt|
+|\pl\oa| \lt|   \pl^{|\al|+1}\bar\pl^{|\ba|}  \oa   \rt| +   \widetilde{\mathcal{I}}^{0,|\al|,|\ba|},
\end{align*}
so that
\begin{align}
\lt\|\sa^{\frac{\iota+|\al|+1}{2}} \lt( {\rm curl}_\eta  \pl^{\al}\bar\pl^{\ba} \oa - \pl^{\al}\bar\pl^{\ba} {\rm curl}_\eta  \oa  \rt)\rt\|_{L^2}^2 \notag  \\
\les
 \sum_{0\le k\le |\ba|-1} \mathfrak{E}_{II}^{0,|\al|,k}(t) + \mathscr{E}(t)   \sum_{0\le j \le |\al|+|\ba|  }  \mathscr{E}_{j}(t) . \label{8.8c}
\end{align}
So, \ef{curl-a} can be derived from \ef{7.26a} and \ef{8.8c}; \ef{curl-b} from \ef{7.26c} and \ef{8.8b}; and \ef{curl-c} from \ef{8.3c} and \ef{8.8b}.
\hfill $\Box$

\section{Proof of Theorem \ref{thm3.1}}
The proof is based on the estimates obtained in Propositions \ref{newedv} and \ref{prop-curl}.
It follows from  \ef{curl} and \ef{7.10}  that for $k=0,1,\cdots, [\iota]+7$,
\begin{align}
 \mathfrak{V}_k(t) \les  &  \sum_{0\le m\le 1, \ n+l=k}  \lt\|\sa^{\frac{\iota+n+1}{2}}  \pl^{n}\bar\pl^{l}{\rm curl}_\eta \pl_t^m \oa \big|_{t=0} \rt\|_{L^2}^2  + \sum_{0\le j\le k-1} \mathfrak{E}_{j}(t)  \notag \\
 & +  (1+t)^{-2} \sum_{0\le  n+l \le  k-1 }  \lt\|\sa^{\frac{\iota+n+1}{2}}  \pl^{n}\bar\pl^{l}{\rm curl}_\eta \pl_t  \oa \big|_{t=0} \rt\|_{L^2}^2 \notag\\
 &   +   \sup_{s\in [0,t]}  \mathscr{E}(s)   \sum_{0\le j \le k } \lt(\sup_{s\in [0,t]} \mathfrak{E}_{j}(s) + \ln (1+t)   \int_0^t     \mathfrak{D}_{j}(s) ds \rt), \notag
\end{align}
which implies that for $k=0,1,\cdots, [\iota]+7$,
\begin{align}
& \int_0^t (1+s)^{-1} \mathfrak{V}_k(s) ds \les \ln(1+t) \sum_{0\le m\le 1, \ n+l=k}  \lt\|\sa^{\frac{\iota+n+1}{2}}  \pl^{n}\bar\pl^{l}{\rm curl}_\eta \pl_t^m \oa \big|_{t=0} \rt\|_{L^2}^2  \notag\\
&\qquad + \sum_{0\le j\le k-1} \int_0^t  \mathfrak{D}_{j}(s) ds  +   \sum_{0\le  n+l \le  k-1 }  \lt\|\sa^{\frac{\iota+n+1}{2}}  \pl^{n}\bar\pl^{l}{\rm curl}_\eta \pl_t  \oa \big|_{t=0} \rt\|_{L^2}^2
 \notag \\
&\qquad + \ln(1+t)  \sup_{s\in [0,t]}  \mathscr{E}(s)   \sum_{0\le j \le k } \lt(\sup_{s\in [0,t]} \mathfrak{E}_{j}(s) + \ln (1+t)   \int_0^t     \mathfrak{D}_{j}(s) ds \rt) . \notag
\end{align}
These, together with \ef{initial}, give that for $j=0,1,\cdots, [\iota]+7$,
\begin{align}
&\sum_{0\le k \le j }\lt(    \mathfrak{V}_k  (t)  + \int_0^t(1+s)^{-1} \mathfrak{V}_k (s)ds \rt)\notag\\
\les & \lt(  \ln(1+t) +1  \rt) \sum_{0\le m\le 1, \ 0\le  n+l\le j}  \lt\|\sa^{\frac{\iota+n+1}{2}}  \pl^{n}\bar\pl^{l}{\rm curl}_\eta \pl_t^m \oa \big|_{t=0} \rt\|_{L^2}^2\notag \\
& + \sum_{0\le k\le j-1} \lt( \sup_{s\in [0,t]} \mathfrak{E}_{k}(s)
+\int_0^t  \mathfrak{D}_{k}(s) ds  \rt)
+ \ea_0^2  \lt(\sup_{s\in [0,t]} \mathfrak{E}_{j}(s) +   \int_0^t     \mathfrak{D}_{j}(s) ds \rt). \label{9.16-3}
\end{align}
We can use \ef{new-hig}, \ef{9.16-3} and the mathematical induction argument to obtain that  for $j=0,1,\cdots, [\iota]+7$,
\begin{align}
& \mathfrak{E}_{j}(t) + \int_0^t \mathfrak{D}_{j}(s)ds \les   \sum_{0\le k \le j }  \mathfrak{E}_k  (0) \notag \\
& \qquad +   \lt( \ln(1+t) +1 \rt)  \sum_{0\le m\le 1, \ 0\le  n+l \le j}  \lt\|\sa^{\frac{\iota+n+1}{2}}  \pl^{n}\bar\pl^{l}{\rm curl}_\eta \pl_t^m \oa \big|_{t=0} \rt\|_{L^2}^2 , \label{9.17}
\end{align}
which, with the aid of \ef{7.10}, implies that
\begin{align*}
 \mathscr{E}(t) \les  \mathscr{E}(0) + \mathfrak{V}_{add}(0) +  \ln(1+t) \mathfrak{V}_{add}(0) .
 \end{align*}
Moreover, it follows from \ef{7.26a}, \ef{7.26c}, \ef{8.3c}, \ef{7.10}, \ef{initial} and \ef{9.17} that
\begin{align*}
 \mathfrak{V}_{add}(t)\les &  \mathfrak{V}_{add}(0) + \sup_{s\in [0,t]}  \mathscr{E}(s)   \sum_{0\le j \le [\iota]+7 } \lt(\sup_{s\in [0,t]} \mathfrak{E}_{j}(s) + \ln (1+t)   \int_0^t     \mathfrak{D}_{j}(s) ds \rt)\\
 \les & \mathfrak{V}_{add}(0) + \ea_0^2   \sum_{0\le j \le [\iota]+7 } \lt(\sup_{s\in [0,t]} \mathfrak{E}_{j}(s) +   \int_0^t     \mathfrak{D}_{j}(s) ds \rt)\\
 \les  &    \mathscr{E}(0) + \mathfrak{V}_{add}(0) +  \ln(1+t) \mathfrak{V}_{add}(0) .
 \end{align*}
This proves \ef{energy} and finishes the proof of Theorem \ref{thm3.1}.

\section*{Acknowledgements}

I would like to thank Professor Tao Luo and  Professor Chongchun Zeng for their interests  in this work and helpful discussions.
This research was supported in part by NSFC  Grants 11822107 and 11671225.


\bibliographystyle{plain}

\begin{thebibliography}{99}


\bibitem{ba} Barenblatt, G.: On one class of solutions of the one-dimensional problem of non-stationary
filtration of a gas in a porous medium, {\it Prikl. Mat. i. Mekh.} {\bf 17}, 739--742  (1953).

\bibitem{6'}  Chandrasekhar, S.:  {\it Introduction to the Stellar
Structure},  University of Chicago Press (1939).




\bibitem{chemin1} Chemin, J.:  Dynamique des gaz  a
 masse totale finie, {\em Asymptotic Anal.} {\bf 3}, 215--220 (1990).

\bibitem{chemin2}  Chemin, J.:  Remarques sur l¡¯apparition de singularit{\rm  {e}}s dans les {\rm  {e}}coulements eul{\rm {e}}riens
compressibles, {\em Comm. Math. Phys.} {\bf 133}, 323--329 (1990).



\bibitem{7} Coutand, D.,  Lindblad, H.,  Shkoller,  S.:  A priori estimates for the free-boundary
3-D compressible Euler equations in physical vacuum,  {\it Comm. Math. Phys.} {\bf 296},
559--587  (2010).




\bibitem{10}  Coutand,  D.,  Shkoller, S.:  Well-posedness in smooth function spaces for the moving-
boundary 1-D compressible Euler equations in physical vacuum, {\it Comm. Pure
Appl. Math.} {\bf 64}, 328--366   (2011).





\bibitem{10'} Coutand,  D., Shkoller, S.:  Well-Posedness in Smooth Function Spaces for the Moving-Boundary Three-Dimensional Compressible Euler Equations in Physical Vacuum, {\it Arch. Ration. Mech. Anal.} {\bf 206},  515--616 (2012).


\bibitem{cox}   Cox. J.,   Giuli, R.:  {\it Principles of stellar structure, I.,II.},   New York: Gordon and Breach, 1968.

\bibitem{DM} Dacorogna, B., Moser, J.: On a partial differential equation involving the Jacobian determinant, {\em Ann. Inst. H. Poincar\`{e} Anal. Non Lin\`{e} aire} {\bf 7}, 1--26 (1990).





\bibitem{Friedrichs} Friedrichs, K.: Symmetric hyperbolic linear differential equations, {\it Comm. Pure Appl. Math.} {\bf 7} 345--392 (1954).









\bibitem{zhenlei} Gu, X.,  Lei, Z.:  Well-posedness of 1-D compressible Euler-Poisson equations with physical vacuum, {\it J. Diff. Equ.} {\bf 252},   2160--2188 (2012).

\bibitem{zhenlei1} Gu, X.,  Lei, Z.: Local Well-posedness of the three dimensional compressible Euler--Poisson equations with physical vacuum,  {\it J. Math. Pures Appl.} {\bf 105},  662--723 (2016).

\bibitem{HaJa1} Had\u{z}i\`{c}, M.,  Jang, J.:  Expanding large global solutions of the equations of compressible fluid mechanics,  {\em Invent. Math.} {\bf 214}, 1205--1266  (2018).

\bibitem{HaJa2} Had\u{z}i\`{c}, M.,  Jang, J.:  A class of global solutions to the Euler-Poisson system, {\em arXiv:1712.00124}.





\bibitem{HL} Hsiao, L., Liu, T. P. :  Convergence to nonlinear diffusion waves for solutions of a system of hyperbolic conservation laws with damping, {\it Comm. Math. Phys.} {\bf 143},  599--605 (1992).

 \bibitem{HMP}   Huang, F., Marcati, P.,   Pan, R.:  Convergence to the Barenblatt solution for the compressible Euler equations with damping and vacuum,  {\it Arch. Ration. Mech. Anal.} {\bf 176}, 1--24   (2005).

 \bibitem{HPW} Huang, H.,   Pan, R., Wang, Z.:   $L^1$ convergence to the Barenblatt solution for compressible Euler equations with damping, {\it Arch. Ration. Mech. Anal.} {\bf 200},  665--689   (2011).






\bibitem{16} Jang, J., Masmoudi, N.: Well-posedness for compressible Euler with physical vacuum
singularity,  {\it Comm. Pure Appl. Math.} {\bf 62}, 1327--1385 (2009).






\bibitem{16'} Jang, J.,   Masmoudi, N.: Well-posedness of compressible Euler equations in a physical vacuum,  {\it Comm. Pure Appl. Math.}  {\bf 68},   61--111 (2015).


\bibitem{Kato} Kato, T.: The Cauchy problem for quasi-linear symmetric hyperbolic systems, {\it Arch. Rational Mech. Anal.}  {\bf 58}, 181--205 (1975).


\bibitem{17}  Kreiss, H.:   Initial boundary value problems for hyperbolic systems, {\it Comm. Pure
Appl. Math.} {\bf 23}, 277--296 (1970).



\bibitem{18'} Kufner, A.,   Maligranda, L.,  Persson, L. E.: {\it The Hardy inequality}, Vydavatelsky Servis, Plzen, 2007. About its history and some related results.






\bibitem{23} Liu, T.-P.:  Compressible flow with damping and vacuum,  {\it Jpn. J. Appl.Math.} {\bf 13}, 25--32  (1996).

\bibitem{LiuSmoller} Liu, T.-P., Smoller, J.: On the vacuum state for isentropic gas dynamics equations, {\em Adv.
Math.} {\bf 1}, 345--359 (1980).




\bibitem{24} Liu, T.-P.,  Yang, T.:  Compressible Euler equations with vacuum, {\it J. Differ. Equ.} {\bf 140}, 223--237  (1997).

\bibitem{25} Liu, T.-P. ,  Yang, T.:  Compressible flow with vacuum and physical singularity. {\it   Methods
 Appl. Anal.} {\bf 7}, 495--310 (2000).







\bibitem{LXZ}Luo, T.,  Xin , Z.,   Zeng, H.:  Well-Posedness  for the Motion  of  Physical Vacuum of  the Three-dimensional Compressible Euler Equations with or without Self-Gravitation,  {\it Arch. Rational Mech. Anal}  {\bf 213}, 763--831 (2014).




\bibitem{LZ}  Luo, T.,   Zeng, H.:   Global Existence of Smooth Solutions and Convergence to Barenblatt Solutions for the Physical Vacuum Free Boundary Problem of Compressible Euler Equations with Damping, {\it Comm. Pure  Appl. Math.} {\bf 69}, 1354--1396 (2016).


\bibitem{MUK}  Makino, T., Ukai, S.: On the existence of local solutions of the Euler-Poisson equation for the evolution of gaseous stars, {\it J. Math. Kyoto Univ.} {\bf 27}, 387--399  (1987).


\bibitem{Makino} Makino, T., Ukai, S., Kawashima, S.: On the compactly supported solution of the compressible Euler equation, {\it Japan J Appl Math} {\bf 3}, 249--257  (1986).

\bibitem{Oliynyk} Oliynyk, T.: Future stability of the FLRW fluid solutions in the presence of a positive
cosmological constant, {\em Commun. Math. Phys.} {\bf 346}, 293--312 (2016).


\bibitem{PHJ} Parmeshwar, S.,  Had\u{z}i\`{c}, M.,  Jang, J.: Global expanding solutions of compressible Euler equations with small initial densities, {\em arXiv:1904.01122}.


\bibitem{CHJ} Rickard, C., Had\u{z}i\`{c}, M.,  Jang, J.: Global existence of the nonisentropic compressible Euler equations with vacuum boundary surrounding a variable entropy state, {\em arXiv:1907.01065}



\bibitem{Rodnianski} Rodnianski, I., Speck, J.: The nonlinear future stability of the FLRWfamily of solutions to
the irrotational EuleršCEinstein system with a positive cosmological constant, {\em J. Eur. Math.
Soc.} {\bf 15}, 2369--2462 (2013).


\bibitem{serre} Serre, D.: Expansion of a compressible gas in vacuum, {\em Bull. Inst. Math. Acad. Sin. Taiwan} {\bf 10}, 695--716 (2015).

\bibitem{ShSi} Shkoller, S.,  Sideris, T.: Global existence of near-affine solutions to the compressible Euler equations, {\em arXiv:1710.08368}.


\bibitem{sideris1} Sideris, T.: Spreading of the free boundary of an ideal fluid in a vacuum, {\em J. Differ. Equ.} {\bf 257}, 1--14 (2014).

\bibitem{sideris2}Sideris, T.: Global existence and asymptotic behavior of affine motion of 3D ideal fluids surrounded by vacuum, {\em Arch. Ration. Mech. Anal.} {\bf 225}, 141--176 (2017).





















\bibitem{38}Xu, C.,  Yang, T.:   Local existence with physical vacuum boundary condition to Euler
equations with damping, {\it J. Differ. Equ.} {\bf 210}, 217--231  (2005).







\bibitem{39}  Yang, T.: Singular behavior of vacuum states for compressible fluids,
{\it J. Comput. Appl. Math.} {\bf 190}, 211--231 (2006).

\bibitem{HZ} Zeng, H.:  Global resolution of the physical vacuum singularity for three-dimensional isentropic inviscid flows with damping in spherically symmetric motions, {\em Arch. Ration. Mech. Anal.} {\bf 226},  33--82  (2017).





















\end{thebibliography}

\noindent Huihui Zeng\\
Department of Mathematics\\
\& Yau Mathematical Sciences Center\\
Tsinghua University\\
Beijing, 100084, China;\\
Email: hhzeng@mail.tsinghua.edu.cn

\end{document}